\documentclass{amsart}

\usepackage{amsmath}
\usepackage{amssymb, color}
\usepackage[margin=1in]{geometry}
\usepackage{graphicx}
\usepackage{subfigure}
\usepackage{float}
\usepackage{hyperref}
\usepackage{epsf}
\usepackage{titletoc}
\usepackage{ulem}
\usepackage{algorithm, algorithmicx}
\usepackage{algpseudocode}
\usepackage{multirow}
\usepackage{pifont}

\newtheorem{remark}{Remark}[section]
\newtheorem{lemma}{Lemma}[section]
\newtheorem{theorem}{Theorem}[section]

\newtheorem{definition}{Definition}[section]

\graphicspath{{figure/}}

\newcommand{\bd}{\boldsymbol}

\renewcommand{\hat}{\widehat}
\newcommand{\vct}[1]{\bd{#1}}

\newcommand{\CalD}{{\mathcal{D}}}

\newcommand{\CalP}{{\mathcal{P}}}
\newcommand{\CalT}{{\mathcal{T}}}
\newcommand{\EE}{\mathbb{E}}
\newcommand{\RR}{\mathbb{R}}

\newcommand{\CC}{\mathbb{C}}
\newcommand{\Or}{\mathcal{O}}

\newcommand{\Cov}{\text{Cov}}
\newcommand{\Cost}{\text{Cost}}

\newcommand{\rd}{\mathrm{d}}

\newcommand{\RN}[1]{%
  \textup{\uppercase\expandafter{\romannumeral#1}}%
}

\title{Exploring the locally low dimensional structure in solving random elliptic PDEs}

\author{Thomas Y. Hou}
\thanks{Applied and Comput. Math, Caltech, Pasadena, CA 91125. {\it Email: hou@cms.caltech.edu.}}
\author{Qin Li}
\thanks{Math, UW-Madison, Madison, WI 53705. {\it Email: qinli@math.wisc.edu.}}
 \author{Pengchuan Zhang}
\thanks{Applied and Comput. Math, Caltech, Pasadena, CA 91125. {\it Email: pzzhang@cms.caltech.edu.}}

\date{\today}

\begin{document}
\maketitle
\begin{abstract}
We propose a stochastic multiscale finite element method (StoMsFEM) to solve random elliptic partial differential equations with a high stochastic dimension. The key idea is to simultaneously upscale the stochastic solutions in the physical space for all random samples and explore the low stochastic dimensions of the stochastic solution within each local patch. We propose two effective methods to achieve this simultaneous local upscaling. The first method is a high order interpolation method in the stochastic space that explores the high regularity of the local upscaled quantities with respect to the random variables. The second method is a reduced-order method that explores the low rank property of the multiscale basis functions within each coarse grid patch. Our complexity analysis shows that compared with the standard FEM on a fine grid, the StoMsFEM can achieve computational saving in the order of $(H/h)^{d}/(\log(H/h))^k$, where $H/h$ is the ratio between the coarse and the fine gird sizes, $d$ is the physical dimension and $k$ is the local stochastic dimension. Several numerical examples are presented to demonstrate the accuracy and effectiveness of the proposed methods. In the high contrast example, we observe a factor of 2000 speed-up. 
\end{abstract}

\section{Introduction}
Many problems arising from various physical and engineering applications have multiple scale features and uncertainties. For example, to simulate flow in heterogeneous porous media, the permeability field is often characterized as a multiple scale random medium. The parametrization of a multiscale random medium requires a large number of random variables, leading to a high dimensional random partial differential equation (PDE), which is challenging to solve numerically. Similarly, in shallow water modeling, the basin topography can contain multiple scales and high dimensional uncertainties. Moreover, these problems are typically solved for many source terms and boundary conditions. These problems can be formulated using an input-output relation as it is typically done in reduced-order modeling. In the case of flow in porous media, the input space consists of the random permeability field, source terms and/or boundary conditions. The output space depends on the quantities of interest and may consist of the mean of coarse-grid solutions or some other statistical quantities with respect to the solution. In many applications, the dimension of the output space is typically smaller than that of the input space. The main objective of this paper is to design an efficient reduced-order method that takes advantage of the effective low dimensional solution space for problems with multiple scales and large uncertainties.

The direct simulation consists of two steps. First of all, we generate a large number of samples of the random coefficient and numerically solve the corresponding deterministic PDE's. Secondly, we apply an appropriate stochastic method (e.g. Monte Carlo, Stochastic Collocation, etc) to compute the statistical quantities of interest. Because of the presence of small scales in the physical space and high dimensional uncertainties in the stochastic space, the direct simulations of these problems are prohibitively expensive. We need to develop an efficient model reduction method by obtaining a low dimensional parametrization of the solution in both the physical space and the stochastic space. In this paper, when we refer to ``stochastic space'' we mean the space of the parametrized random variables, and ``stochastic dimension" means the number of the parametrized random variables.

There are a number of multiscale methods that use multiscale basis to represent the multiscale solutions in the physical space; see, e.g., \cite{hou_multiscale_1997, hou_PetrovGalerkin_2004, strouboulis2001generalized, owhadi_metric-based_2007, babuska_optimal_2011, owhadi_localized_2011, efendiev_generalized_2013, owhadi_polyharmonic_2014, maalqvist2014localization, ci_multiscale_2014, owhadi2015multi}. Naive applicatin of these multiscale methods to each sample of multiscale random PDE provides little computational saving because a low dimensional representation needs to be recomputed for every sample. There have been a lot of research activities that explore the low dimensional representation of the solutions of the corresponding random PDEs in the stochastic space. In particular, the generalized polynomial chaos methods (gPC)~\cite{ghanem2003stochastic, xiu_wieneraskey_2002, babuska_galerkin_2004, frauenfelder2005finite, babuska_stochastic_2007, nobile_sparse_2008, nobile_anisotropic_2008, doi:10.1137/040615201} have received a lot of attention in the last decade. These methods are very effective when the stochastic dimension is small. However, their performance deteriorates dramatically when the stochastic dimension increases due to {\it the curse of dimensionality}~\cite{cohen_convergence_2010,cohen_analytic_2011}.

It is important to point out that for the problems with high stochastic dimensions their solutions typically have multiple scales in the spatial domain. For example, a random permeability field with a short correlation length has a high stochastic dimension, and at the same time it has multiple spatial scales ranging from the size of the physical domain to the correlation length of the random permeability field. If we use a traditional method to solve these determinstic problems, we need to use a fine grid mesh that is finer than the correlation length to obtain accurate numerical solutions. Many existing stochastic methods that are used to solve high dimensional stochastic PDEs use standard finite element methods with linear nodal basis, see e.g. ~\cite{nobile_sparse_2008, nobile_anisotropic_2008, bieri2009sparse_FEM, gittelson2013adaptive, doostan2009least, sraj2015coordinates, doostan2011non, yan2012stochastic, cohen_convergence_2010,cohen_analytic_2011}.  The computational cost of these methods could be very expensive for every sample. The total computational cost can be tremendous since we need to simulate many sample solutions.

In this paper, we propose a reduced-order method that performs model reduction in the physical space for all samples simultaneously by using a local parametrization of the random coefficients. Our method can significantly speed up the existing {\it non-intrusive} stochastic methods. By ``non-intrusive stochastic methods'', we mean those methods that can call a deterministic PDE solver as a blackbox, e.g., Monte Carlo, multilevel Monte Carlo~\cite{giles2008multilevel, barth2011multi, cliffe2011multilevel}, (sparse grid) stochastic collocation~\cite{babuska_stochastic_2007, nobile_sparse_2008, nobile_anisotropic_2008, doi:10.1137/040615201}, least-squares methods~\cite{doostan2009least, sraj2015coordinates} and compressed sensing methods~\cite{doostan2011non, yan2012stochastic}. Our method is based on the following observation: most deterministic model reduction methods only require solving local problems, e.g. \cite{hou_multiscale_1997, hou_PetrovGalerkin_2004, owhadi_polyharmonic_2014, maalqvist2014localization, owhadi2015multi}, and the local problems often have much lower stochastic dimensions. To be more specific, the random coefficients restricted to a local subdomain can be parametrized by a much smaller number of parameters, which depends only on the ratio between the subdomain size and the correlation length of the random coefficients. Therefore, the local upscaling (equivalent to deterministic model reduction in our paper) results in low stochastic dimensional problems locally in the physical space and can be efficiently precomputed by the gPC like methods in the offline stage. Based on this observation, we propose a stochastic multiscale finite element method (StoMsFEM) to solve the random PDEs that have high stochastic dimension globally but low stochastic dimension locally. This method inherits almost all the advantages of the deterministic model reduction methods, but removes the limitation that the model reduction process needs to be recomputed for every sample. In this paper, we use the following elliptic equation with heterogeneous random coefficients as an example to illustrate the main idea of our approach:
\begin{equation} \label{eqn}
\begin{cases}
-\nabla_x\cdot\left( \kappa (x,\omega)\nabla_x u(x,\omega) \right) = b(x),\quad x\in D, \omega \in \Omega, \\
u(x,\omega) = 0,\quad x \in\partial D
\end{cases} \quad \text{a.s. } P.
\end{equation}
Here, $D \in \RR^d$ is a bounded spatial domain, and $(\Omega, \mathcal{F}, P)$ is a probability space. The random coefficient $\kappa(x,\omega)$ is of high stochastic dimension and has multiscale features. We assume that $\kappa(x,\omega)$ is a symmetric, positive definite matrix satisfying $\lambda_{min} \ge \alpha > 0$, for a.e., $x\in D$ and $\omega \in \Omega$, where $\lambda_{min}$ is the smallest eigenvalue of $\kappa(x,\omega)$. For such coefficients, the solutions are only H\H{o}lder continuous. If $\kappa(x,\omega)$ has multiple scales, the solution will have multiscale features as well. For simplicity, we assume that the forcing function $b(x)$ is deterministic.

Our StoMsFEM method consists of three steps. The first two steps are in the offline stage and the third step is in the online stage. In the first step, we parametrize the random coefficient $\kappa(x,\omega)$ by exploring the locally low dimensional property of the random media. This can be done by several approaches, including the local KL expansion of the random coefficient, sparse PCA~\cite{Zou_PCA_06, dAspremont_sparsePCA, vu2013fantope, ozolins_compressed_2013, lai2014density} and the intrinsic sparse mode decomposition~\cite{hou_LocalModes_2014}. 
In the second step, we apply a deterministic local upscaling method to obtain a {\it parametric upscaled system}. We provide two methods to do the parametric upscaling: random interpolation method and reduced basis method. The random interpolation method takes advantage of the fact that the local upscaled coefficients are analytic functions of the local stochastic parameters, and builds an interpolation scheme for each upscaled coefficient at the coarse-grid level. The random interpolation method can be viewed as a local reduced-order method in the stochastic space. The reduced basis method makes use of the low rank property of the solutions for the local upscaling problems, and prepares a small set of spatial basis functions for each local upscaling problem. The reduced basis method can be viewed as a local reduced-order method in the physical space. In the online stage, for each sample of the random parameters, we either interpolate the upscaled coefficients in the random interpolation setting, or solve the small reduced-order systems to obtain the upscaled coefficients. A numerical coarse-grid solution for this sample can be obtained by solving the upscaled system. 


We have performed a careful computational complexity analysis of our method. The computational cost of the StoMsFEM consists of the offline and online costs. The offline cost is equivalent to solving the random PDE for $N_{\text{off}}$ samples on the fine grid. In the online stage, the computational cost consists of solving the upscaled system $N_{\text{on}}$ times. Our complexity analysis shows that $N_{\text{off}} \ll N_{\text{on}}$ and the offline computational cost of the StoMsFEM is negligible compared with the online cost. Moreover, we show that the ratio between the online cost for the StoMsFEM and the cost of the standard FEM on the fine grid is of the order $(h/H)^{d}(\log(H/h))^k$. Here $H/h$ is the ratio between the coarse and the fine gird sizes, $d$ is the physical dimension and $k$ is the local stochastic dimension. Therefore, the StoMsFEM gives a speed-up of order $(H/h)^d/(\log(H/h))^k$ over the standard FEM method on a fine grid for a single query problem. We have applied the StoMsFEM to solve several random elliptic PDEs with varying degrees of difficulty to demonstrate the accuracy and computational saving of the StoMsFEM. In the high contrast example, we observe a factor of 2000 speed-up over the naive application of the MsFEM.

We remark that the MsFEM achieves computational saving only for multiple queries. For a multi-query problem, the StoMsFEM can reuse the parametric upscaled system that we obtain in the offline stage and thus there is no offline cost for additional source term or boundary condition. The computational saving for the StoMsFEM is even more significant in a multi-query setting.

There are several other methods that share a philosophy similar to that of StoMsFEM. In~\cite{efendiev_generalized_2013}, GMsFEM has been applied to solve parametric PDEs with multiple scales. GMsFEM assumes the coefficients are already parametrized, while in StoMsFEM we need to first build a locally low dimensional parametrization of the random coefficients. In addition, the StoMsFEM can be implemented with any locally upscaling method including GMsFEM. 
In~\cite{doi:10.1137/140970100}, the authors also observed that the random inputs have low stochastic dimensions locally and used the local KL expansion to parametrize the random coefficients. They proposed to combine the deterministic domain decomposition method (DDM) with the local gPC expansions and Monte Carlo sampling to achieve computational saving. A major difference between the StoMsFEM and the method proposed in ~\cite{doi:10.1137/140970100} is that their method does not deal with the multiscale feature, which contributes to the stochastic high dimensionaility. As a result, they have not explored the low-dimensional structure in the physical space. We will compare the StoMsFEM with these methods as we get into details of StoMsFEM. 

The rest of the paper is organized as follows. In Section 2, we give a brief review of several locally low dimensional parametrization methods and the MsFEM. In Section 3, we present our StoMsFEM method that uses either the random interpolation method or the reduced basis method. We also perform complexity analysis for our method. In Section 4, we show how to combine our methodology with existing stochastic methods, e.g., the Monte Carlo (MC), Multi-Level Monte Carlo (MLMC) and sparse grid stochastic collocation (SC) method. In Section 5, we demonstrate the accuracy and efficiency of our method through several numerical examples. Finally, some concluding remarks are given in Section 6.

\section{Preliminaries}
The StoMsFEM we propose relies on two building blocks: locally low dimensional parametrization of a random field and deterministic local upscaling methods. Although these methods are not the focus of the current paper, we give a brief review below for completeness.

\subsection{Locally low dimensional parametrization}\label{sec:parametrization}
Typically for a random coefficient $\kappa(x,\omega)$ in Eqn.~\eqref{eqn}, its parametrization is not known, but rather, its mean and covariance are given:
\begin{equation}\label{eqn:covariance}
\bar{\kappa}(x) = \EE[\kappa(x,\omega)], \quad \Cov{(x,y)} = \EE{\left[(\kappa(x,\omega)-\bar{\kappa}(x))(\kappa(y,\omega)-\bar{\kappa}(y))\right]}.
\end{equation}
In order to solve the random PDE, one first needs to parametrize the random coefficient with some random parameters first. The KL expansion~\cite{karhunen_uber_1947,loeve_probability_1977} is the most popular method in parametrizing the random media. However, the eigenfunctions of the covariance function (also called the KL modes) are global in nature. As a result, the local stochastic dimension is the same as the global stochastic dimension. In this section, we will briefly review a few methods to get a locally low dimensional parametrization, including the local KL expansion, the Intrinsic Sparse Mode Decomposition and the sparse PCA approach. For a detailed description and comparison of these methods, please refer to a companion paper~\cite{hou_LocalModes_2014}.

\subsubsection{Local KL expansion}
The local KL expansion is a natural way to construct a locally low-dimensional parametrization of the random medium (also used in \cite{doi:10.1137/140970100}). Let $D$ be divided into a set of non-overlapping subdomains $\{P_m\}_{m=1}^M$, called patches,
\begin{equation} \label{eqn:domainpar}
	\overline{P} = \cup_{m=1}^M \overline{P_m}, \qquad P_m\cap P_n = \emptyset \text{ for } m \neq n. 
\end{equation}
Let $\Cov_m: P_m \times P_m \rightarrow \RR$~ be the global covariance function $\Cov(x,y)$ restricted to the $m$-th patch:
\begin{equation} \label{eqn:localCov}
\Cov_m(x, y) = \Cov(x, y), \quad    x,y\in P_m.
\end{equation} 
Similar to the standard KL expansion, we can define a local KL expansion as follows:
\begin{definition} [Local KL expansion of $\kappa(x,\omega)$]
Perform KL expansion in each subdomain $P_m$:
\begin{equation} \label{eqn:localKL}
\int_{P_m} \Cov_m(x,y)f_{k,m}(y)\rd{y} = \lambda_{k,m}f_{k,m}(x)\,,\quad \xi_{k,m}(\omega) = \frac{1}{\sqrt{\lambda_{k,m}}} \int_{P_m} \left(\kappa(x,\omega)-\bar{\kappa}(x)\right)f_{k,m}(x)\rd{x}.
\end{equation}
Arrange $\lambda_{k,m}$ in a descending order, and truncate the expansion at the $K_{m}$-th mode. Then, we obtain a local parametrization as follows:
\begin{equation}\label{eqn:local_KL_truncation}
\kappa(x,\omega) \approx \EE[\kappa(x,\cdot)] + \sum_{k=1}^{K_m}\sqrt{\lambda_{k,m}}\xi_{k,m} f_{k,m}(x), \quad x \in P_m.
\end{equation}
\end{definition}

\subsubsection{Intrinsic sparse mode decomposition} \label{sec:localrepresentation}
In~\cite{hou_LocalModes_2014} the authors proposed the intrinsic sparse mode decomposition (ISMD) that decomposes a symmetric positive semidefinite matrix into several sparse rank-one components. We assume that the coviance matrix, $\Cov$, can be decomposed into a finite number of sparse modes, i.e. $\Cov = \sum_{k=1}^K g_k g_k^T$. ISMD looks for a patch-wise sparse decomposition by minimizing the total local dimension, i.e.,
\begin{equation}\label{opt:minsparseness}
	\boxed{\min_{g_1, \dots, g_K}\quad \sum_{m=1}^{M} d_m\,, \quad \text{subject to} \quad \Cov = \sum_{k=1}^K g_k g_k^T\, ,}
\end{equation}
where $d_m$ is the number of nontrivial modes among $\{g_k\}_{k=1}^K$ on the local patch $P_m$, defined as
\begin{equation*}
	d_m = \#\{k ~:~ g_k|_{_{P_m}} \neq \vct{0}\}\,.
\end{equation*}
Under certain non-degenerate assumptions on the covariance $\Cov$ and the partition $\CalP$, we proved that the ISMD exactly produces one minimizer of the minimization problem~\eqref{opt:minsparseness}, see Theorem 3.5 in~\cite{hou_LocalModes_2014}. After projecting the random field $\kappa(x,\omega)$ onto the sparse modes $\{g_k\}_{k=1}^K$, we get a parametrization with $K$ random parameters:
\begin{equation} \label{eqn:sparseKL}
	\kappa(x,\omega) = \bar{\kappa}(x) + \sum_{k=1}^K g_k(x) \xi_k(\omega) ,
\end{equation} 
where the random variables $\{\xi_k\}_{k=1}^K$ are normalized (with center zero and variance one) and uncorrelated. Moreover, the parametrization~\eqref{eqn:sparseKL} achieves the minimal total local stochastic dimension, as desired.

It is worth mentioning that there are several other methods that are able to achieve locally low dimensional parametrization, for example, the sparse PCA~\cite{Zou_PCA_06, dAspremont_sparsePCA, vu2013fantope, LaiLuOsher_2014} and the sparse operator compression~\cite{hou_sparseOC1_2016, hou_sparseOC2_2016}.

When applying the proposed StoMsFEM, the most important factor in choosing a parametrization method is the global stochastic method, which we will discuss in Section~\ref{sec:globalstochastic}. If one wants to use the MC type methods, we recommend to use the local KL expansion. This is the typical case because the StoMsFEM is aiming at stochastically high dimensional problems, where the gPC type methods would have difficulties to deal with. If the global stochastic dimension is within the range of the (sparse grid) SC method, and if one wants to use the (sparse grid) SC method to save computational cost, one should choose the ISMD or the sparse PCA that would parametrize the random coefficients more effectively. The ISMD is recommended when the random parameters are required to be uncorrelated and a high accuracy parametrization is desired, e.g., for the synthetic porous media in Section~\ref{sec:patchstudy} and \ref{exam:highcontrast}. The Sparse PCA and many other matrix factorization methods are good at parametrizing random coefficients whose covariance matrix has continuously decaying eigenvalues, e.g., the Gaussian kernel $\exp(-|x-y|^2/l^2)$ in Section~\ref{subsec:guassiandecay}. We give a more detailed comparison between the ISMD and the sparse PCA in our companion paper~\cite{hou_LocalModes_2014}.

\subsubsection{Nonlinear transformations}\label{sec:nonlineartrans}
All the parametrization methods above, including the KL expansion, are affine with respect to the random parameters. In some applications, the use of nonlinear transformations may reduce the number of parameters significantly. 
For example, the following nonlinear transformation has been widely used for parametrize a positive random field $\kappa(x,\omega)$:
\begin{equation}\label{eqn:nonlinear1}
\kappa(x,\omega) = \kappa_{min} + \exp(\beta(x,\omega))\,.
\end{equation}
The expression has strict positive lower bound $\kappa_{min}$, and in practice we an apply affine parametrization to $\beta(x,\omega)$.

If the random field has both lower and upper bound ($\kappa_{min}$ and $\kappa_{max}$ respectively), the following nonlinear transformation are usually used:
\begin{equation}\label{eqn:nonlinear2}
\kappa(x,\omega) = \frac{\kappa_{max} + \kappa_{min}}{2} + \frac{\kappa_{max} - \kappa_{min}}{2}\tanh(\beta(x,\omega))\,.
\end{equation}

\subsection{Multiscale finite element method}\label{sec:basis}
Model reduction methods based on local upscaling is the other building block in our StoMsFEM. There have been a number of such local upscaling methods for elliptic equations with heterogeneous diffusion coefficients; see, e.g., \cite{hou_multiscale_1997, strouboulis2001generalized, hou_PetrovGalerkin_2004, owhadi_metric-based_2007, babuska_optimal_2011, owhadi_localized_2011, efendiev_generalized_2013, owhadi_polyharmonic_2014, maalqvist2014localization, ci_multiscale_2014, owhadi2015multi}. In this paper, we will use the multiscale finite element method (MsFEM) developed in~\cite{hou_multiscale_1997, hou_PetrovGalerkin_2004} for the local upscaling. We point out that per user's preference, the MsFEM can be replaced by other local upscaling methods with minor modifications. In the following, we briefly review MsFEM applied on Eqn.~\eqref{eqn} with a specific sample media, denoted as $\kappa(x,\omega)$. Note that here the media is fixed and Eqn.~\eqref{eqn} is deterministic.

Suppose that the physical domain $D$ is partitioned into a finite set of compact triangles or quadrilaterals $\{D_H^m, 1 \le m \le M\}$, which forms a triangulation $\CalT_H$ with mesh size $H$. We assume that the coarse grid mesh size $H$ is much larger than the small scale $\epsilon$ in the rough coefficient $\kappa(x,\omega)$, i.e. $H \gg \varepsilon$. In block $D^m_H$, we compute the following cell problem:
\begin{equation} \label{eqn:MsFEB}
\begin{cases}
-\nabla\cdot(\kappa(x,\omega)\nabla\phi^{ml}(x,\omega))=0\,,\quad x\in D^m_H\, ,\\
\phi^{ml}(x,\omega) = p^{ml}(x)\,,\quad x\in\partial D^m_H\,,
\end{cases}\qquad l = 1,\cdots,L\,.
\end{equation}
Here $L$ is the number of nodes on $D^m_H$ and $p^{ml}$ is defined on the boundary $\partial D^m_H$ playing the role of Dirichlet boundary conditions. In our computation we could choose $p^{ml}$ as linear basis (for triangles) or standard bilinear basis (for quadrilaterals) that takes value $1$ at node $l$ and $0$ for all the other $L-1$ nodes in the patch. In practice, we solve the local cell problem on a fine mesh $\CalT_h$ that resolves the small scales in $\kappa(x,\omega)$. Thereafter, the local upscaled stiffness matrix and the load vector can be computed as:
\begin{equation}\label{eqn:stiffness_local}
\begin{split}
\vct{S}_{ll'}^{m}(\omega) &= \int_{D^m_H}\kappa(x,\omega)\nabla\phi^{m l}(x,\omega)\cdot\nabla\phi^{m l'}(x,\omega)\rd{x}\,, \\
\vct{b}_{l}^{m}(\omega) &= \int_{D^m_H} b(x) \phi^{m l}(x,\omega) \rd{x}\,.
\end{split}
\end{equation}
The standard assembling procedure can be utilized to assemble $\vct{S}$ and $\vct{b}$ by looping over all the coarse grid elements. After solving the upscaled system
\begin{equation}\label{eqn:upscaleSUF}
	\vct{S} \vct{U} = \vct{b}\,,
\end{equation}
we obtain the multiscale solution
\begin{equation} \label{eqn:MsSolution}
	u_H(x,\omega) = \sum_m U_m(\omega)\phi_m(x,\omega)\,.
\end{equation}

When the boundary conditions $p^{m l}$ are linear, the following convergence theorem was proved in~\cite{hou_convergence_1999}:
\begin{theorem}\label{thm:MsFEM}
Let $\kappa^\varepsilon(x) = \kappa(\frac{x}{\varepsilon})$ be a smooth periodic medium and $u^\varepsilon(x)$ the solution to~\eqref{eqn}. Denote $u^\varepsilon_H(x)$ the multi-scale finite element approximation obtained from the space spanned by the multiscale basis with linear boundary conditions. Then we have:
\begin{equation*}
\|u^\varepsilon - u^\varepsilon_H\|_{H^1}\leq C(H+\varepsilon)\|f\|_{L^2} + C\sqrt{\frac{\varepsilon}{H}}\|u_0\|_{H^2}\,,
\end{equation*}
where $u_0$ is the solution to the homogenized equation of Eqn.~\eqref{eqn}.
\end{theorem}

\begin{remark}\label{remark:oversampling}
Theorem~\ref{thm:MsFEM} implies that the MsFEM captures the correct homogenized results for small $\varepsilon$. However, the method may produce a large error if $\varepsilon\sim H$. This is called the resonance error between the coarse grid scale $H$ and media small-scale parameter $\varepsilon$. The oversampling technique proposed in~\cite{hou_multiscale_1997} successfully reduces the resonance error, based on the observation that the boundary layer typically gets damped out quickly within a width of $\varepsilon$. This suggests one compute Eqn.~\eqref{eqn:MsFEB} in a larger domain of order $H+\mathcal{O}(\varepsilon)$ in each dimension, and utilizing the interior information for basis construction. This significantly reduces the resonance error and gives more accurate results, as demonstrated in numerical examples in~\cite{hou_multiscale_1997}.


We remark that the over-sampling technique results in a non-conforming FEM. To reduce the non-conforming FEM error, we use the Petrov-Galerkin MsFEM formulation~\cite{hou_PetrovGalerkin_2004}. In the Petrov-Galerkin MsFEM, the local upscaled quantities are computed from
\begin{equation}\label{eqn:stiffness_local_petrov}
\begin{split}
\vct{S}_{ll'}^{m}(\omega) &= \int_{D^m_H}\kappa(x,\omega)\nabla\phi^{m l}(x,\omega)\cdot\nabla\phi_{t}^{m l'}(x)\rd{x}\,, \\
\vct{b}_{l}^{m} &= \int_{D^m_H} b(x) \phi_{t}^{m l}(x) \rd{x}\,,
\end{split}
\end{equation}
with $\{\phi_{t}^{m l}(x)\}_{l=1}^L$, the test functions, being linear (triangular grid) or bilinear (quadrilateral grid) locally. One added benefit of using the Petrov-Galerkin MsFEM is that the test functions are deterministic, and thus the local load vector $\vct{b}^m$ is independent of random samples.
\end{remark}

It is important to point out that there is some computational overhead in computing cell problems~\eqref{eqn:MsFEB} and assembling the upscaled system~\eqref{eqn:stiffness_local}. However, the upscaled stiffness matrix $\vct{S}$ can be reused for different source terms $b(x)$ and the MsFEM achieves no computational saving except in the multi-query setting. This is true for most model reduction methods since there is an overhead in the offline stage in constructing the reduced order models. In the random setting, with a single source term in Eqn.~\eqref{eqn}, if we naively apply the MsFEM for every sample media, the reduced models need to be recomputed for every sample coefficient, which leads no computational saving. In the multi-query case, since random samples may be different for different queries, we are still not able to gain the full power of the MsFEM.

We introduce our stochastic multiscale finite element method (StoMsFEM) below. In this new method, we prepare the upscaled quantities~\eqref{eqn:stiffness_local} for all samples simultaneously in the offline stage. In the online stage, for every sample we only need to assemble and solve the upscaled system~\eqref{eqn:upscaleSUF}. The computational saving is achieved even in the single-query setting.

\section{Stochastic Multiscale Finite Element Method}\label{sec:algorithm}
The StoMsFEM consists of two stages: the offline stage, in which we construct reduced-order models and prepare upscaled quantities, and the online stage, in which we sample media and compute the upscaled system. We discuss offline preparation stage in this section and leave the online global computation to the next section.

The offline stage also consists of two steps. In the first step, we parametrize the random coefficients $\kappa(x,\omega)$ using a local parameterization method, as reviewed in Section~\ref{sec:parametrization}. More discussion is found in the companion paper~\cite{hou_LocalModes_2014}. In this paper, we assume that the random coefficient $\kappa(x,\omega)$ has already been parametrized, denoted as $\kappa(x,\vct{\xi})$, where $\vct{\xi} \equiv [\xi_1, \xi_2, \cdots, \xi_K]^T$ is the collection of $K$ random variables. Note that we do {\it not} assume the affine structure with respect to $\vct{\xi}$ in the parametrization of $\kappa(x,\vct{\xi})$. Our StoMsFEM works for any parametrization with locally low stochastic dimensions. 

With this parametrization, the cell problem~\eqref{eqn:MsFEB} becomes a parametrized PDE:
\begin{equation} \label{eqn:MsFEB2}
\begin{cases}
-\nabla\cdot(\kappa(x,\vct{\xi}_m)\nabla\phi^{ml}(x,\vct{\xi}_m))=0\,,\quad &x\in D^m_H\,,\\
\phi^{ml}(x,\vct{\xi}_m) = p^{ml}(x)\,,\quad & x\in\partial D^m_H\,,
\end{cases} \qquad\forall \vct{\xi}_m\in\Gamma_m\,, l = 1,\cdots,L\,.
\end{equation}
Here $\vct{\xi}_m \in\Gamma_m\subset \RR^{K_m}$ are the local effective parameters with $\Gamma_m$ being its range, and $K_m$ is the local stochastic dimension on patch $D^m_H$. Following the same assembling procedure~\eqref{eqn:stiffness_local}, the local upscaled quantities $\vct{S}^{m}$ and $\vct{b}^{m}$ are functions of the local parameters, i.e., $\vct{S}^{m}(\vct{\xi}_m)$ and $\vct{b}^{m}(\vct{\xi}_m)$. 

In the second step, we construct reduced-order models for the upscaled quantities $\vct{S}^{m}(\vct{\xi}_m)$ and $\vct{b}^{m}(\vct{\xi}_m)$. Two methods are proposed, i.e., the random interpolation method and the reduced basis method. The random interpolation method makes use of the fact that $\phi^{ml}(x,\vct{\xi}_m)$ (and thus $\vct{S}^{m}(\vct{\xi}_m)$ and $\vct{b}^{m}(\vct{\xi}_m)$) are smooth with respect to $\vct{\xi}_m$, while the reduced basis method makes use of the low dimensional structure of $\phi^{ml}(x,\vct{\xi}_m)$ in the physical space. The random interpolation method works for all parametrized local coefficients $\kappa(x, \vct{\xi}_m)$ that are smooth with respect to the local parameters $\vct{\xi}_m$, while the reduced method basis is only recommended when $\kappa(x, \vct{\xi}_m)$ is affine with respect to $\vct{\xi}_m$.


\subsection{The random interpolation method}\label{sec:randomInter}
In the random interpolation method, we propose to compute the cell problem~\eqref{eqn:MsFEB2} and the local upscaled quantities $\vct{S}^{m}(\vct{\xi}_m)$ and $\vct{b}^{m}(\vct{\xi}_m)$ on multiple {\it deterministic} collocation points and use them to build interpolants of the upscaled quantities in terms of the parameters $\vct{\xi}_m$. These interpolants allow us to efficiently approximate the upscaled quantities for any given sample in the online stage.

We take patch $D_H^m$ as an example to illustrate the idea. If $\kappa(x,\vct{\xi}_{m})$ is smooth with respect to $\vct{\xi}_m$, the local multiscale basis functions $\phi^{ml}(x,\vct{\xi}_{m})$ are smooth with respect to $\vct{\xi}_{m}$, see~\cite{babuska_galerkin_2004, babuska_stochastic_2007, cohen_analytic_2011, hoang2011n}. Assume the range for all these random variables lies in the interval $[-1,1]$ (other bounded ranges can be rescaled accordingly), and in the simplest scenario, we sample $\nu+1$ collocation points along each dimension, and the entire collocation set is the tensor product:
\begin{equation*}
\Gamma_\text{1d} = \{-1 \le \xi_{0}<\xi_1<\cdots<\xi_{\nu} \le 1\}\,,\quad\Gamma_{m,c} = \prod_{k=1}^{K_m} \Gamma_\text{1d}\in\Gamma_m\,.\footnote{To lighten the notations, we choose the same number of collocation points along different directions. However, in general we use different numbers of collocation points along different directions.}
\end{equation*}
The ``c" in the subscript stands for collocation. In total $N_c = (\nu+1)^{K_m}$ collocation points are sampled. If the joint distribution of $\vct{\xi}_{m}$ is known, we can choose the 1d collocation nodes $\Gamma_\text{1d}$ in the same way as stochastic collocation methods~\cite{babuska_stochastic_2007}. For example, if $\vct{\xi}_{m}$ are i.i.d. uniformly distributed on $[-1,1]$, then $\Gamma_\text{1d}$ should be the zeros of the Legendre polynomials. If the joint distribution of $\vct{\xi}_{m}$ is unknown, we can simply use the Chebyshev nodes.

For each collocation point $\vct{\xi}_{m,c}\in\Gamma_{m,c}$, we solve the cell problem \eqref{eqn:MsFEB2} for the local basis functions, denoted as $\phi^{ml}(x,\vct{\xi}_{m,c})$. We then assemble the local stiffness matrix and the local load vector, denoted as $\vct{S}^{m}(\vct{\xi}_{m,c})$ and $\vct{b}^{m}(\vct{\xi}_{m,c})$, according to \eqref{eqn:stiffness_local}. Here $\vct{S}^{m}(\vct{\xi}_{m,c})$ is a $L\times L$ matrix and $\vct{b}^{m}(\vct{\xi}_{m,c})$ is a $L$-dimensional vector. We then construct the interpolants of $\vct{S}^{m}$ and $\vct{b}^m$ in terms of the random variables $\vct{\xi}_m$. Such interpolants are constructed for each element of $\vct{S}^m$ and $\vct{b}^m$. For example for each element in $\vct{S}^{m}$ and $\vct{b}^m$, we construct the Lagrange polynomial approximation, denoted as $\hat{\vct{S}}^{m}$ and $\hat{\vct{b}}^m$.


The interpolants $\hat{S}^m$ and $\hat{b}^m$ will be used to obtain the approximation for every sample in the online stage, as Algorithm~\ref{alg:randominterp_online} depicts. We will discuss how to determine the sample set to solve in the online stage in Section~\ref{sec:globalmethods}.

\begin{remark}
Several remarks are in order:
\begin{itemize}
\item When the local stochastic dimensions are moderate, sparse grids and dimensional adaptive grids~\cite{gerstner2003dimension, spdoc} can be used to reduce the number of collocation points. We use these techniques in our numerical examples.
\item To determine the collocation points $\Gamma_m$, we do not need to know the density of the parameters $\vct{\xi}_m$. We only need to know their ranges. This is very different from the standard gPC method, which requires the joint density function of the random parameters. It is possible that the range of some random parameter is an unbounded domain and in this case we truncate its range to a bounded domain that is large enough to cover the parameter range with high probability. For parameter configurations that fall outside this bounded domain, we directly compute $\phi^{ml}(x,\vct{\xi}_{m})$, $\vct{S}^m(\vct{\xi}_{m})$ and $\vct{b}^{m}(\vct{\xi}_{m})$ from \eqref{eqn:MsFEB2} and \eqref{eqn:stiffness_local}.
\end{itemize}
\end{remark}

Finally, we summarize the offline and online implementation of random interpolation method in Algorithm~\ref{alg:randominterp_offline} and Algorithm~\ref{alg:randominterp_online}.
\begin{algorithm}
\caption{The offline stage of random interpolation method}\label{alg:randominterp_offline}
\begin{algorithmic}[1]
\State Partition the physical domain $D$ into coarse grid blocks $\mathcal{T}_{H}$
\For{each patch $D^m_H$}
\State Determine the set of interpolation nodes $\Gamma_{m,c}$
\State Solve local multiscale basis $\{\phi^{ml}(x,\vct{\xi}_{m,c})\}_{c=1}^{N_c}$ according to \eqref{eqn:MsFEB}
\State Assemble the local upscaled quantities $\vct{S}^m$ and $\vct{b}^m$ according to \eqref{eqn:stiffness_local} or \eqref{eqn:stiffness_local_petrov}
\State Build the Lagrange polynomial interpolants for $\vct{S}^m$ and $\vct{b}^m$
\EndFor
\end{algorithmic}
\end{algorithm}

\begin{algorithm}
\caption{The online stage of random interpolation method for a specific configuration $\vct{\xi}$}\label{alg:randominterp_online}
\begin{algorithmic}[1]
\For{each patch $D^k_H$}
\State Determine the values of local variables $\vct{\xi}_{m}$
\State Interpolate the local upscaled stiffness matrix $\vct{S}^m$ and local upscaled loading vector $\vct{b}^m$
\EndFor
\State Assemble and solve the upscaled system \eqref{eqn:upscaleSUF} to obtain $\hat{u}_H(x,\vct{\xi})$
\end{algorithmic}
\end{algorithm}

\subsubsection{Accuracy of the random interpolation}\label{sec:RandomInter_theory}
In this subsection, we estimate the interpolation error $\hat{\vct{S}}^{m} - \vct{S}^{m}$. If the media $\kappa(x, \vct{\xi}_m)$ smoothly depend on $\vct{\xi}_m$, we can prove that the solution to the cell problem~\eqref{eqn:MsFEB2} is also smooth with respect to $\vct{\xi}_m$. So do $\vct{S}^{m}$ and $\vct{b}^{m}$. The strong regularity of $\vct{S}^{m}(\vct{\xi}_m)$ and $\vct{b}^{m}(\vct{\xi}_m)$, combined with high order approximation method, is the key for the success of the random interpolation method.

It is worth mentioning that the regularity problem has been well-studied in the literature, see~\cite{babuska_galerkin_2004, babuska_stochastic_2007, cohen_analytic_2011, hoang2011n}. In particular, with small modification of Lemma 3.2 in~\cite{babuska_stochastic_2007} we have the following lemma.
\begin{lemma}\label{lem:smoothness_basis}
Let $\phi^{ml}$ be the multiscale basis function in the $m$-th cell problem~\eqref{eqn:MsFEB2}. We use the notation $\vct{\xi}_m = [\xi_{m,1},\xi_{m,2},\cdots,\xi_{m,K_m}]$ for the list of the random parameters effective in patch $D^m_H$. We assume that $\vct{\xi}_m \in \Gamma_m \equiv [-1,1]^{K_m}$. If the local parametrization on patch $D_H^m$, i.e. $\kappa(x, \vct{\xi}_m)$, is infinitely differentiable with respect to $\vct{\xi}_m$ and there exists some $c_k > 0$ such that
\begin{equation}\label{eqn:kappagrowth}
	\left\|\frac{\partial^{n}_{\xi_{m,k}} \kappa(\vct{\xi}_m)}{\kappa(\vct{\xi}_m)} \right\|_{L^{\infty}(D^m_H)} \le c_k^n n! \qquad \forall n \ge 0\,,
\end{equation}
then we have for every $\vct{\xi}_m \in \Gamma_m$ and $n\ge 0$
\begin{equation}\label{eqn:phigrowth}
	\left\| \sqrt{\kappa(\vct{\xi}_m)} \nabla \partial^{n}_{\xi_{m,k}} \phi^{ml}(\vct{\xi}_m) \right\|_{L^\infty(\Gamma_m; L^2(D_H^m))} \le C_{\phi} \sqrt{\kappa_{max}} n! (2 c_k)^n,
\end{equation}
where $C_{\phi}$ only depends on the local domain $D^m_H$ and the deterministic boundary condition $p^{ml}$ for $\phi^{ml}$.
\end{lemma}
The assumption~\eqref{eqn:kappagrowth} holds true for most parametrization methods. For example, for a linear parametrization 
\begin{equation*}
	\kappa(x,\vct{\xi}_m) = \bar{\kappa}(x) + \sum_{k=1}^{K_m} g_{m,k}(x) \xi_{m,k},
\end{equation*}
provided that such expansion guarantees $\kappa(x,\vct{\xi}_m) \ge \kappa_{min}$ for almost every $x \in D$ and $\vct{\xi}_m \in \Gamma_m$, we can take $c_k = \|g_k\|_{L^{\infty}(D)}/\kappa_{min}$. When a linear parametrization is combined with the exponential transformation~\eqref{eqn:nonlinear1}, i.e.,
\begin{equation*}
	\kappa(x,\vct{\xi}_m) = \kappa_{min} + \exp\left(\bar{\beta}(x) + \sum_{k=1}^{K_m} g_{m,k}(x) \xi_{m,k}\right),
\end{equation*}
we can take $c_k = \|g_k\|_{L^{\infty}(D)}$.

The regularity of $\phi^{ml}(\vct{\xi}_m)$ in Lemma~\ref{lem:smoothness_basis} could be extended to that of $\vct{S}^{m}$ in a straightforward manner.
\begin{theorem}\label{thm:smoothness_S}
Under the same assumptions in Lemma~\ref{lem:smoothness_basis}, for any $l, l' \in \{1, 2, \dots, L\}$, $\vct{S}_{ll'}^{m}$ as a function of $\xi_{m,k}$, $\vct{S}_{ll'}^{m}: [-1,1]\to C([-1,1]^{K_m-1})$ admits an analytic extension to the complex domain
\begin{equation}\label{eqn:sigmadef}
	\Sigma([-1,1]; \tau_k) \equiv \{z\in \CC, \text{dist}(z, [-1,1]) \le \tau_k\}
\end{equation}
with $0 < \tau_k < 1/(2 c_k)$.
\end{theorem}

Finally, we have the following theorem that guarantees the accuracy of the random interpolation method. 
\begin{theorem}\label{thm:approx_S}
Under the same assumptions in Lemma~\ref{lem:smoothness_basis}, for any $l, l' \in \{1, 2, \dots, L\}$ there exists positive constants $r_k$, $k=1, \dots, K_m$, and $C_s$, independent of $\nu\equiv (\nu_1, \nu_2, \dots, \nu_{K_m})$, such that 
\begin{equation*}
\| \vct{S}_{ll'}^{m} - \hat{\vct{S}}_{ll'}^{m}  \|_{C(\Gamma_m)} \le C_s \prod_{k=1}^{K_m} \left( \frac{2}{\pi}\log(\nu_k + 1) + 1 \right) \sum_{k=1}^{K_m} \exp(- r_k \nu_k)\,,
\end{equation*}
where $\hat{S}_{ll'}^{m}$ is the Chebyshev interpolation of $\vct{S}_{ll'}^{m}$ with $\nu_k+1$ collocation points in the $\xi_{m,k}$ direction, $r_k = \log\left[\tau_k \left( 1 + \sqrt{1 + 1/\tau_k^2} \right) \right]$ and $\tau_k$ is any positive constant that is strictly smaller than the distance between the real line segment $[-1,1]$ and the nearest singularity in the complex plane, as defined in Theorem~\ref{thm:smoothness_S}.
\end{theorem}
The proof is the same as the proof of Theorem~4.1 in~\cite{babuska_stochastic_2007}. The only difference is that we are considering the interpolation in the continuous function space. Therefore, we have the Lebesgue constant of the Chebyshev interpolation, i.e. $\frac{2}{\pi}\log(\nu_k + 1) + 1$, which appears in our error estimation. The regularity and approximation accuracy of $\vct{b}_m$ can be analyzed similarly. 

\begin{remark}\label{rem:unbounded}
Several remarks are in order:
\begin{itemize}
\item The estimates above are based on simple energy estimate and is far from being sharp. More dedicated analysis has been carried in~\cite{cohen_analytic_2011} for improved results and could be easily carried over. Detailed regularity analysis is not the focus of the paper and we omit it from here.  
\item In the case when the random variables have infinite range (e.g. the Gaussian variable), we can sample the collocation points in a range that is large enough to cover the random variable with high probability. For example, denote $\Gamma_{m,0}$ as the finite domain that is large enough to cover a very large portion of $\Gamma_m$ such that: $P(\vct{\xi}_m \notin\Gamma_{m,0})<\varepsilon$, then we can build an interpolant of $\vct{S}^m$ such that $\vct{S}^{m} - \hat{\vct{S}}^{m}$ is small point-wisely on $\Gamma_{m,0}$. For $\vct{\xi}_m \notin\Gamma_{m,0}$, we directly compute $\vct{S}^m(\vct{\xi}_{m})$ and $\vct{b}^{m}(\vct{\xi}_{m})$ from \eqref{eqn:MsFEB2} and \eqref{eqn:stiffness_local}. The computed upscaled quantities, denoted as $\hat{\vct{S}}^{m}$ and $\hat{\vct{b}}^{m}$, give a very accurate approximation of its true values $\vct{S}^m$ and $\vct{b}^{m}$.
\end{itemize}
\end{remark}

\subsubsection{Complexity analysis}\label{sec:RandomInter_complexity}
We summarize the computational cost in this subsection. Without loss of generality, we assume that the diameter of the physical domain is 1. The coarse mesh size is denoted by $H$, and thus the number of coarse grid elements is $M\sim 1/H^d$. In each coarse grid element, we use a fine mesh of size $h$ to solve the cell problem. Given a sample $\kappa(x,\omega)$, we assume that the computational cost to solve a deterministic PDE~\eqref{eqn} is 
\begin{equation}\label{eqn:defmu}
	\mu = (1/h)^{d \gamma}, \quad \gamma \ge 1,
\end{equation}
where $\gamma = 1$ corresponds to the multigrid method (neglect the logarithmic factor). In the same manner, we assume that the computational cost to solve a upscale system~\eqref{eqn:upscaleSUF} is $M^{\gamma} \sim (1/H)^{d \gamma}$, and that the computational cost to solve a local cell problem is $(\eta H/h)^{d\gamma}$, where $\eta$ is the oversampling ratio. We also assume that the number of random variables in each coarse grid element is about the same, denoted by $K_m$. By an appropriate choice of $H$, $K_m$ could be as small as $2$ or $3$. The number of offline collocation points is denoted as $N_c=(\nu+1)^{K_m}$ for each coarse grid element.

The computational cost of the random interpolation method consists of the offline cost and the online cost. In the offline stage, see algorithm~\ref{alg:randominterp_offline}, we need to solve $L-1$ local cell problems and assemble local stiffness matrices for each collocation point on each coarse grid element.
\begin{itemize}
\item For each collocation point within a coarse grid element, we need to construct $L-1$ basis functions by solving an elliptic equation with a total number of $(\eta H/h)^d$ fine grid points, where $\eta$ is the oversampling ratio. The computational cost for this step is given by:
\begin{equation*}
	\Cost_\text{basis} = (L-1) M N_c (\eta H/h)^{d \gamma}\,.
\end{equation*}
The oversampling ratio $\eta$ is typically taken to be $\eta= 2$. 
\item The second step in the offline stage is to assemble the local stiffness matrix. Such procedure is performed on every collocation point within each coarse grid element and the computational cost is:
\begin{equation*}
	\Cost_\text{assemble} =\frac{(L+1)L}{2} M N_c (H/h)^d\,.
\end{equation*}
Here $(L+1)L/2$ is due to the fact that each local stiffness matrix has $(L+1)L/2$ different entries (other $(L-1)L/2$ elements are determined by symmetry), and the factor $(H/h)^d$ comes from evaluating the $l_2$ inner product defined in ~\eqref{eqn:stiffness_local} over a coarse grid element by using $(H/h)^d$ number of fine grid points.
\end{itemize}
In total, we have
\begin{equation}\label{eqn:costRI_offline}
\frac{\Cost_\text{offline}}{\mu} \sim (L-1)N_c H^{d(\gamma-1)} + \frac{(L+1)L}{2} N_c h^{d(\gamma-1)}.
\end{equation}
In the extreme case of $\gamma=1$, we have $\frac{\Cost_\text{offline}}{\mu} \sim O(N_c)$. If $\gamma > 1$, $\frac{\Cost_\text{offline}}{\mu} $ is even smaller. This is to say that compared with the multigrid method on the fine grid ($\gamma = 1$), the offline computational cost is the same as solving the original equation~\eqref{eqn} for about $N_c$ times. Here $N_c = (\nu + 1)^{K_m}$ is the number of local collocation points and is much smaller than the number of samples $N_{\text{on}}$ that is required to solve in the online stage. Therefore, the computational overhead of the random interpolation method is quite reasonable.

In the online stage, we need to interpolate the stiffness matrix and then solve the coarse grid system. If $N_{\text{on}}$ samples are computed, the computational cost is
\begin{equation*}
\Cost_\text{online} = N_{\text{on}}\left( \frac{(L+1)L}{2} M N_c + M^{\gamma} \right)\,.
\end{equation*}
Here the first term comes from the stiffness matrix interpolation and the second is to solve the upscaled linear system~\eqref{eqn:upscaleSUF}. For every sample in the online stage, we have
\begin{equation}\label{eqn:costRI_online0}
\frac{\Cost_\text{online}}{\mu} \sim \frac{(L+1)L}{2}N_c H^{d(\gamma-1)} (h/H)^{\gamma d} + (h/H)^{\gamma d},
\end{equation}
which is of the order $N_c (h/H)^{d}$ in the extreme case $\gamma = 1$ and is much smaller if $\gamma > 1$. The computational saving $(h/H)^{\gamma d}$ comes from the usage of MsFEM. However, since we need to do interpolation to get the upscaled system, we have to pay a factor of $N_c$ as the interpolation cost.

Since the cell problem is solved on the fine mesh $\CalT_h$, there is $\Or(h)$ error in the upscaled system due to this spatial discretization. Therefore, as long as the interpolation error in the stochastic space is smaller than $\Or(h)$, it will not influence the accuracy of the computed upscaled system. Due to the exponential decay of the interpolation error, see Theorem~\ref{thm:approx_S}, it is sufficient to choose the degree of interpolation polynomial $\nu \sim \log(H/h)$. Therefore, the online cost would be
 \begin{equation}\label{eqn:costRI_online}
\frac{\Cost_\text{online}}{\mu} \sim \frac{\left(\log(H/h)\right)^{K_m}}{(H/h)^{\gamma d}},
\end{equation}
which implies that we obtain a significant computational saving in the online stage if the local dimension $K_m$ is small. However, the computational saving quickly decreases as $K_m$ increases. In this case, one should use a sparse grid interpolation instead. Finally, the total computational cost for this random interpolation method is
 \begin{equation}\label{eqn:costRI}
	\frac{\Cost_\text{StoMsFEM}}{\mu} = N_{\text{off}} + R N_{\text{on}},
\end{equation}
where $N_{\text{off}} = \Or(N_c)$ is the effective number of samples we solve in the offline stage, and $R = \Or(N_c (h/H)^{\gamma d})$ is the online computational saving achieved by the random interpolation method.

\subsection{The reduced basis method} \label{sec:reducedbasis}
Besides exploring the regularity of the upscaled quantities $\vct{S}^m$ as what is done in the random interpolation method, another idea is to make use of the low rank property of multiscale basis $\phi^{m}(x,\omega)$. This leads to the design of the reduced basis method. As in the last section, we suppress the super-index $l$ in $\phi^{ml}$ in Eqn.~\eqref{eqn:MsFEB2} when no confusion arises.

\subsubsection{Reduced basis construction via KL expansion}
To obtain the reduced basis, we apply the KL expansion to $\phi^{m}(x,\omega) \in L^2(\Omega; H_0^1(D^m_H;\bar{\kappa}))$
\begin{equation} \label{eqn:KLphi}
	\phi^{m}(x,\omega) = \bar{\phi}^{m}(x) + \sum_{q=1}^{\infty}\sqrt{\lambda_q}\tau_q(\omega) \zeta^{m,q}(x)\,,  \quad \lambda_1 \ge \lambda_2 \ge \lambda_3 \ge ... 
\end{equation}
where $\bar{\phi}^{m}(x) = \int\phi^{m}(x,\omega)\rd{P}(\omega)$ is the mean of $\phi^{m}(x,\omega)$. Notice that this KL expansion is performed in the Hilbert space $L^2(\Omega; H_0^1(D^m_H;\bar{\kappa}))$ (which is isometric to  $L^2(\Omega) \times H_0^1(D^m_H;\bar{\kappa})$) to guarantee that the reduced basis method is accurate in $H_0^1(D^m_H;\bar{\kappa})$\footnote{$H_0^1(D^m_H;\bar{\kappa})$ is the function space $\{v \in H_0^1(D^m_H): \int_{D^m_H} \bar{\kappa} |\nabla v|^2 \rd x \le \infty\}$ with inner product $\langle u, v \rangle = \int_{D^m_H} \bar{\kappa} \nabla u \cdot \nabla v \rd x$.}. Thanks to the local low stochastic dimensionality, the energies $\{\lambda_q\}_{q=1}^{\infty}$ in \eqref{eqn:KLphi} decay exponentially fast. We take the first $Q$ KL modes $\{\zeta^{m,q}\}_{q=1}^Q$ as the reduced basis functions, and expand the solution to Eqn. \eqref{eqn:MsFEB} as
\begin{equation} \label{eqn:localRB}
	\phi_{rb}^{m}(x,\omega) = \bar{\phi}^{m}(x) + \sum_{q=1}^Q c_q^m(\omega) \zeta^{m,q}(x).
\end{equation}
Using the Galerkin method, we solve the following linear system to obtain the coefficients $\vct{c}^m$:
\begin{equation} \label{eqn:localRBsystem}
	\vct{A}^m(\omega) \vct{c}^m(\omega) = \vct{F}^m(\omega),
\end{equation}
where $\vct{A}^m(\omega)$ is a $Q\times Q$ symmetric positive definite matrix with entries $\vct{A}_{qq'}^m(\omega) = \kappa(\zeta^{m,q}, \zeta^{m,q'}; \omega)$, $1\le q,q' \le Q$, and $\vct{F}^m(\omega)$ is the load vector with entries $\vct{F}_q^m(\omega) = - \kappa(\bar{\phi}^{m}, \zeta^{m,q}; \omega)$, $1\le q \le Q$. 
Since the number of the reduced basis $Q$ is much smaller than the number of fine grid points, $(H/h)^d$, Eqn. \eqref{eqn:localRBsystem} can be solved very efficiently. Finally we use $\{\phi_{rb}^{ml}\}$ to build an approximation of the local stiffness matrix $\vct{S}^m$, denoted as $\hat{\vct{S}}^m$, as in \eqref{eqn:stiffness_local}.

To perform the KL expansion \eqref{eqn:KLphi}, we apply the stochastic collocation method to estimate the mean and covariance. Due to the locally low dimensionality, the stochastic collocation method requires much smaller number of samples than the MC method does, and thus accelerate the offline computation significantly. We have the following theorem that guarantees the accuracy and efficiency of the reduced basis method.
\begin{theorem} \label{thm:phi_S_RB}
Suppose we take the first $Q$ KL modes $\{\zeta^{m,q}(x)\}_{q=1}^Q$ as the reduced basis functions, and use the Galerkin method to obtain the reduced-basis solution of Eqn.~\eqref{eqn:MsFEB}, denoted as $\phi_{rb}^{m}(x,\omega)$ in ~\eqref{eqn:localRB}.
Assume that
\begin{enumerate}
\item there exists $C_1>0$ and $\beta>1$ such that $\lambda_j \le C_1 \beta^{-j}$, and
\item there exists a constant $C_2$ such that $\kappa(x,\omega) \le C_2\bar{\kappa}(x)$ for all realizations $\omega\in \Omega$.
\end{enumerate}
Then we have for any $\varepsilon>0$,
\begin{enumerate}
\item
\begin{equation}\label{eqn:phiaccuracy}
P\left[ \|\phi^{m}(x,\omega) - \phi_{rb}^{m}(x,\omega)\|_{H_0^1(D_H^m;\kappa)} \ge \varepsilon \right] \le \frac{C_1 C_2^2 \beta^{-Q}}{(\beta -1)\varepsilon^2}\,;
\end{equation}
\item for any $l, l' \in \{1,2, \dots, L\}$
\begin{equation}\label{eqn:Saccuracy}
	P\left[ |\vct{S}_{ll'}^{m}(\omega) - \hat{\vct{S}}_{ll'}^{m}(\omega)| \ge 2 C_{\phi} \sqrt{\kappa_{max}} \varepsilon \right] \le \frac{2 C_1 C_2^2 \beta^{-Q}}{(\beta -1)\varepsilon^2}.
\end{equation}
\end{enumerate}
Here, $H_0^1(D;\kappa)$ is the Hilbert space with norm $\int_D \kappa(x,\omega) |\nabla f(x)|^2 \rd x$; $C_{\phi}$ is a constant that only depends on the local domain $D^m_H$ and the deterministic boundary condition $p^{ml}$ for $\phi^{ml}$.
\end{theorem}
We point out that the first assumption in Theorem~\ref{thm:phi_S_RB} holds true in general. Moreover, under the same assumptions in Lemma~\ref{lem:smoothness_basis}, we can prove that the smaller the $c_k$ is, the bigger $\beta$ is. This exponential decay is also observed in our numerical examples, see Figure~\ref{fig:stiffness_99_RB}. The second assumption is also valid in general. We will demonstrate that this assumption is satisfied in our numerical example with a high contrast random medium, see Section~\ref{exam:highcontrast}. We will not present the proof of Theorem~\ref{thm:phi_S_RB} in this paper and refer the interested reader to ~\cite{thesis_zhang_2016}.

Theorem~\ref{thm:phi_S_RB} guarantees that for any pre-specified $\varepsilon > 0$ and $\delta > 0$, with only $Q=\Or(\log(1/\varepsilon)+\log(1/\delta))$ reduced basis functions, our reduced basis approximation $\hat{\vct{S}}^{m}$ is $\Or(\varepsilon)$ accurate with probability at least $1-\delta$. There are two ways to deal with the rare event when our approximation is not guaranteed to be accurate. In the first approach, with appropriate {\it a posteriori} error estimate, see~\cite{rozza2008reduced}, we are able to efficiently detect the small-probability failure samples, and recompute these samples directly to make sure they are accurate. In the second approach, we do not care about this rare event at all because this small probability error only introduces a small error in estimating statistical properties of the solution $u_H$. In this paper, we take the second approach, see also \cite{jakeman2010numerical, chen2013flexible}.
\begin{remark} \label{rem:unknowndist}
In the case when the distribution of local parameters $\vct{\xi}_m$ is unknown, we need to choose an auxiliary distribution, its density denoted as $\hat{\rho}_m(\vct{\xi}_m)$, to do the KL expansion~\eqref{eqn:KLphi}. In practice, one can just the use uniform distribution for bounded variables. For unbounded variables, one can choose a sufficient large square domain, which covers the range of local parameters $\vct{\xi}_m$ with high probability, and then use the uniform distribution on the square domain. Theoretically, we can prove that when there exists a constant $C>0$ such that $\rho_m/\hat{\rho}_m \le C$, our reduced basis approximation $\hat{\vct{S}}^{m}$ is still $\Or(\varepsilon)$ accurate with probability at least $1-\delta$ with only $Q=\Or(\log(1/\varepsilon)+\log(1/\delta))$ reduced basis functions. In practice, one can simply take $Q = 3 K_m$ where $K_m$ is the local stochastic dimension.
\end{remark}

\subsubsection{Exploring the affine structure of the coefficient for further speedup}
In the online stage, the reduced basis method described above still evaluates $\phi^{m}(x,\omega)$ point-wisely according to \eqref{eqn:KLphi} and assembles local stiffness matrix \eqref{eqn:stiffness_local} on the fine grid. Even in the case when $\phi^{m}(x,\vct{\xi}_m)$ is prepared by gPC expansions as in \cite{doi:10.1137/140970100}, the evaluation of $\phi^{m}(x,\omega)$ and numerical integration are still performed on the fine grid, which offers little computational saving compared with the multigrid method. To make the upscaling step more efficient, we assume that the parametrization of $\kappa(x, \vct{\xi}_m)$ is affine with respect to the parameters $\vct{\xi}_m$. With this assumption, we can pre-compute the essential part of the stiffness matrix in the offline stage, which leads to considerable saving in assembling the stiffness matrix for each sample in the online stage. Specifically, we assume that the local random coefficient $\kappa(x, \vct{\xi}_m)$ can be expressed as follows:
\begin{equation}\label{eqn:RBaffine}
	\kappa(x, \vct{\xi}_m) = \sum_{k=1}^{K_m} \xi_{m,k} \kappa_{m,k}(x), \quad x \in D^m_H.
\end{equation}
By applying the affine structure of the coefficient, we obtain
\begin{equation} \label{eqn:localAF}
	\vct{A}^m(\omega) = \sum_{k=1}^{K_m} \xi_{m,k}(\omega) \vct{A}_k^m, \quad \vct{F}^m(\omega) = \sum_{k=1}^{K_m} \xi_{m,k}(\omega) \vct{F}_k^m,
\end{equation}
where the deterministic coefficients $\vct{A}_k^m$ and $\vct{F}_k^m$ are given by $\vct{A}_{k,q q'}^m = \kappa_{m,k}(\zeta^{m,q}, \zeta^{m,q'})$, $1\le q,q' \le Q$ and $\vct{F}_{k,q}^m = - \kappa_{m,k}(\bar{\phi}^{m}, \zeta^{m,q})$, $1\le q \le Q$. We can precompute $\vct{A}_k^m$ and $\vct{F}_k^m$ and efficiently assemble the stiffness matrix $\vct{A}^m(\omega)$ and load vector $\vct{F}^m(\omega)$ for each sample. We remark that the affine structure also simplifies the assembling of local stiff matrix $\vct{S}^m$ and loading vector $\vct{b}^m$. 

Finally, we summarize the offline and online implementation of random interpolation method in Algorithm~\ref{alg:reducedbasis_offline} and Algorithm~\ref{alg:reducedbasis_online}.
\begin{algorithm}
\caption{The offline stage of reduced basis method}\label{alg:reducedbasis_offline}
\begin{algorithmic}[1]
\State Partition the physical domain $D$ into coarse grid blocks $\mathcal{T}_{H}$
\For{each patch $D^m_H$}
\State Solve local multiscale basis \eqref{eqn:MsFEB} with affine coefficient \eqref{eqn:RBaffine} by stochastic collocation method to obtain samples $\{\phi^{ml}(x,\vct{\xi}_{m,c})\}_{c=1}^{N_c}$
\State Apply KL expansion to get reduced basis $\{\zeta^{ml,q}\}_{q=1}^Q$
\EndFor
\end{algorithmic}
\end{algorithm}

\begin{algorithm}
\caption{The online stage of reduced basis method for specific parameter configuration $\vct{\xi}$}\label{alg:reducedbasis_online}
\begin{algorithmic}[1]
\For{each patch $D^k_H$}
\State Determine the values of local variables $\vct{\xi}_{m}$
\State Assemble and solve local reduced systems \eqref{eqn:localRBsystem}
\State Assemble the local upscaled stiffness matrix $\hat{\vct{S}}^m$ and loading vector $\hat{\vct{b}}^m$
\EndFor
\State Assemble and solve the upscaled system \eqref{eqn:upscaleSUF} to obtain $\hat{u}_H(x,\vct{\xi})$
\end{algorithmic}
\end{algorithm}

\subsubsection{Complexity analysis}\label{sec:costanalysisRB}
Using the same notations and assumptions as in Section~\ref{sec:RandomInter_complexity}, we analyze the computational cost of the StoMsFEM with local reduced basis in this section. As for the random interpolation method, the computational cost consists of offline and online parts.

The offline cost consists of three parts: obtaining samples for local cell problems \eqref{eqn:MsFEB}, performing the KL expansion to get reduced basis \eqref{eqn:KLphi} and assembling upscaled quantities.
At the coarse grid level, we have about $M$ quadrilateral coarse grid elements and on each element we solve $N_c$ samples to do the KL expansion \eqref{eqn:KLphi}. The cost of obtaining these solution samples is $M N_c (\eta H/h)^{d \gamma}$. The cost of obtaining the first $Q$ KL modes in \eqref{eqn:KLphi} is of order $N_f N_c \log(Q) + Q^2 (N_f + N_c)$, see~\cite{halko2011finding, cheng_data-driven_2013}. Since $N_f = (\eta H/h)^{d}$, this part of the cost is $\Or\left( (N_c \log(Q) + Q^2) (\eta H/h)^{d} + Q^2 N_c \right)$. 
Finally, the cost of assembling the upscaled stiffness matrix and the loading vector is about $(\sum_{m=1}^M K_m Q^2) (\eta H/h)^{d}$. Assuming that $Q = \Or(K_m)$, $\eta = \Or(1)$ and that all the local dimensions are about the same, we have:
\begin{equation*}
\Cost_\text{offline} \sim N_c H^{d(\gamma - 1)} h^{-d\gamma} + (N_c\log(K_m) + K_m^2) h^{-d} + K_m^2 N_c H^{-d} + K_m^3 h^{-d}.
\end{equation*}
Therefore, the ratio between $\Cost_\text{offline}$ and $\mu$ is bounded by
\begin{equation}\label{eqn:costRB_offline}
\frac{\Cost_\text{offline}}{\mu} \underset{\sim}{<} N_c H^{d(\gamma - 1)} + (N_c\log(K_m) + K_m^2 + K_m^3) h^{d(\gamma-1)} + K_m^2 N_c H^{d(\gamma-1)}(H/h)^{-d\gamma}.
\end{equation}
For $\gamma = 1$, we can see that the ratio is of the order $(N_c\log(K_m) + K_m^3) $. 

In the online stage, for a given configuration of the random parameters $\vct{\xi}$, the cost also consists of 3 parts, assembling and solving the reduced basis system \eqref{eqn:localRBsystem}, assembling the local upscaled stiffness matrix and loading vector,
 and finally globally assembling and solving the upscaled system \eqref{eqn:upscaleSUF}. For each sample, the computational cost is
\begin{equation*}
\Cost_\text{online} \sim M (K_m Q^2 + Q^2) + M K_m Q^2 + M + M^{\gamma}\,.
\end{equation*}
Here the first term comes from assembling and solving the reduced basis system \eqref{eqn:localRBsystem}; the second term is for assembling the local upscaled stiffness matrix and loading vector; 
 the third and forth term is from globally assembling and solving the upscaled system \eqref{eqn:upscaleSUF}. In practice, we observe that $Q = \Or(K_m)$, and then the ratio between $\Cost_\text{online}$ and $\mu$ is give by
\begin{equation}\label{eqn:costRB_online}
\frac{\Cost_\text{online}}{\mu} \underset{\sim}{<} \frac{1+K_m^3 H^{d(\gamma-1)}}{ (H/h)^{d\gamma}}.
\end{equation}
For $\gamma=1$, the ratio is about $K_m^3 (h/H)^{d}$ where $(h/H)^{d}$ comes from the usage of the MsFEM and the factor $K_m^3$ comes from assembling and solving the reduced basis systems. Again, the computational saving is more significant for $\gamma>1$. Finally, when we solve $N_{\text{on}}$ samples in the online stage, the total computational cost for this reduced basis method is
 \begin{equation}\label{eqn:costRB}
	\frac{\Cost_\text{StoMsFEM}}{\mu} = N_{\text{off}} + R N_{\text{on}},
\end{equation}
where $N_{\text{off}} = \Or(N_c\log(K_m) + K_m^3)$ is the effective number of samples that we solve in the offline stage, and $R = \Or(K_m^3 (h/H)^{\gamma d})$ is the online computational saving achieved by the reduced basis method.

\section{Global Stochastic Methods}\label{sec:globalstochastic}
The StoMsFEM is designed to compute an approximate solution $\hat{u}_H(x, \vct{\xi})$ for every parameter configuration $\vct{\xi}$ efficiently, as described in Algortithm~\ref{alg:randominterp_online} and Algorithm~\ref{alg:reducedbasis_online}. It is straightforward to combine StoMsFEM with any non-intrusive stochastic method, which determines the sample set to solve in the online stage, to finally estimate the statistical properties of the coarse grid solution. In subsection~\ref{sec:globalmethods}, we combine StoMsFEM with the MC method and the (sparse grid) stochastic collocation method. In subsection~\ref{sec:globalerror}, we show that to achieve the same level of estimation error, compared with the standard FEM on fine grid, StoMsFEM indeed offers significant computational saving by optimally balancing the spatial discretization error from MsFEM and the stochastic sampling error from the global stochastic methods.

\subsection{Global stochastic methods}\label{sec:globalmethods}
\subsubsection{Global Monte Carlo method} \label{sec:gMonteCarlo}
The Monte Carlo method estimates statistical properties by ensemble average, i.e.
\begin{equation}\label{eq:MCestimate}
	 \EE[f(\omega)] \approx \mathcal{M} [f(\omega)] \equiv \frac{1}{N_{on}} \sum\limits_{\omega_i \in \mathcal{S}} f(\omega_i),
\end{equation}
where $\omega_i$ is the $i$-th sample and $N_{on}$ is the total number of independent samples $\mathcal{S}$. This can be used to approximate the mean value or the variance of $u_H(x,\omega)$ as
\begin{equation*}
\begin{split}
	\EE[u_{H}(x,\omega)] &\approx \mathcal{M} \left[ u_H(x,\omega) \right],\\
	\mathrm{var} [u_{H}(x,\omega)] &\approx \mathcal{M} \left[ (u_H(x,\omega))^2 \right] - \left\{ \mathcal{M} \left[ u_H(x,\omega) \right]\right\}^2.
\end{split}
\end{equation*}
For sample $\omega_i$, Algorithm~\ref{alg:randominterp_online} or Algorithm~\ref{alg:reducedbasis_online} can be applied to compute the solution $\hat{u}_H(x,\omega_i)$ on the coarse mesh $\CalT_H$.

Due to the probabilistic nature of the MC estimator~\eqref{eq:MCestimate}, we use the mean square error (MSE) to quantify its performance. For example, to estimate $\EE[u(x,\omega)]$ by $\mathcal{M}[\hat{u}_H]$, simple calculation shows that the MSE can be written as
\begin{equation}\label{eqn:umeanMSE}
\EE_{\mathcal{S}}\left[ \| \mathcal{M}[\hat{u}_H] -  \EE[u(x,\omega)\|_2^2 \right] = \left\| \EE[\hat{u}_H(x,\omega)] - \EE[u(x,\omega)] \right\|_2^2 + \frac{1}{N_{on}} \int_D \text{var}[ \hat{u}_H(x,\omega) ] \rd x.
\end{equation}
Here, $\EE_{\mathcal{S}}$ is the expectation taking w.r.t. to the random ensemble $\mathcal{S}$, the first part is the spatial discretization error introduced by StoMsFEM, and the second part is the sampling error introduced by the MC method. 

\begin{remark}\label{rem:MLMC}
We can also consider the following two-level MC estimator~\cite{giles2008multilevel, barth2011multi, cliffe2011multilevel}
\begin{equation}\label{eq:MCmean2level}
	\mathcal{M}^{(2)}[\hat{u}_H] = \frac{1}{N_{on,H}} \sum_{i=1}^{N_{on,H}} \hat{u}_H(x,\omega_{i,H}) + \frac{1}{N_{on,h}} \sum_{i=1}^{N_{on,h}} \left[ u_h(x,\omega_{i,h}) - \hat{u}_H(w,\omega_{i,h}) \right],
\end{equation}
where $\{\omega_{i,H}\}_{i=1}^{N_{on,H}}$ and $\{\omega_{i,h}\}_{i=1}^{N_{on,h}}$ are independent samples. Simple calculation shows that its MSE is
\begin{equation}\label{eqn:mseML}
	\left\| \EE[u_h(x,\omega)] - \EE[u(x,\omega)] \right\|_2^2 + \frac{1}{N_{on,H}} \int_D \text{var}[ \hat{u}_H(x,\omega) ] \rd x + \frac{1}{N_{on,h}} \int_D \text{var}[ \hat{u}_H(x,\omega) - u_h(x,\omega) ] \rd x.
\end{equation}
Compared with~\eqref{eqn:umeanMSE}, the two-level MC estimator is able to reduce its MSE to the order of fine grid discretization error by properly choosing $N_{on,H}$ and $N_{on,h}$. At the same time, its computational cost, including computing $N_{on,H}+N_{on,h}$ coarse grid solutions and $N_{on,h}$ fine grid solutions, can be significantly smaller than that of the standard MC method. The complexity analysis and comparison with the standard MC method are provided in the last paragraph of Section~\ref{sec:globalerror}.
\end{remark}

\subsubsection{Global stochastic collocation methods} \label{sec:globalSC}
When the density $\rho(\vct{\xi})$ is known, the stochastic collocation (SC) methods~\cite{babuska_stochastic_2007, nobile_sparse_2008, nobile_anisotropic_2008, doi:10.1137/040615201} may have better convergence rate when approximating the expectation (multivariate integral) in some cases. To illustrate the idea, we assume that the random variables $\vct{\xi}$ are independent. In this case, their joint density $\rho$ factorizes as $\rho(\xi_1, \dots, \xi_K) = \prod_{k=1}^K \rho_k(\xi_k)$.

Similar to the Monte Carlo estimator \eqref{eq:MCestimate}, the global SC method solves the parametrized problem~\eqref{eqn} on a {\it deterministic} set of collocation points, denoted as $\text{CN}(K)$, and then approximates $\EE[f(\vct{\xi})]$ by some {\it deterministic} numerical quadrature rule, i.e.,
\begin{equation}\label{eqn:sparseintegration}
	\EE[f(\vct{\xi})] \approx \mathcal{I} [f(\vct{\xi})] \equiv \sum\limits_{\vct{\xi} \in \text{CN}(K)} w(\vct{\xi}) f(\vct{\xi}). 
\end{equation} 
This can be used to approximate the mean value or the variance of $u_H(x,\omega)$ as
\begin{equation*}
\begin{split}
	\EE[u_{H}(x,\omega)] &\approx \mathcal{I} \left[ u_H(x,\vct{\xi}) \right],\\
	\mathrm{var} [u_{H}(x,\omega)] &\approx \mathcal{I} \left[ (u_H(x,\vct{\xi}))^2 \right] - \left\{ \mathcal{I} \left[ u_H(x,\vct{\xi}) \right]\right\}^2.
\end{split}
\end{equation*}
For the standard SC method~\cite{babuska_stochastic_2007}, $\text{CN}(K)$ is a tensor product grid of all the one-dimensional collocation points. In our case when the global stochastic dimension $K$ is large, the sparse grid SC method~\cite{nobile_sparse_2008, nobile_anisotropic_2008} is preferred, where $\text{CN}(K)$ is a high-dimensional sparse grid. 

The locally low dimensionality can offer a huge computational saving for the SC method. The key observation is that the global collocation points $\text{CN}(K)$ repeatedly use the local collocation points. For example, the tensor product grid collocation points $\prod\limits_{k=1}^K \Gamma_k$, where $\Gamma_k$ is the collocation points in the $\xi_k$ direction, reuse the local collocation points $\prod\limits_{\xi_k \in \vct{\xi}_m} \Gamma_k$, which is the local tensor product grid with the same degree. In the same manner, the local collocation points of a global sparse grid is still a sparse grid of the low dimensional local parameter space. 

Therefore, the global (sparse grid) SC method also contains the offline and online stage. In the offline stage, we take the local interpolation nodes $\Gamma_m$ to be the local collocation points corresponding to the global (sparse grid) collocation points $\text{CN}(K)$, and run Algorithm~\ref{alg:randominterp_offline}. Our estimate in Eqn.~\eqref{eqn:costRI_offline} implies that
\begin{equation}\label{eqn:SCofflinecost}
\frac{\Cost_\text{offline}}{\mu} \sim (L-1)N_c H^{d(\gamma-1)} + \frac{(L+1)L}{2} N_c h^{d(\gamma-1)},
\end{equation}
where $N_c$ is the number of local collocation points. In the extreme case of $\gamma=1$, we have $\frac{\Cost_\text{offline}}{\mu} \sim O(N_c)$.

The algorithm in the online stage is almost the same with Algorithm~\ref{alg:randominterp_online}, but we simply do search-and-plug-in instead of interpolation for every collocation point $\vct{\xi} \in \text{CN}(K)$. Since searching cost is typically negligible, the online computational cost for {\it each collocation point} only contains assembling and solving the upscaled system~\eqref{eqn:upscaleSUF}, i.e.
\begin{equation} \label{eqn:SConlinecost}
	\frac{\Cost_\text{online}}{\mu} \approx \frac{M^{\gamma}}{\mu} \sim (h/H)^{\gamma d}\,.
\end{equation}
Compared with Eqn.~\eqref{eqn:costRI_online0}, we do not have the interpolation cost when the global (sparse grid) SC solver is utilized. This is a big difference between the global SC solver and the global MC solver when they are combined with the local random interpolation method.

We still use $N_{on}$ to denote the number of samples to solve in the online stage, which is $|\text{CN}(K)|$ in the global SC solver. Therefore, the total computational cost for the StoMsFEM with the global SC solver is 
\begin{equation}\label{eqn:costSC}
	\frac{\Cost_\text{StoMsFEM}}{\mu} \approx N_{c} + R N_{on},
\end{equation}
where the online saving factor $R = (h/H)^{\gamma d}$. On the other hand, the total computational cost for the standard FEM on the fine grid $\CalT_h$ with a sparse grid collocation is $\mu N_{on}$. Notice that $N_{c}$ is the number of local sparse grid collocation points, which is nearly negligible compared with the number of global collocation points $N_{on}$, thanks to the locally low dimensionality. Therefore, we get a computational saving with nearly a factor of $R = (h/H)^{\gamma d}$.

%

The estimation error of the SC method is determined by the error of the numerical quadrature, i.e. $\mathcal{I} [f(\vct{\xi})] - \EE[f(\vct{\xi})]$. For example, to estimate $\EE[u(x,\omega)]$ by $\mathcal{I}[\hat{u}_H]$, the estimation error can be bounded as follows:
\begin{equation}\label{eqn:umeanSCerror}
\| \mathcal{I}[\hat{u}_H] -  \EE[u(x,\vct{\xi})\|_2^2 \le 2 \left\| \EE[\hat{u}_H(x,\omega)] - \EE[u(x,\omega)] \right\|_2^2 + 2 \| \mathcal{I}[\hat{u}_H] -  \EE[\hat{u}_H(x,\omega)]\|_2^2.
\end{equation}
Here, the first part is the spatial discretization error introduced by StoMsFEM, and the second part is the sampling error introduced by the SC method. 

\begin{remark}
In a general multivariate problem, if the random variables $\vct{\xi}$ are not independent, the density $\rho$ does not factorize, i.e., $\rho(\xi_1, \dots, \xi_K) \neq \prod_{k=1}^K \rho_k(\xi_k)$. To this end, we first introduce an auxiliary probability density function $\hat{\rho}: \RR^K \to \RR$ that can be seen as the joint density of $K$ independent random variables, i.e., it factorizes as $\hat{\rho}(\xi_1, \dots, \xi_K) = \prod_{k=1}^K \hat{\rho}_k(\xi_k)$ and satisfies $\frac{\rho(\vct{\xi})}{\hat{\rho}(\vct{\xi})} \le C$ for a positive constant $C$. For each dimension $k = 1, 2, \dots, K$, the 1d collocation nodes $\mathcal{V}_k^{i}$ can be the Gaussian abscissas of $\hat{\rho}_k$ or nested abscissas associated with $\hat{\rho}_k$. The auxiliary density $\hat{\rho}$ should be chosen as close to the true density $\rho$ as possible, so that the quotient $\rho/\hat{\rho}$ remains bounded.
\end{remark}

\subsection{Global error analysis}\label{sec:globalerror}
The estimation error of both the MC method and the SC method consists of the spatial discretization error from StoMsFEM and the sampling error from the corresponding global stochastic methods, see Eqn.~\eqref{eqn:umeanMSE} and \eqref{eqn:umeanSCerror}. We should balance these two kinds of errors to achieve the optimal estimate within our budget of computing resources. To further analyze the estimation error, we assume the following estimates:
\begin{eqnarray}
\left\| \EE[u_h(x,\omega)] - \EE[u(x,\omega)] \right\|_2^2 \underset{\sim}{<} h^{\beta}, \quad  \left\| \EE[\hat{u}_H(x,\omega)] - \EE[u(x,\omega)] \right\|_2^2 \underset{\sim}{<} H^{\beta},\label{eqn:errorFEM}\\
\| \mathcal{I}(u_h) - \EE[u_h] \|_2 \underset{\sim}{<} N_{on}^{-\zeta}, \quad \| \mathcal{I}(\hat{u}_H) - \EE[\hat{u}_H] \|_2 \underset{\sim}{<} N_{on}^{-\zeta}, \label{eqn:errorSto}\\
\int_D \text{var}[ \hat{u}_H(x,\omega) ] \rd x \approx \int_D \text{var}[u_h(x,\omega)] \rd x \approx \int_D \text{var}[ u(x,\omega) ] \rd x = c_1,\label{eqn:variance0}
\end{eqnarray}
The rate $\beta$ in~\eqref{eqn:errorFEM} characterizes the discretization error from the standard FEM on fine mesh $\mathcal{T}_h$ and MsFEM on coarse mesh $\mathcal{T}_H$, and $\beta \approx 4$ in our case. The rate $\zeta$ in~\eqref{eqn:errorSto} characterizes the sampling error from the (sparse grid) SC method, and it is typically very small in our high stochastic dimension case. For some problems with moderate stochastic dimensions, $\zeta$ can be relatively large. For example, in our high contrast example the SC with the sparse Clenshaw-Curtis formulas, we observe $\zeta \approx 5$, see Figure~\ref{fig:errorcompare}. We assume that $\int_D \text{var}[ u(x,\omega) ] \rd x = \Or(1)$ in~\eqref{eqn:variance0}. Error analysis of standard FEM gives $\|u_h - u\|_{L^2(D\times \Omega)} = \Or(h^2)$. For any successful upscaling method, we expect $\|u_H - u\|_{L^2(D\times \Omega)} = \Or(H^2)$. For example, Theorem~\ref{thm:MsFEM} validates this for MsFEM on periodic random coefficients with period $\epsilon \ll H$. Other local upscaling methods~\cite{owhadi_polyharmonic_2014, maalqvist2014localization, owhadi2015multi} satisfy this assumption on much richer set of random coefficients, and the StoMsFEM can be adapted to work with them. Therefore, we have $\int_D \text{var}[ u_h(x,\omega) ] - \int_D \text{var}[ u(x,\omega) ] \rd x = \Or(h^2)$ and $\int_D \text{var}[ u_H(x,\omega) ] - \int_D \text{var}[ u(x,\omega) ] \rd x = \Or(H^2)$, and thus validate the assumption~\eqref{eqn:variance0}.

In~\eqref{eqn:errorFEM}, \eqref{eqn:errorSto} and \eqref{eqn:variance0}, we assume that $\hat{u}_H - u_H$ is negligible. This is reasonable since we can easily drive the error $\hat{u}_H - u_H$ below other errors due to its exponential decay implied by Theorem~\ref{thm:approx_S} in the random interpolation setting and Theorem~\ref{thm:phi_S_RB} in the reduced basis setting. Our numerical examples also validate this assumption.

Combing the above assumptions and Eqn.~\eqref{eqn:umeanMSE} and \eqref{eqn:umeanSCerror}, we need $N_{on} = \Or(H^{-\beta})$ for the Monte Carlo method and $N_{on} = \Or(H^{-\beta/\zeta})$ for the SC method to achieve $\Or(H^{\beta})$ estimation error. Notice that the number of samples required keeps the same for the standard FEM on fine grid when it aims to achieve the same $\Or(H^{\beta})$ estimation error. Since in high stochastic dimensional problems the decay rate of the physical discretization error is much faster than that of the sampling error, i.e., $\beta$ and $\beta/\zeta$ is large, the number of samples to be solved in the online stage is huge. For example, the MC method requires about 100,000,000 samples when we take $H = 0.01$ for a physical domain with $\Or(1)$ size. Compared with this huge number, the effective number of samples $N_{\text{off}}$ in the StoMsFEM offline stage, which is roughly equal to the number of local interpolation nodes, is negligible. Therefore, to achieve $\Or(H^{\beta})$ estimation error, the total computational cost ratio between StoMsFEM and the standard FEM on fine mesh, i.e. $N_{\text{off}}/N_{\text{on}} + R$, is nearly $R$. As we derived in the previous sections, $R$ is $\Or(N_c (h/H)^{\gamma d})$ for the random interpolation method, $\Or(K_m^3 (h/H)^{\gamma d})$ for the reduced basis method and $\Or((h/H)^{\gamma d})$ for the global SC method.

If we want to reduce the estimation error to the level of $\Or(h^{\beta})$,  we can combine StoMsFEM with the two-level MC estimator~\eqref{eq:MCmean2level}. Similar to the Multi-Level Monte Carlo method (MLMC), we reduce the variance part in~\eqref{eqn:mseML} to $\Or(h^{\beta})$ while optimally distributing computing resources to the coarse and fine grid computations. If we assume that $\int_D \text{var}[ \hat{u}_H(x,\omega) - u(x,\omega) ]\rd x  \underset{\sim}{<} H^{\alpha}$ that characterizes the variance reduction effect of $\hat{u}_H$, the ratio of total computation cost between this two-level MC estimator and the MC based on the standard FEM on the fine grid is $\Or((R^2 + H^{\alpha/2})^{1/2})$, where $R$ is the cost ratio as before.

\section{Numerical examples}\label{sec:numerics}
In this section we demonstrate the accuracy and efficiency of the proposed StoMsFEM. All our computations are performed using MATLAB R2015a (64-bit) on an Intel(R) Core(TM) i7-3770 (3.40 GHz).

\subsection{Patch study of a synthetic 2d example} \label{sec:patchstudy}
\begin{figure}
\centering
\includegraphics[width = 0.3\textwidth]{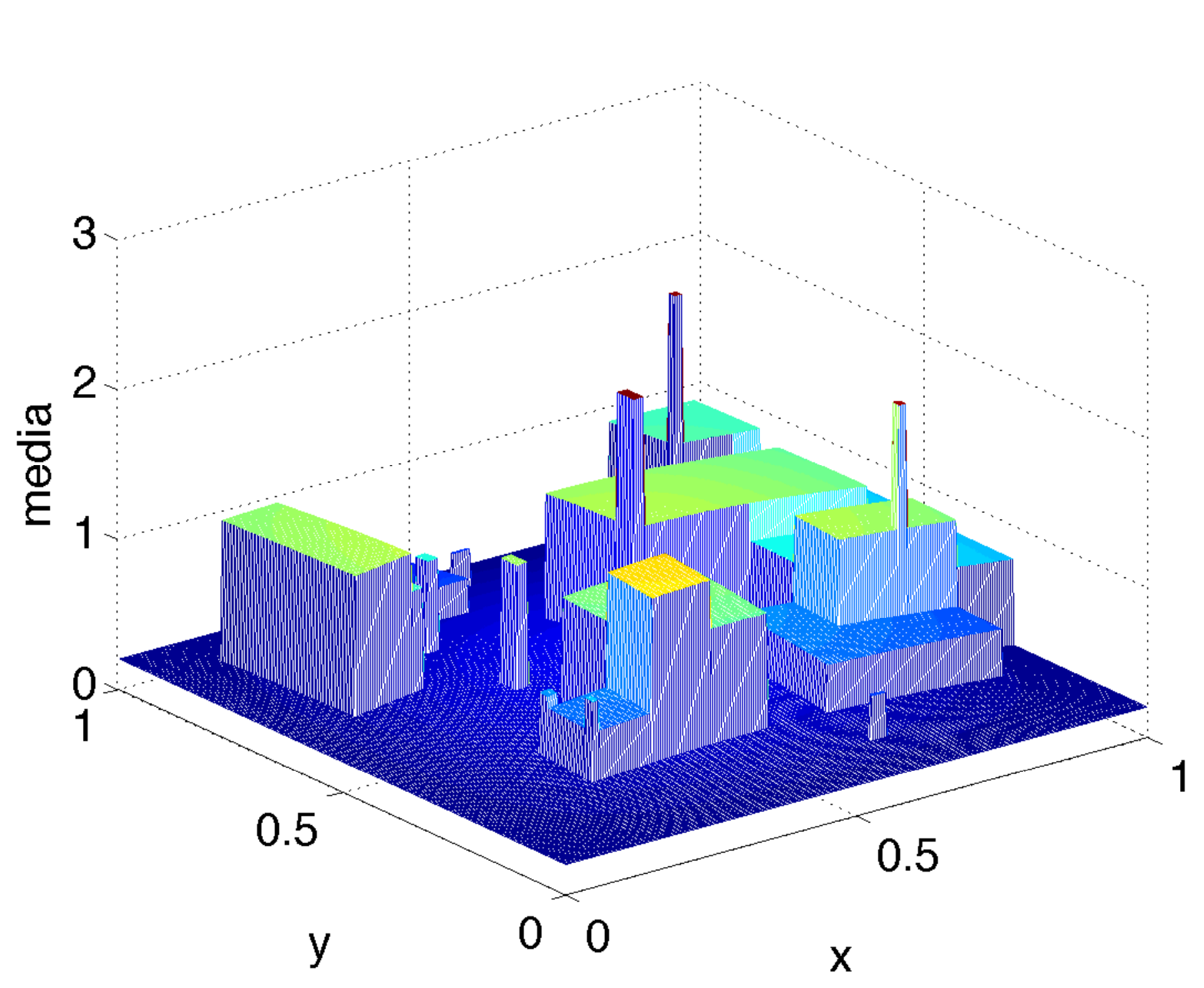}
\includegraphics[width = 0.3\textwidth]{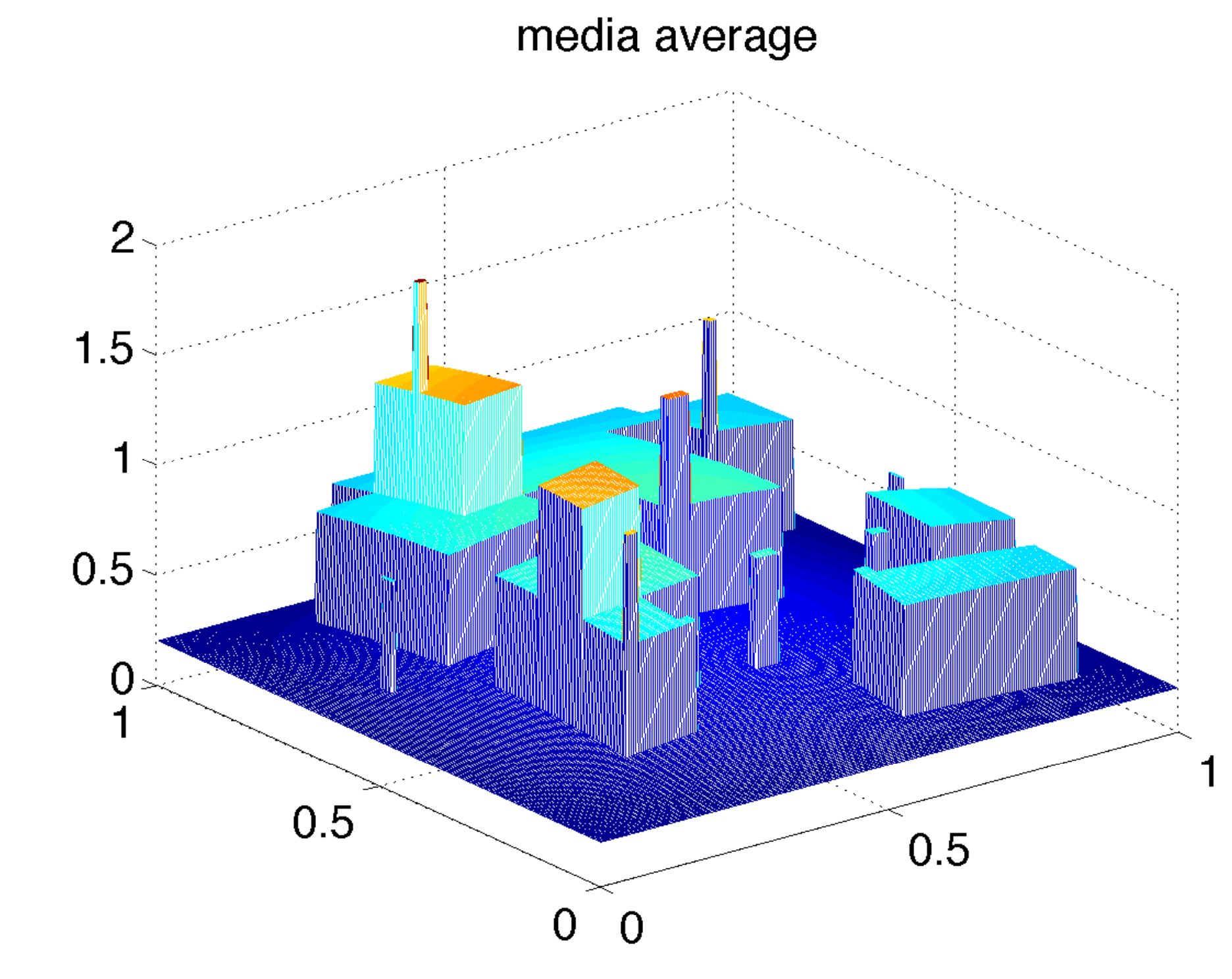}
\includegraphics[width = 0.3\textwidth]{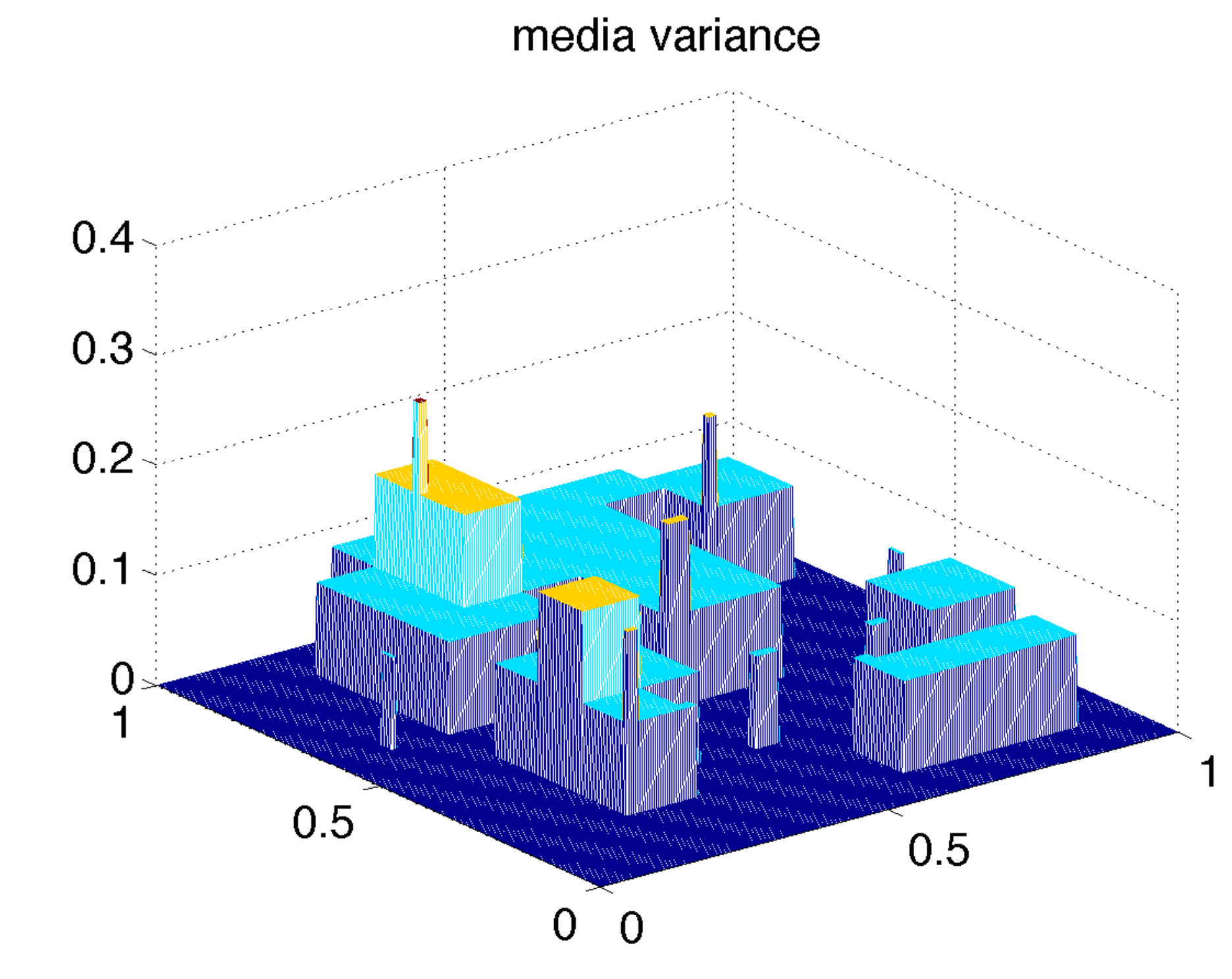}
\caption{The three figures respectively show one sample of the media, the mean and the variance of the media.}\label{fig:media}
\end{figure}

\begin{figure}
\centering
\includegraphics[width = 0.4\textwidth]{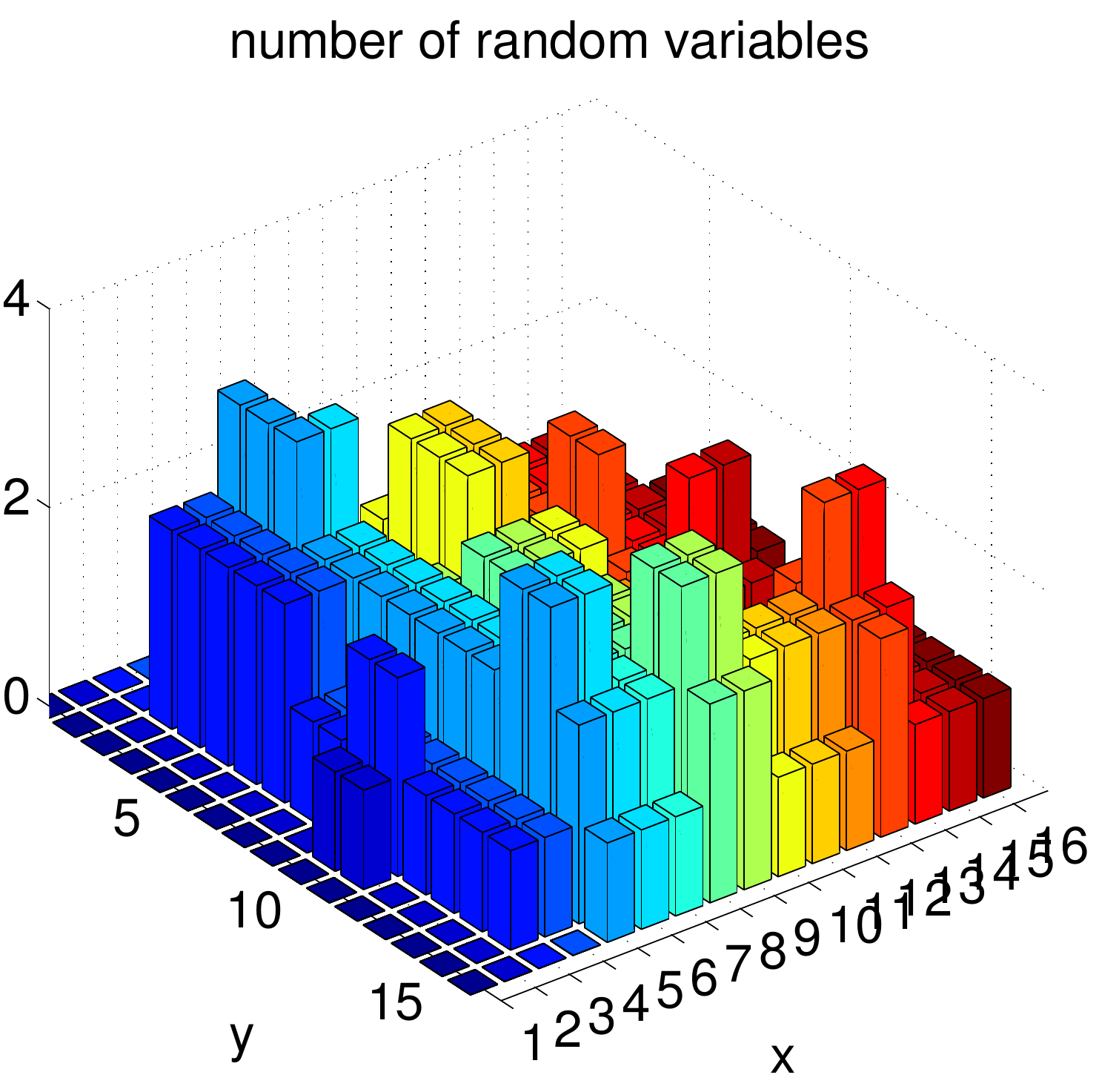}
\caption{The number of random variables seen in each patch. The domain is decomposed into $16\times 16$ coarse grid elements.}\label{fig:random_distribution}
\end{figure}

This synthetic example is adopted from problems with porous media~\cite{galvis_domain_2010} where the medium contains some channels and inclusions:
\begin{equation*}
\kappa(x,\omega) = 0.2+0.2\sin(\pi x)\sin(\pi y) + \sum_{k=1}^{20}\kappa_m(x)\xi_m\,.
\end{equation*}
Here the first two terms give the background of the medium, and in the summation $\kappa_m(x)$ are the characteristic functions representing the channels/inclusions and $\xi_m$ are the associated random variables. In our computation we set them uniformly distributed in $[0,1]$. We plot one sample of the medium in Figure~\ref{fig:media} and show the mean and the variance of the medium. It is easy to see that the medium contains many small sized inclusions, making the multi-scale treatment necessary.
For this two dimensional problem, we first decompose it into $16\times 16$ coarse grid elements. The oversampling ratio is chosen as $\eta = 2$, meaning each patch is enlarged by $2$ in each dimension for the oversampling. In Figure~\ref{fig:random_distribution} we plot the number of random variables in each patch. As shown in Figure~\ref{fig:random_distribution}, in each patch, there are about $2$ to $3$ random variables.

\begin{figure}
\centering
\includegraphics[width = 0.5\textwidth]{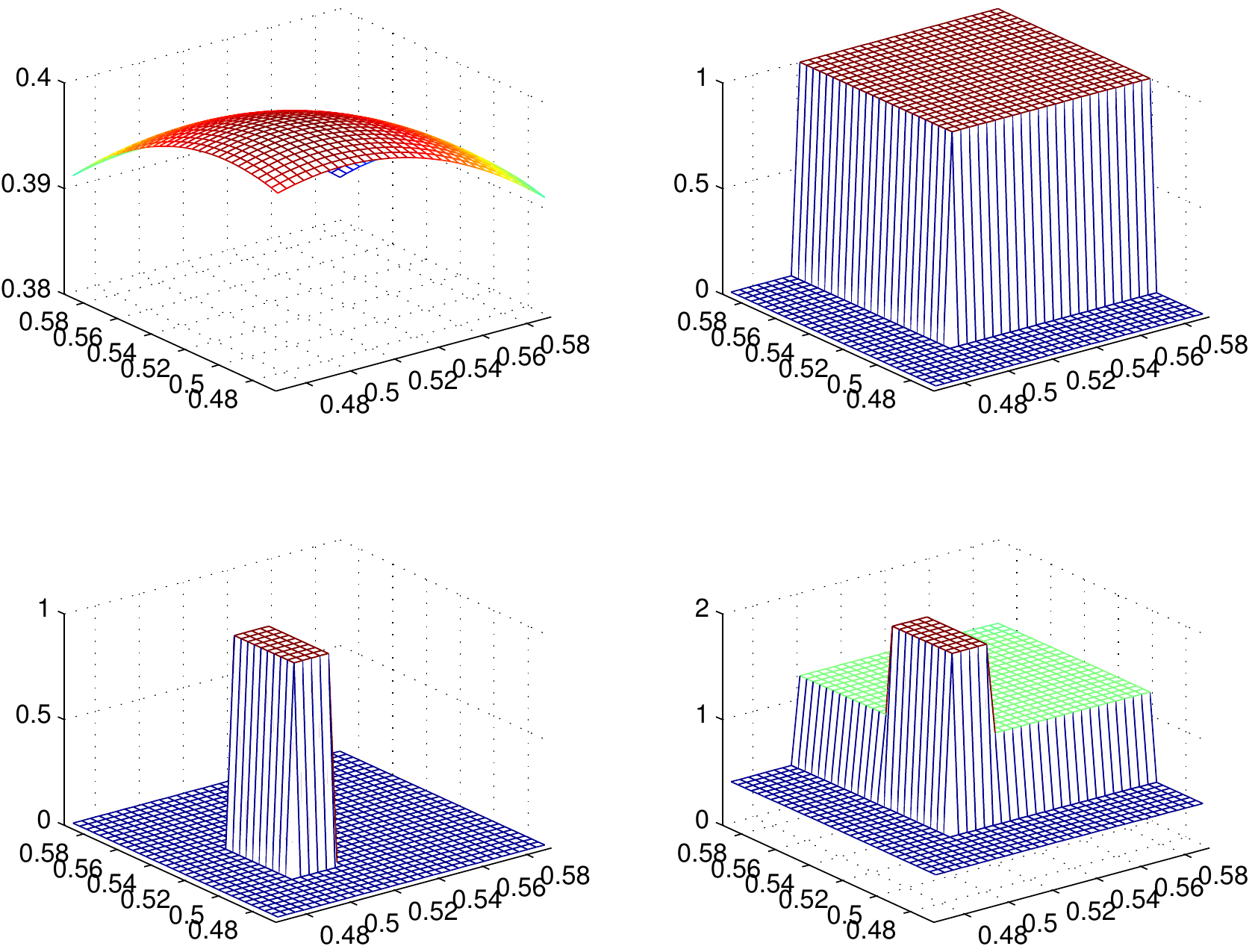}
\caption{Here we show the medium confined in patch $(9,9)$. The upper left corner presents the background media, and in the upper right and the lower left two corners we plot the two physical modes associated with two random variables. The lower right figure shows the sample media confined in this patch.}\label{fig:patch_99_media}
\end{figure}

In this example, we only show how to apply the random interpolation method and the reduced basis method on a local patch and study their performances. We will show the full process of the StoMsFEM on more realistic examples later. Let us pick patch $(9,9)$ for example, two random variables are present in this patch as shown in Figure~\ref{fig:patch_99_media}, and thus the local stiffness matrix and the local load vector are functions of only two random variables.

\begin{figure}
\centering
\includegraphics[width = 0.3\textwidth]{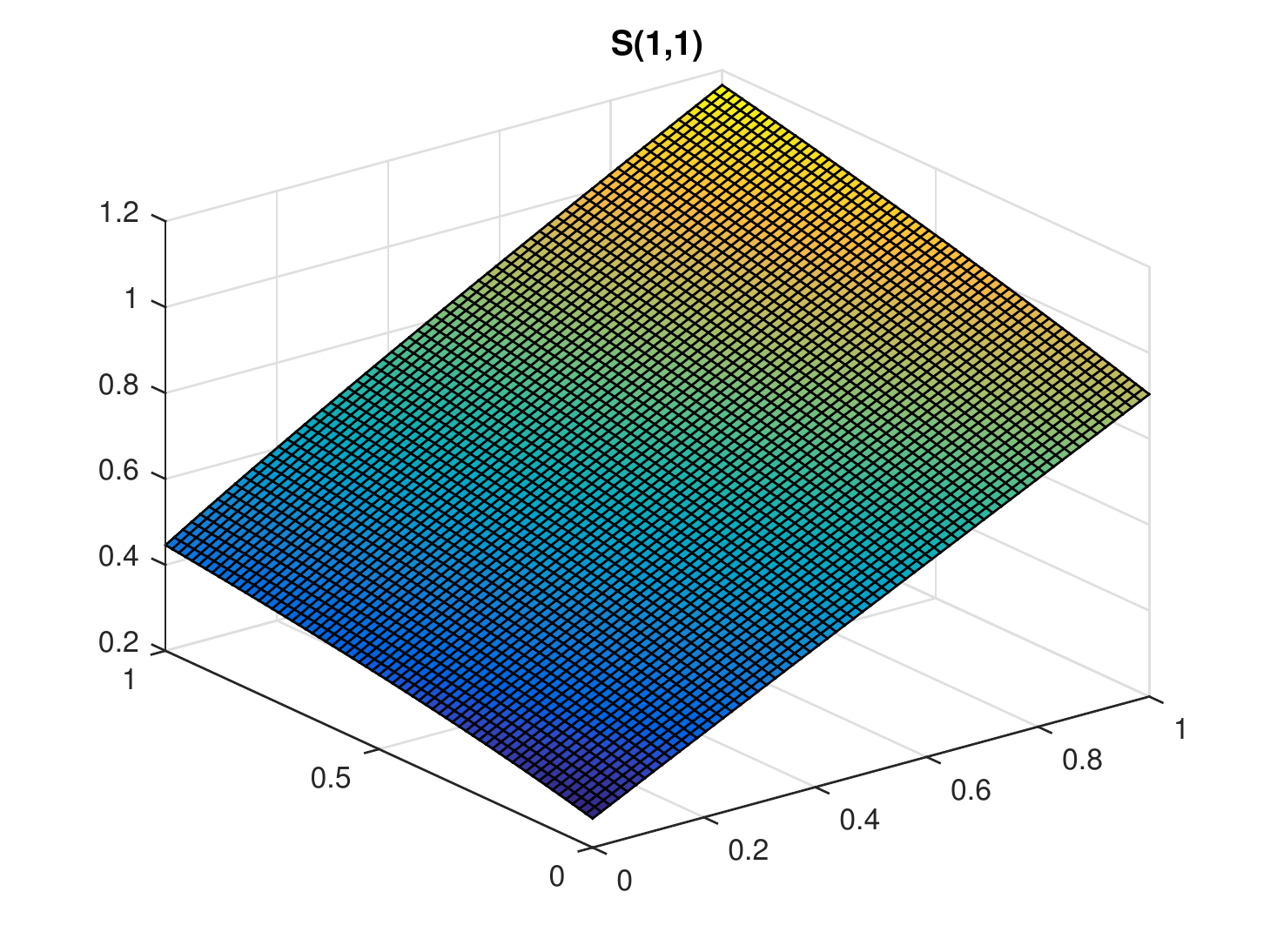}
\includegraphics[width = 0.3\textwidth]{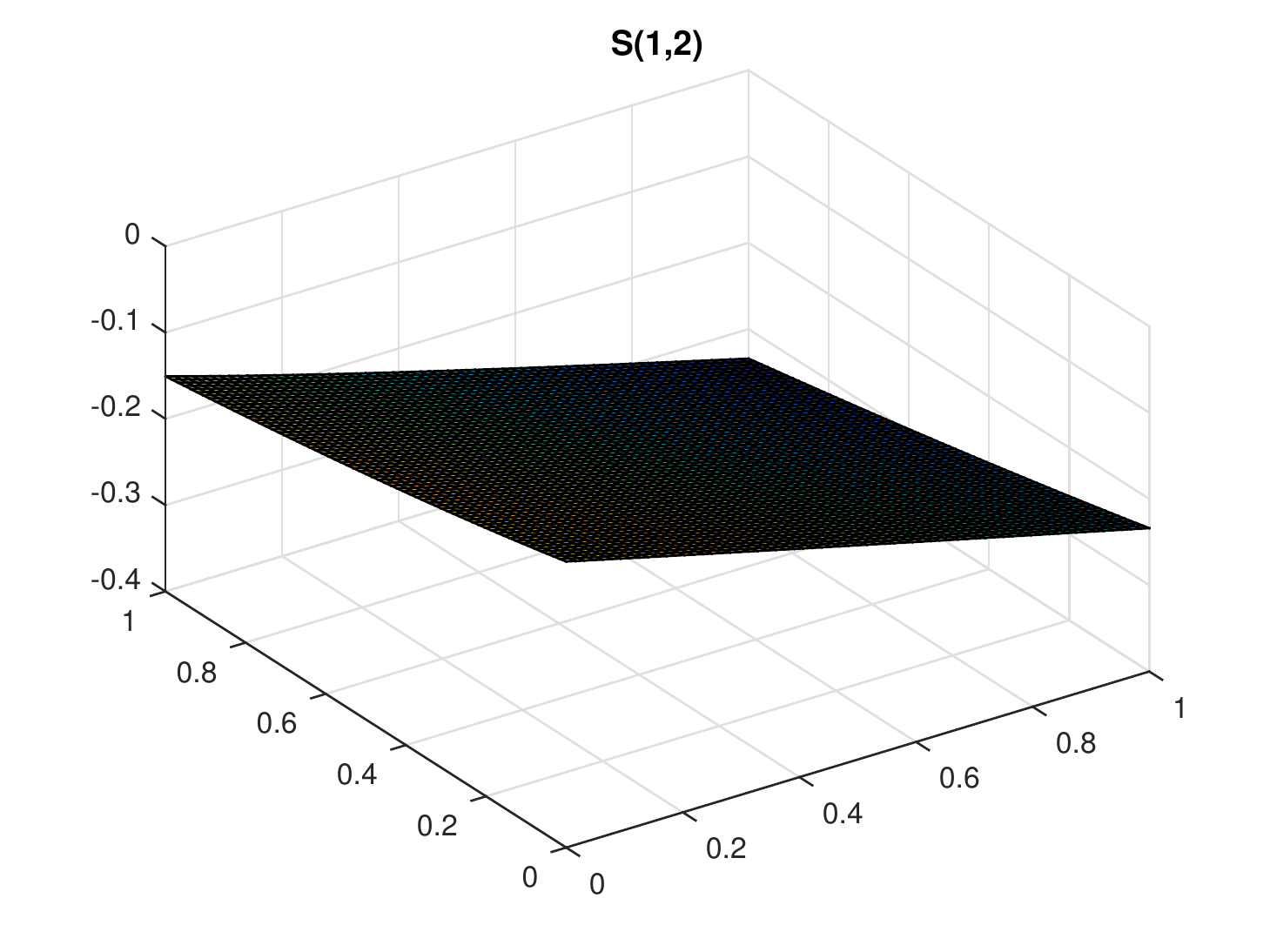}
\includegraphics[width = 0.3\textwidth]{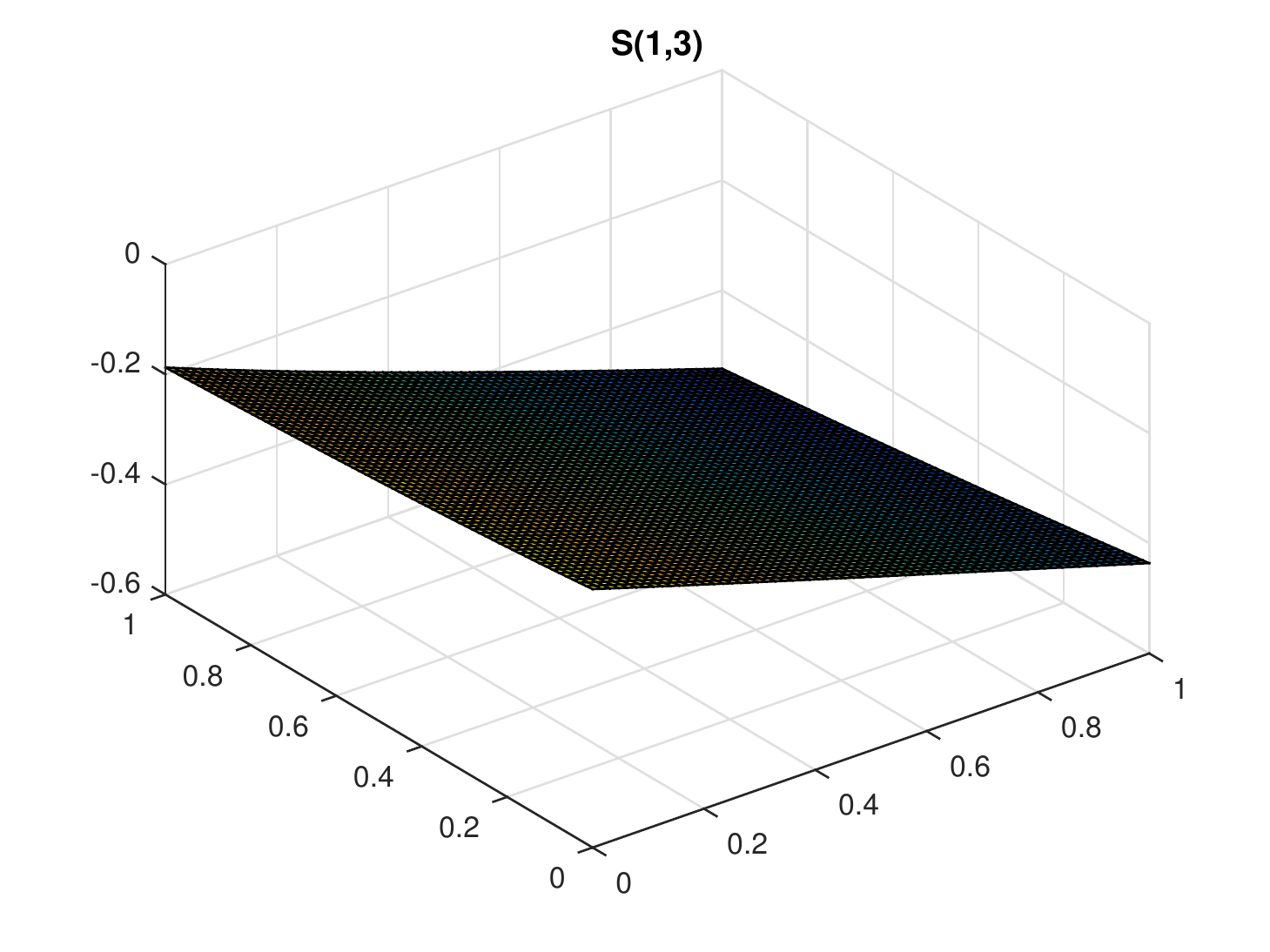}
\caption{The three figures demonstrate the dependence on the random variables of the $(1,1)$, $(1,2)$ and $(1,3)$ entries in the local stiffness matrix, when confined in patch $(9,9)$. Both random variables are uniformly distributed in $[0,1]$.}\label{fig:stiffness_99}
\end{figure}

\begin{figure}
\centering
\includegraphics[width = 0.3\textwidth]{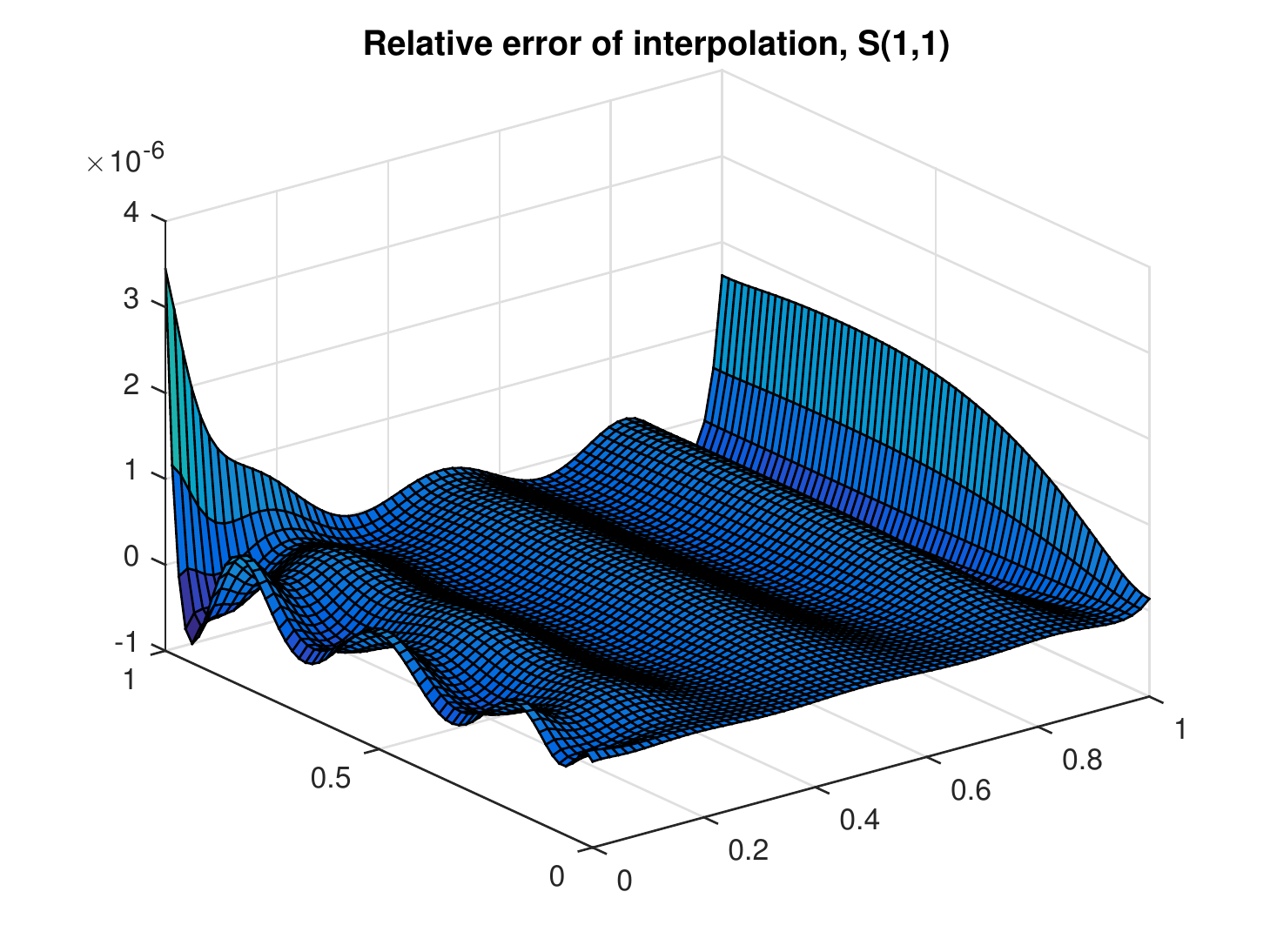}
\includegraphics[width = 0.3\textwidth]{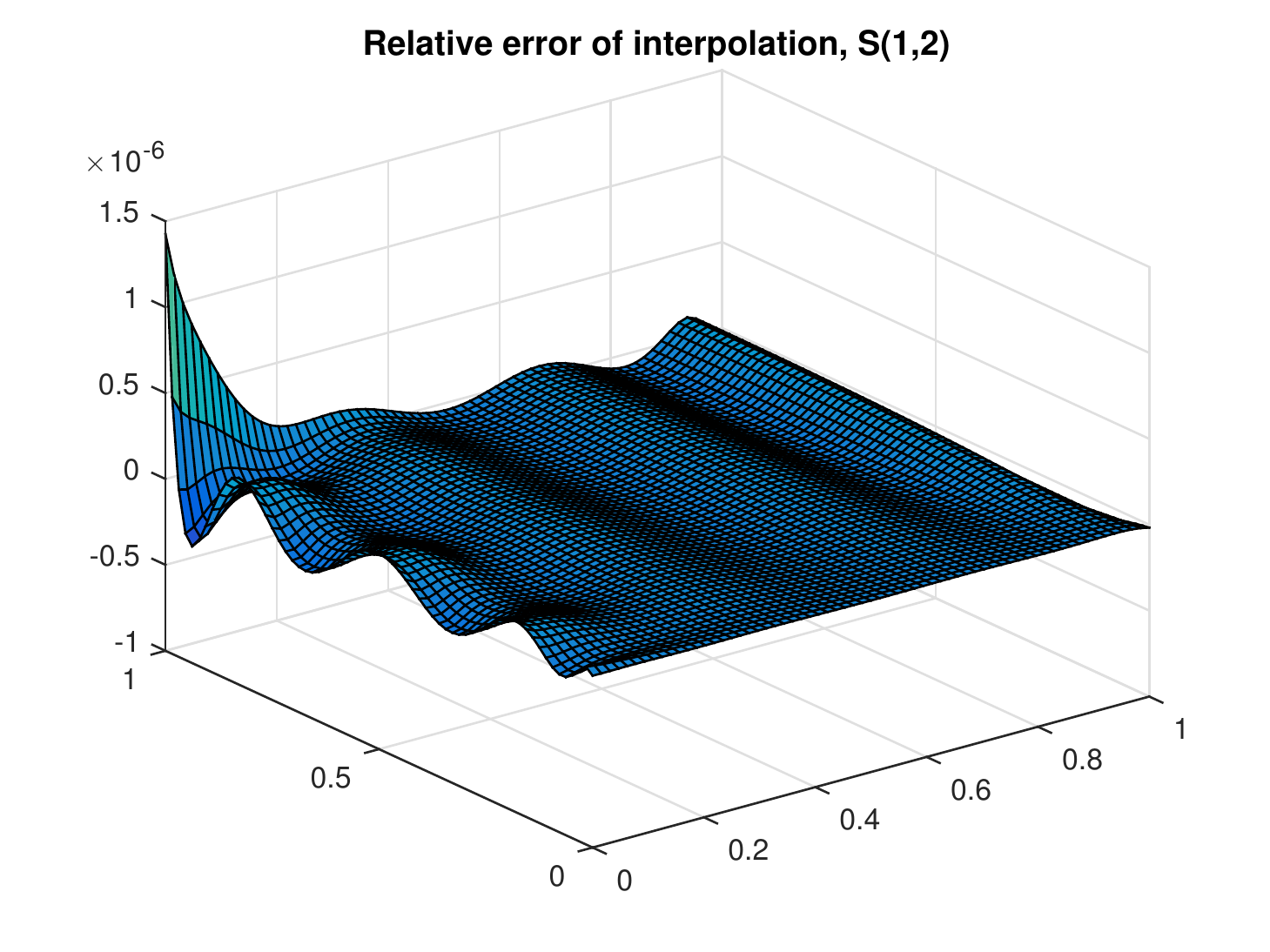}
\includegraphics[width = 0.3\textwidth]{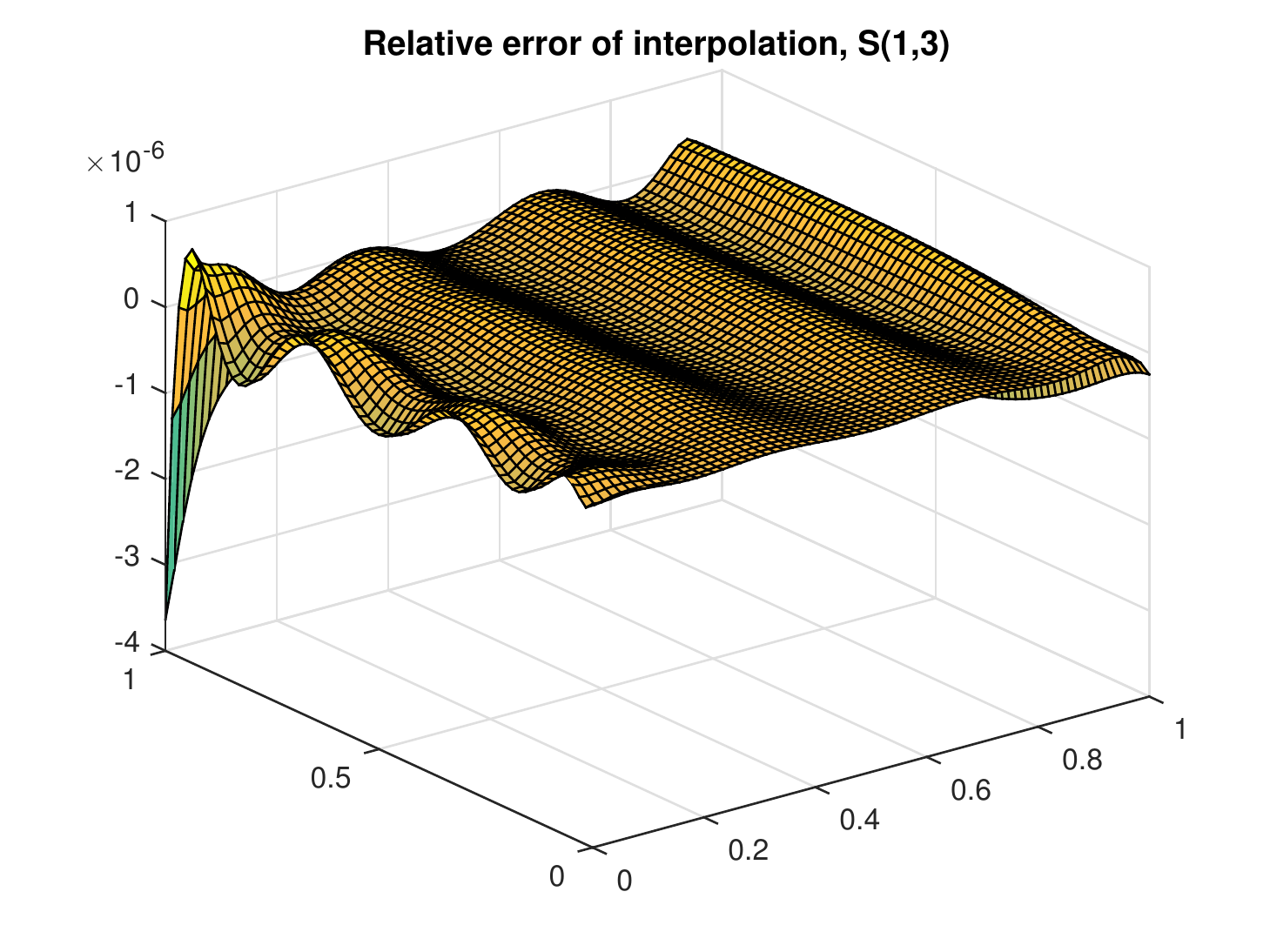}
\caption{The relative interpolation error for the $(1,1)$, $(1,2)$ and $(1,3)$ entries of the local upscaled stiffness matrix is in the order of $10^{-6}$. Its contribution of the final estimation error is negligible compared with the spatial discretization error and sampling error.}\label{fig:stiffness_99_interp}
\end{figure}

\begin{figure}
\centering
\includegraphics[width = 0.4\textwidth]{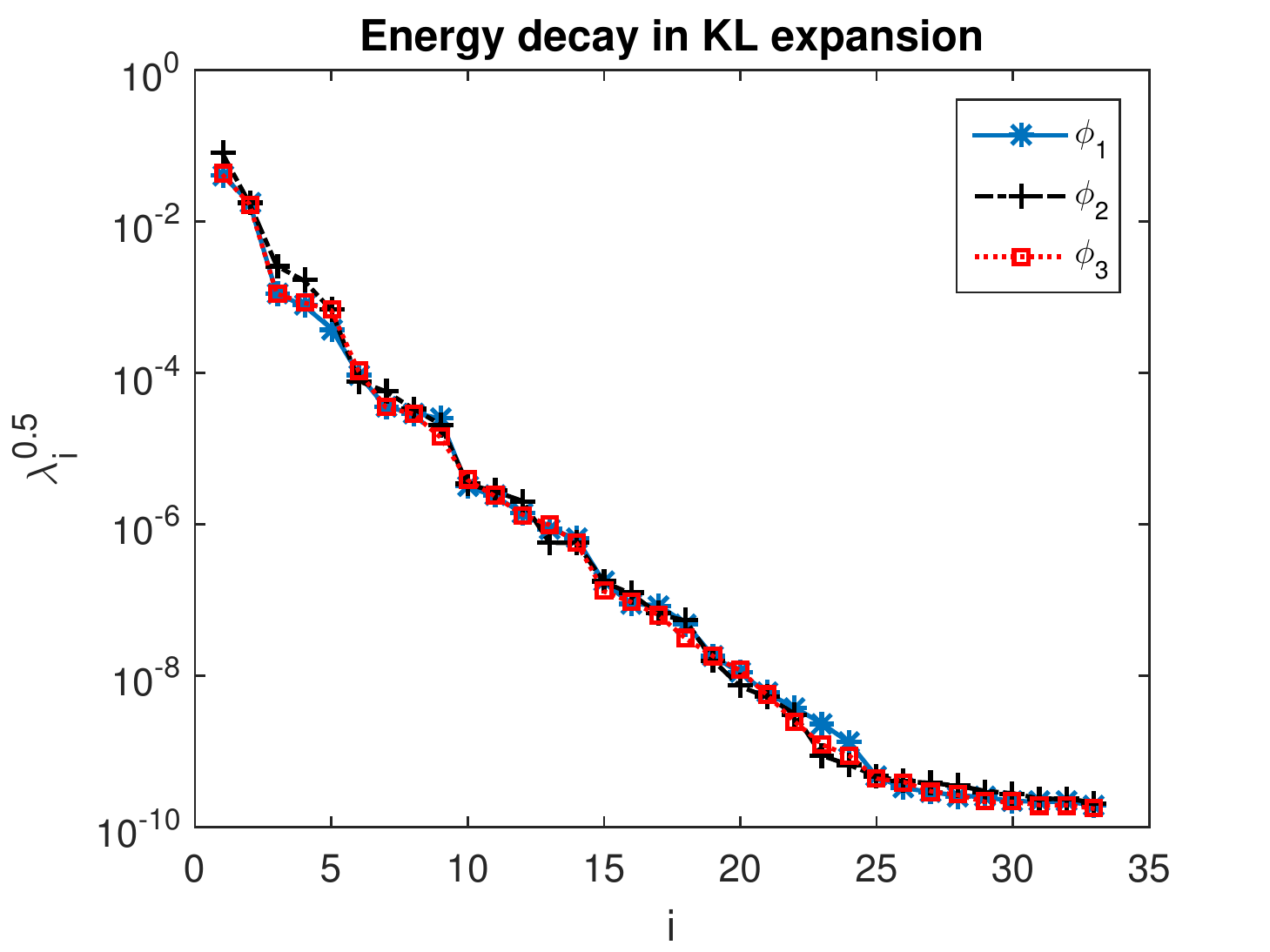}
\caption{The figure shows the fast KL energy decay for the solutions of the local cell problems~\eqref{eqn:MsFEB}, when confined to patch $(9,9)$. We truncate at $10^{-6}$ to obtain the reduced basis functions.}\label{fig:KLenergy_99_RB}
\end{figure}

\begin{figure}
\centering
\includegraphics[width = 0.3\textwidth]{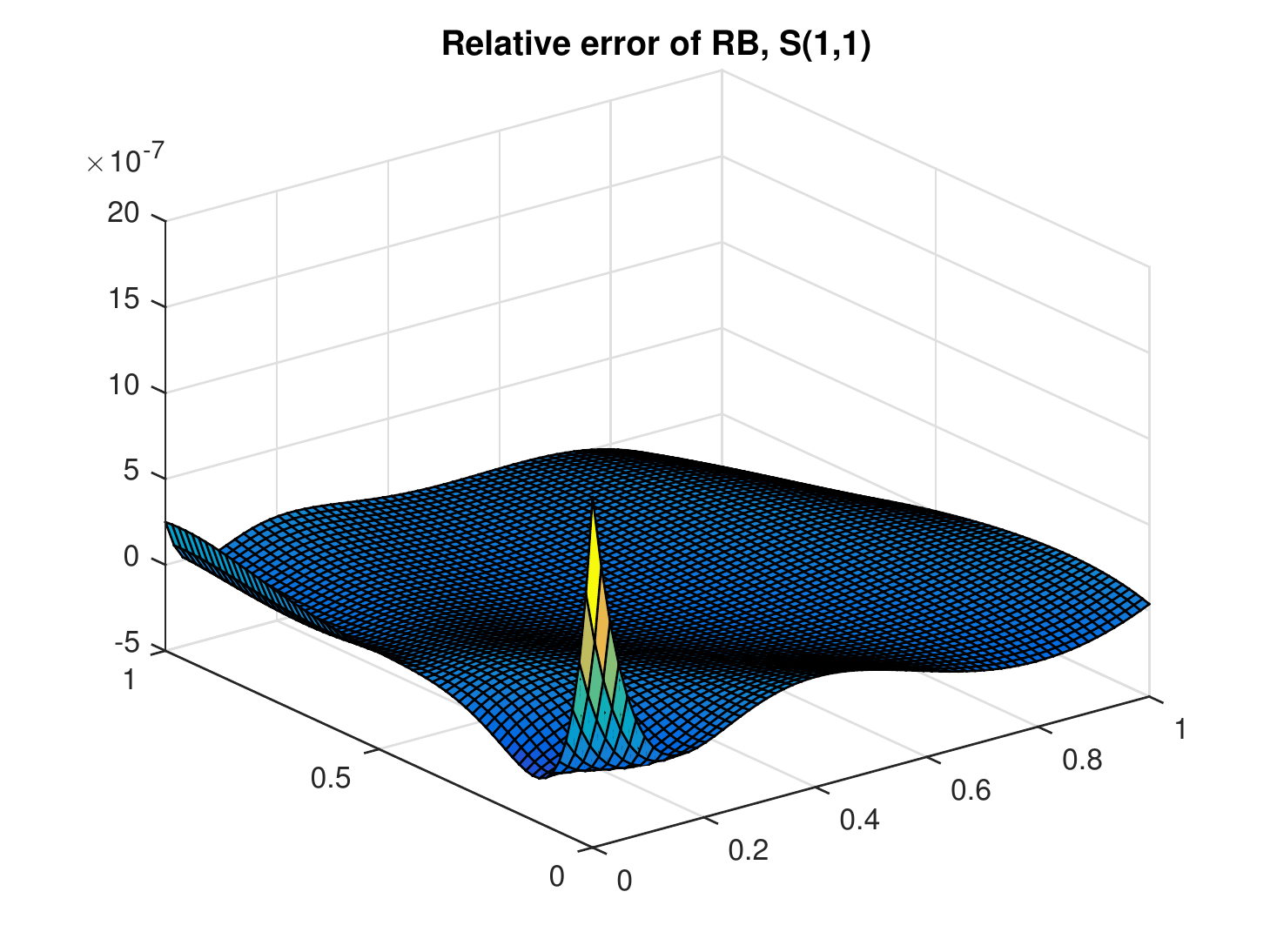}
\includegraphics[width = 0.3\textwidth]{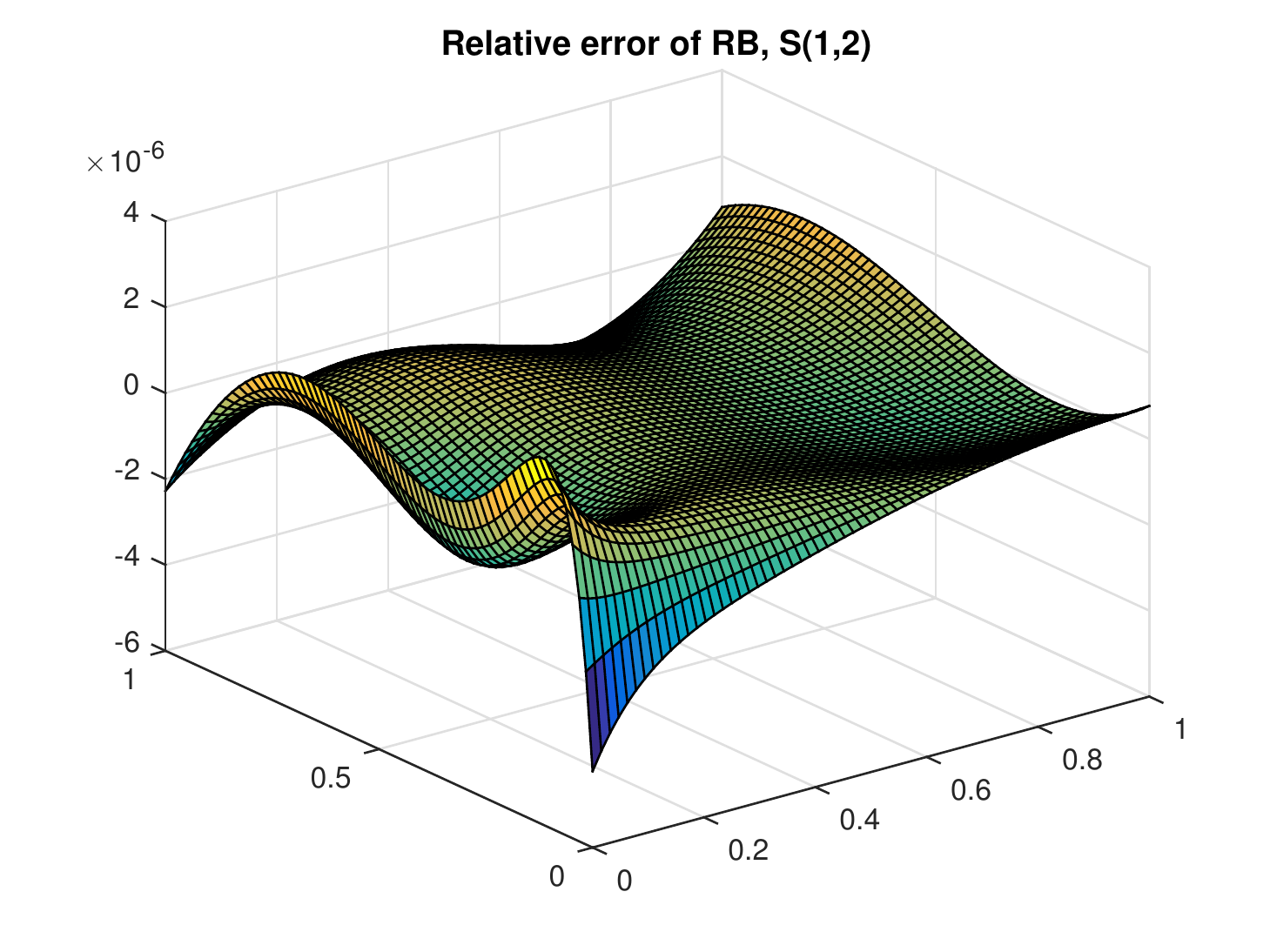}
\includegraphics[width = 0.3\textwidth]{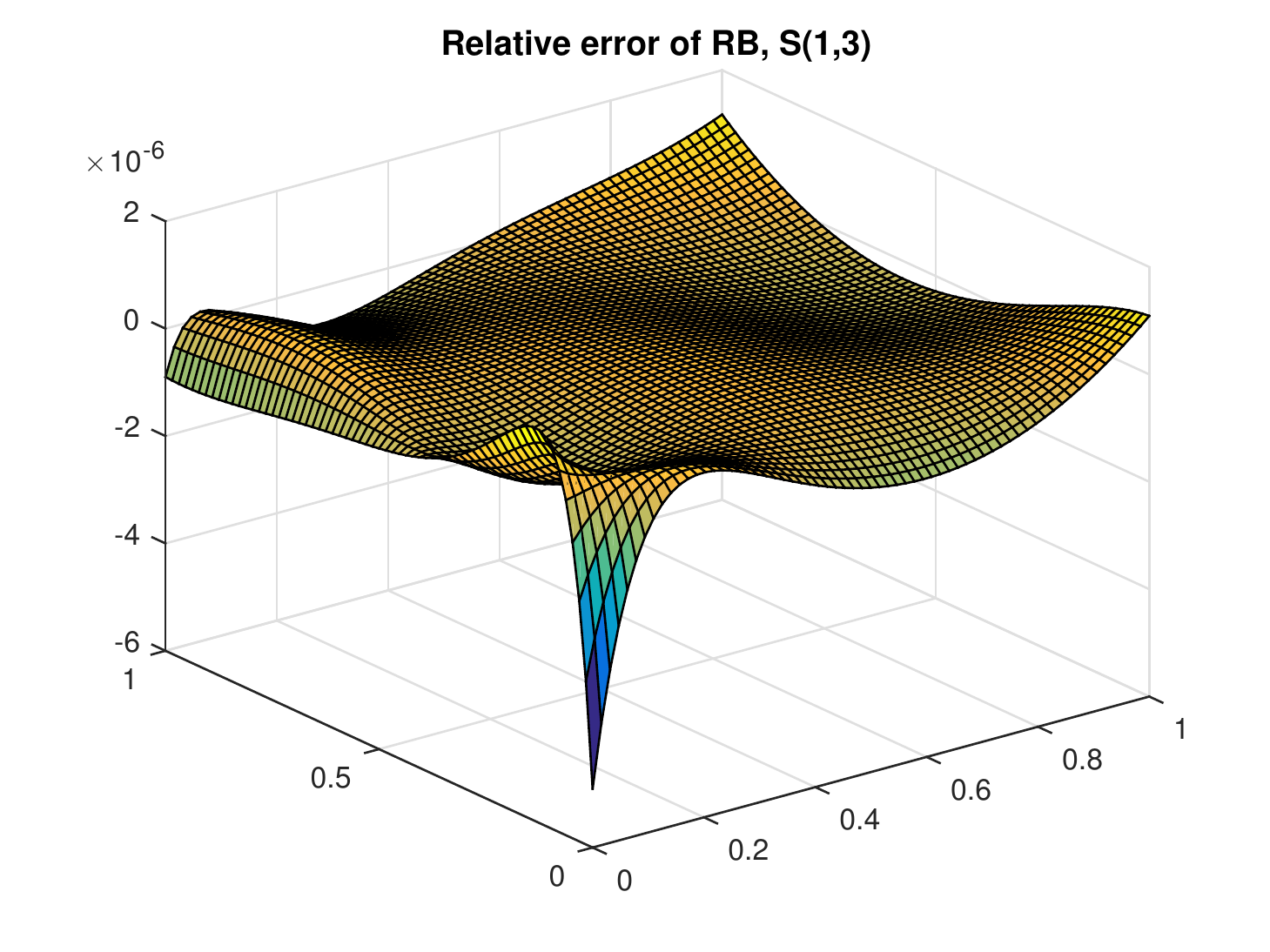}
\caption{The relative error for the $(1,1)$, $(1,2)$ and $(1,3)$ entries of the local upscaled stiffness matrix is in the order of $10^{-6}$.}\label{fig:stiffness_99_RB}
\end{figure}

In the offline step, we use both the random interpolation method and the reduced basis method to construct approximations for the upscaled local stiffness matrices. In Figure~\ref{fig:stiffness_99} we plot the $(1,1)$, $(1,2)$ and $(1,3)$ entries of the stiffness matrix's dependence on the two random variables.
For the random interpolation method, we take $9$ Chebyshev nodes along each dimension and take their tensor products. In Figure~\ref{fig:stiffness_99_interp} we plot the relative interpolation error for the $(1,1)$, $(1,2)$ and $(1,3)$ entries. The relative interpolation error is in the order of $10^{-6}$, whose contribution to the final estimation error is negligible compared with the spatial discretization error and sampling error. 
For the reduced basis method, we perform the KL expansion of the three basis functions, obtain their reduced basis functions, and precompute the relevant quantities. In Figure~\ref{fig:KLenergy_99_RB}, we show the fast energy decay in the KL expansion of these basis functions. We truncate the KL expansion at $\sqrt{\lambda_Q/\lambda_1} < 10^{-6}$, resulting in 15, 13 and 15 basis functions for $\phi^1, \phi^2$ and $\phi^3$ respectively. In Figure~\ref{fig:stiffness_99_RB} we see that the relative error to compute the $(1,1)$, $(1,2)$ and $(1,3)$ entries is also in the order of $10^{-6}$. 

It is worth mentioning that the oversampling domain is four times bigger than the original patch. If the effective region for a random variable falls in the boundary layer region, the associated random variable shows limited impact on the stiffness matrix. This anisotropic property suggests that we can do interpolation on dimension-adaptive grids to reduce the number of interpolation nodes when the local dimension grows. We will use dimensional-adaptive grids in our next two examples. We point out that the local reduced basis approach automatically detects this anisotropic property and always gives the most important basis functions for the local cell problems~\eqref{eqn:MsFEB}.

\subsection{An example with high contrast random medium}\label{exam:highcontrast}
\begin{figure}[htp]
\centering
     \includegraphics[width=0.3\textwidth]{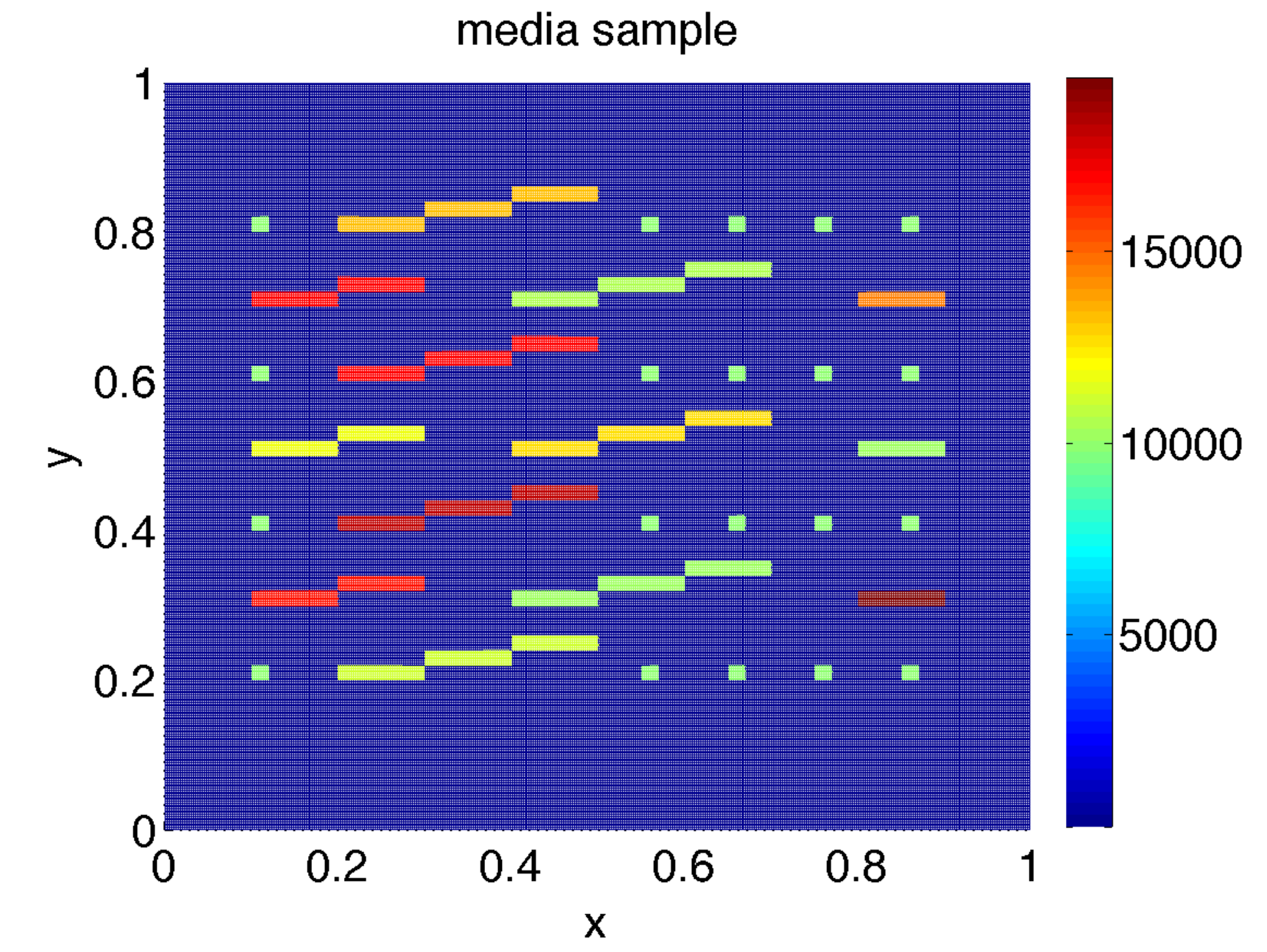}
     \includegraphics[width=0.3\textwidth]{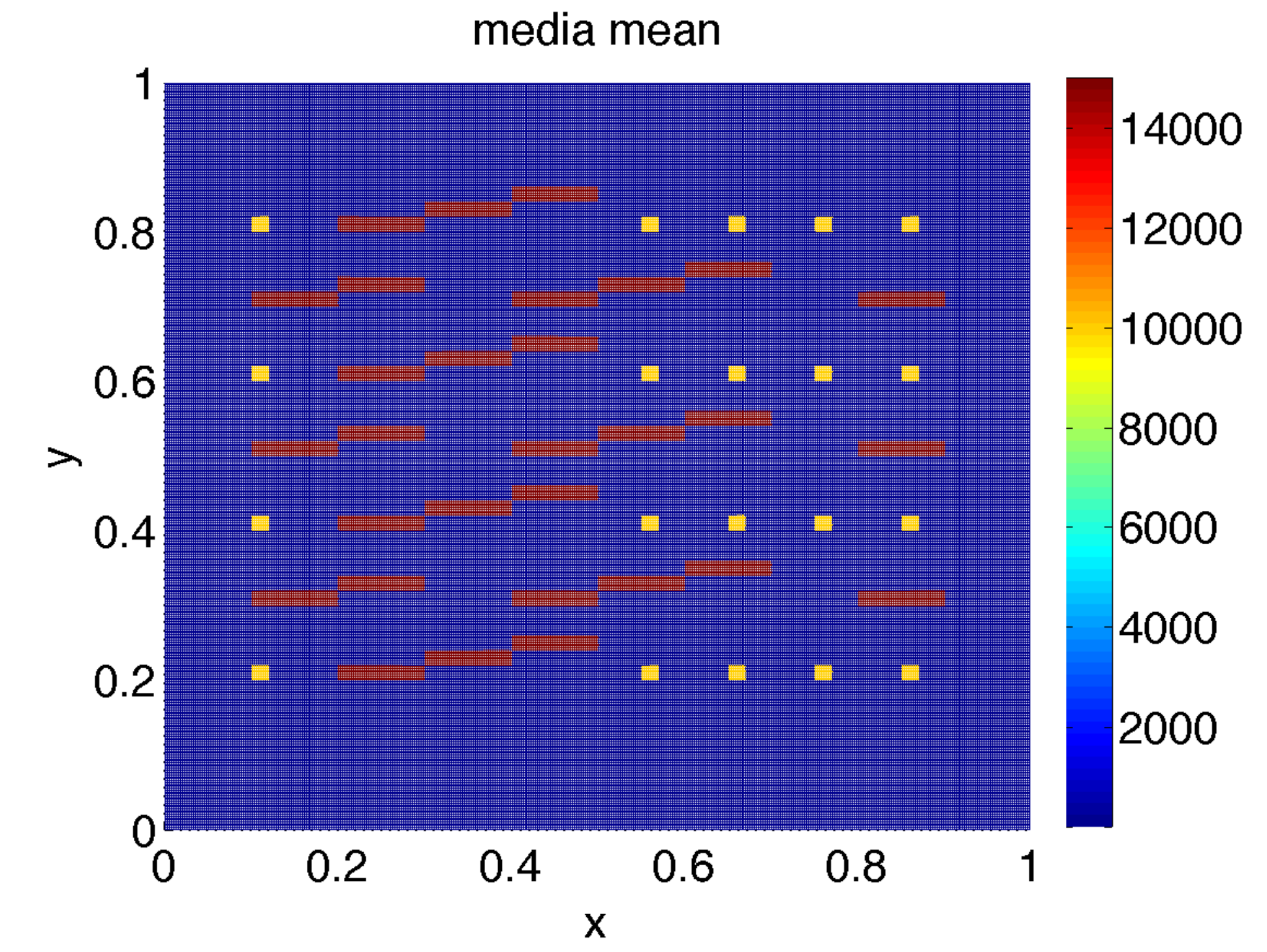}
     \includegraphics[width=0.3\textwidth]{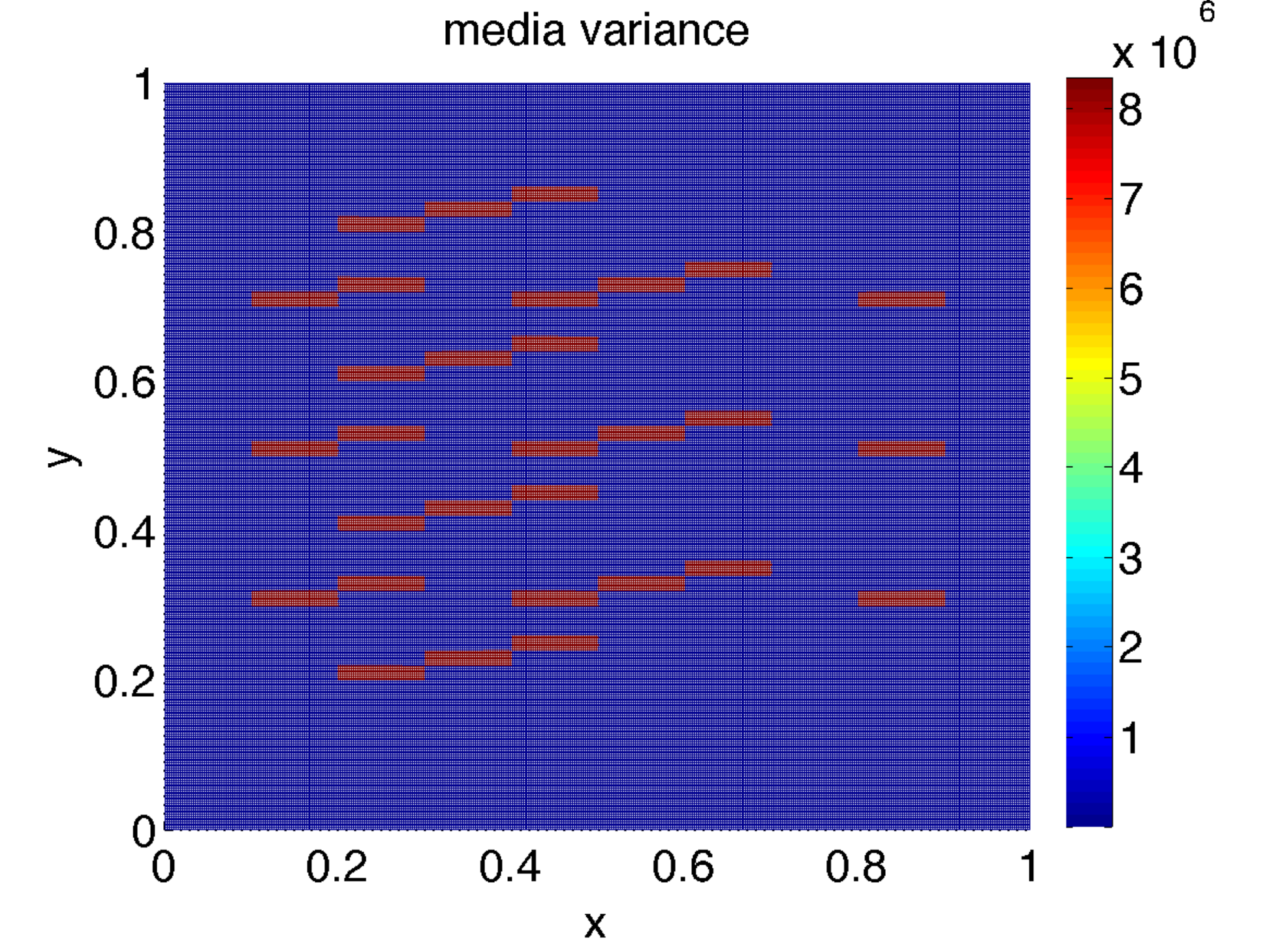}
\caption{Left: one sample medium; middle: medium mean; right: medium variance. There are 13 high permeability (of order $10^{4}$) channels in the x direction and a few high permeability inclusions. The background permeability is of order 1.}\label{fig:HighContrastMedia}
\end{figure}

The random medium of this example contains a non-constant global background, channels with high permeability and localized inclusions. One sample and statistical properties of the random media are shown in Figure~\ref{fig:HighContrastMedia}. We can see that there are several high permeability channels in the x-direction and some high permeability inclusions. Utilizing the ISMD presented in Section~\ref{sec:localrepresentation}, we parametrize this random medium as Eqn.~\eqref{eq:highcontrastmedia}.
\begin{equation}\label{eq:highcontrastmedia}
\kappa(x,y,\omega) = f_0(x,y) + \xi_0(\omega) + \sum_{k=1}^{13} f_k(x,y)\xi_k(\omega).
\end{equation}
Here, $\xi_0$ is the global random variable uniformly distributed in $[0,1]$ corresponding to the low permeability background, and $\{\xi_k\}_{k=1}^{13}$ are independent random variables uniformly distributed in $[10^4,2\times 10^4]$ corresponding to the high permeability channels. In this problem, we incorporate the StoMsFEM with both the Monte Carlo method and the sparse grid SC method. This demonstrates that our StoMsFEM can be easily combined with most non-intrusive global stochastic methods.

We use the MsFEM in~\cite{hou_multiscale_1997}, with oversampling and linear boundary conditions. We point out that the MsFEM is not the best local upscaling method for high contrast coefficients, and it is not guaranteed to have error estimate $\|u_H - u\|_{H_1}$ small. Local upscaling methods specifically designed for high contrast problems can be found at~\cite{chu2010new, efendiev2011multiscale, owhadi_localized_2011}. In this paper, we focus on the accuracy of the proposed random interpolation and the reduced basis method, i.e., $\hat{u}_H - u_H$, instead of the accuracy of the upscaling method, i.e., $u_H - u$.

In the physical domain $[0,1]^2$, we have a uniform coarse mesh $\CalT_H$ with mesh size $H_x = H_y = 0.05$ and a fine mesh $\CalT_h$ with $h_x = h_y =0.0025$. Due to the high contrast permeability, we take a relatively large oversampling ratio $\eta = 3$. Thanks to ISMD, the local stochastic dimensions of the parametrization~\eqref{eq:highcontrastmedia} are small, typically 2 or 3, on these oversampling local coarse grid elements. As in Example~\ref{sec:patchstudy}, we use both the random interpolation method and the reduced basis method to construct approximations for local upscaled stiffness matrices. Due to the locally low dimensionality, we achieve negligible errors when approximating the local upscaled matrix $\vct{S}^m$ with a small number of interpolation nodes or reduced basis functions. For example, when we approximate $\vct{S}_{1,1}^m$ on patch (14,9), the random interpolation method achieves $\Or(10^{-6})$ relative error with Chebyshev interpolation on a $5\times 16$ Chebyshev grid, while the reduced basis method also achieves $\Or(10^{-6})$ relative error with only 7 reduced basis functions for all $\phi^{ml}$'s. The error plots look similar to Figure~\ref{fig:stiffness_99_interp} and Figure~\ref{fig:stiffness_99_RB} in our patch-study example, and we do not show them here any more.

In the global Monte Carlo solver, since the basis functions constructed from oversampling are nonconforming, we apply the Petrov-Galerkin MsFEM formulation~\cite{hou_PetrovGalerkin_2004} with the standard bilinear basis on the coarse mesh as test functions.

\begin{figure}[htp]
 \centering
     \includegraphics[width=0.3\textwidth]{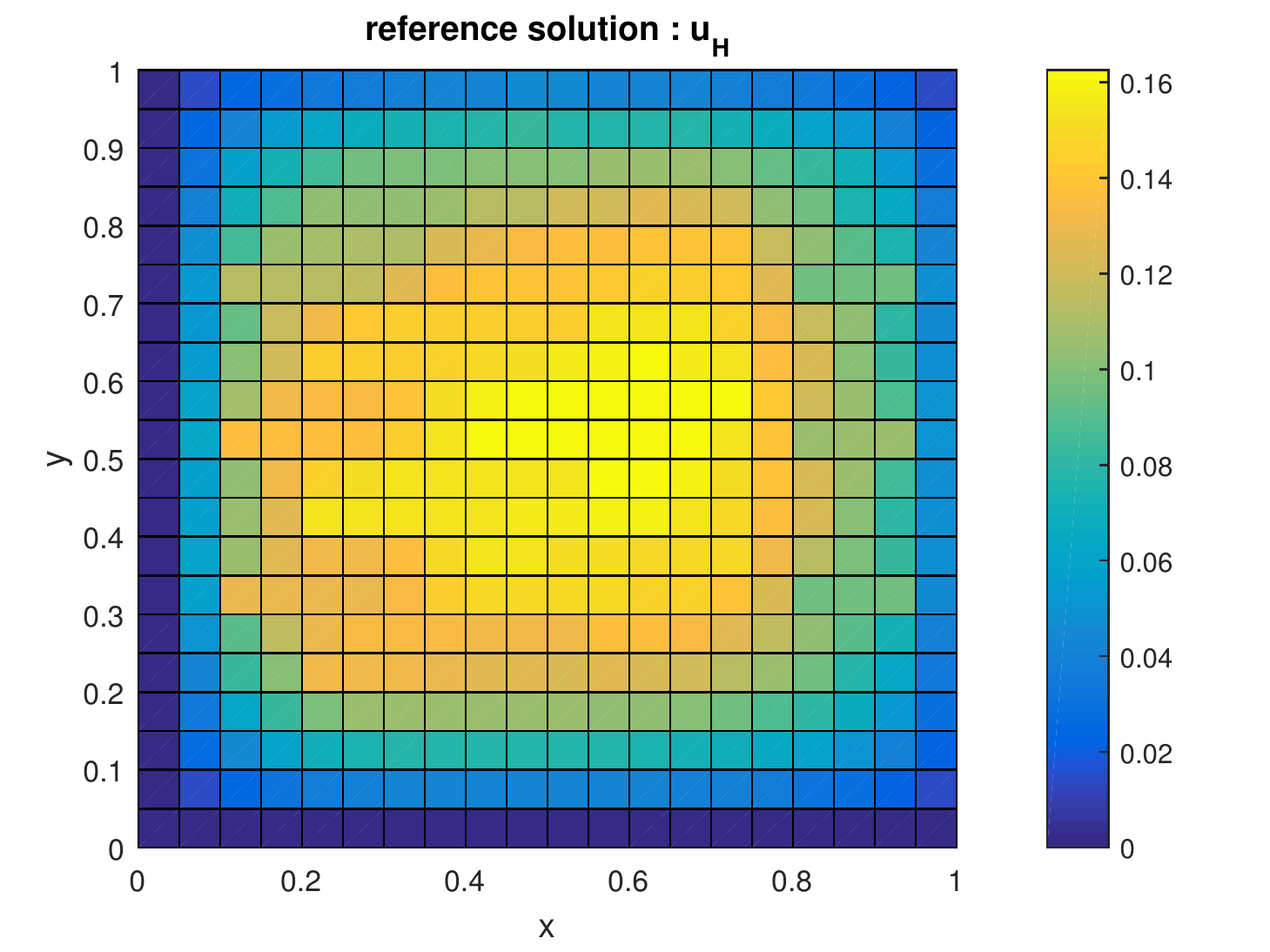}
      \includegraphics[width=0.3\textwidth]{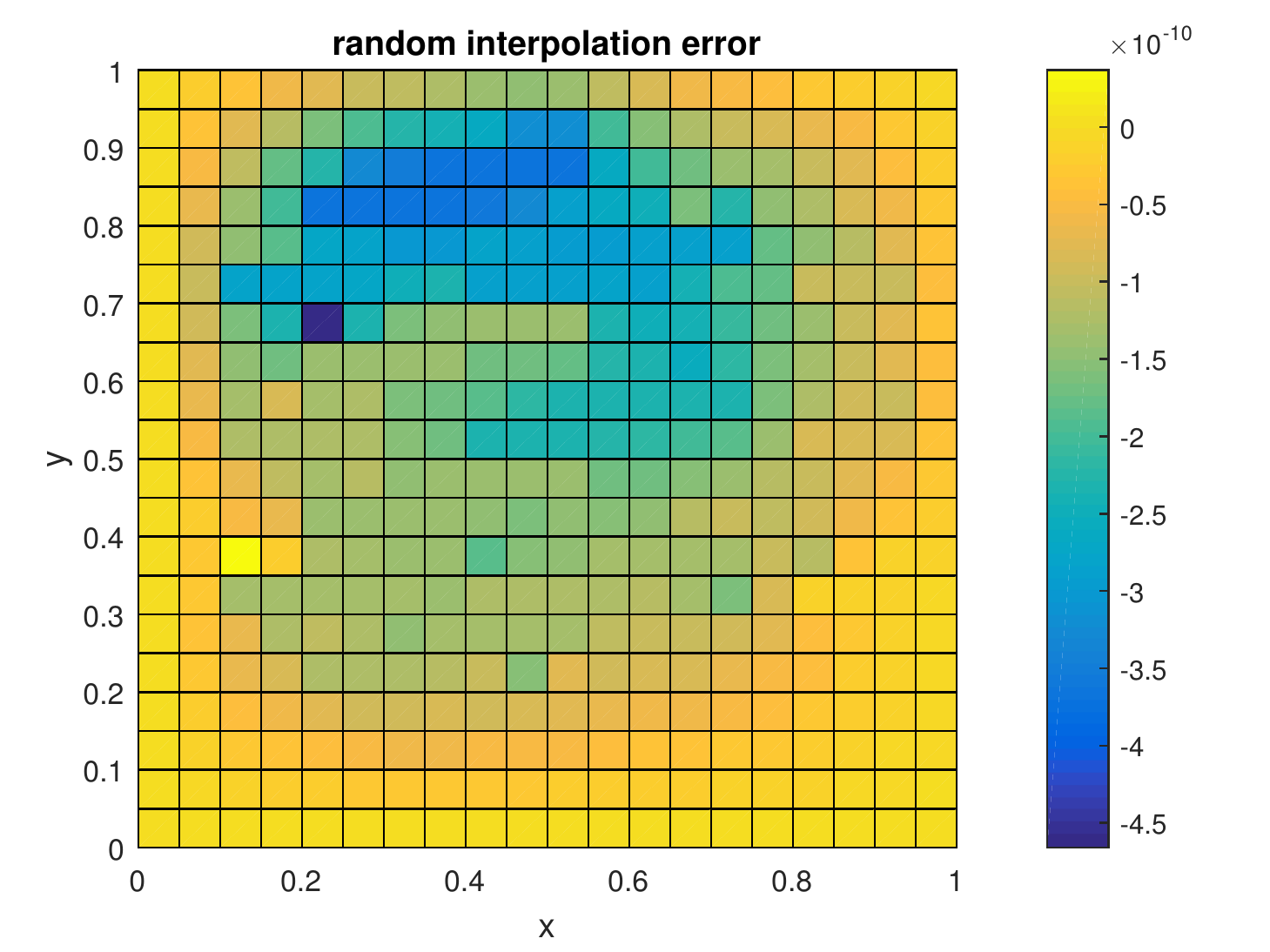}
      \includegraphics[width=0.3\textwidth]{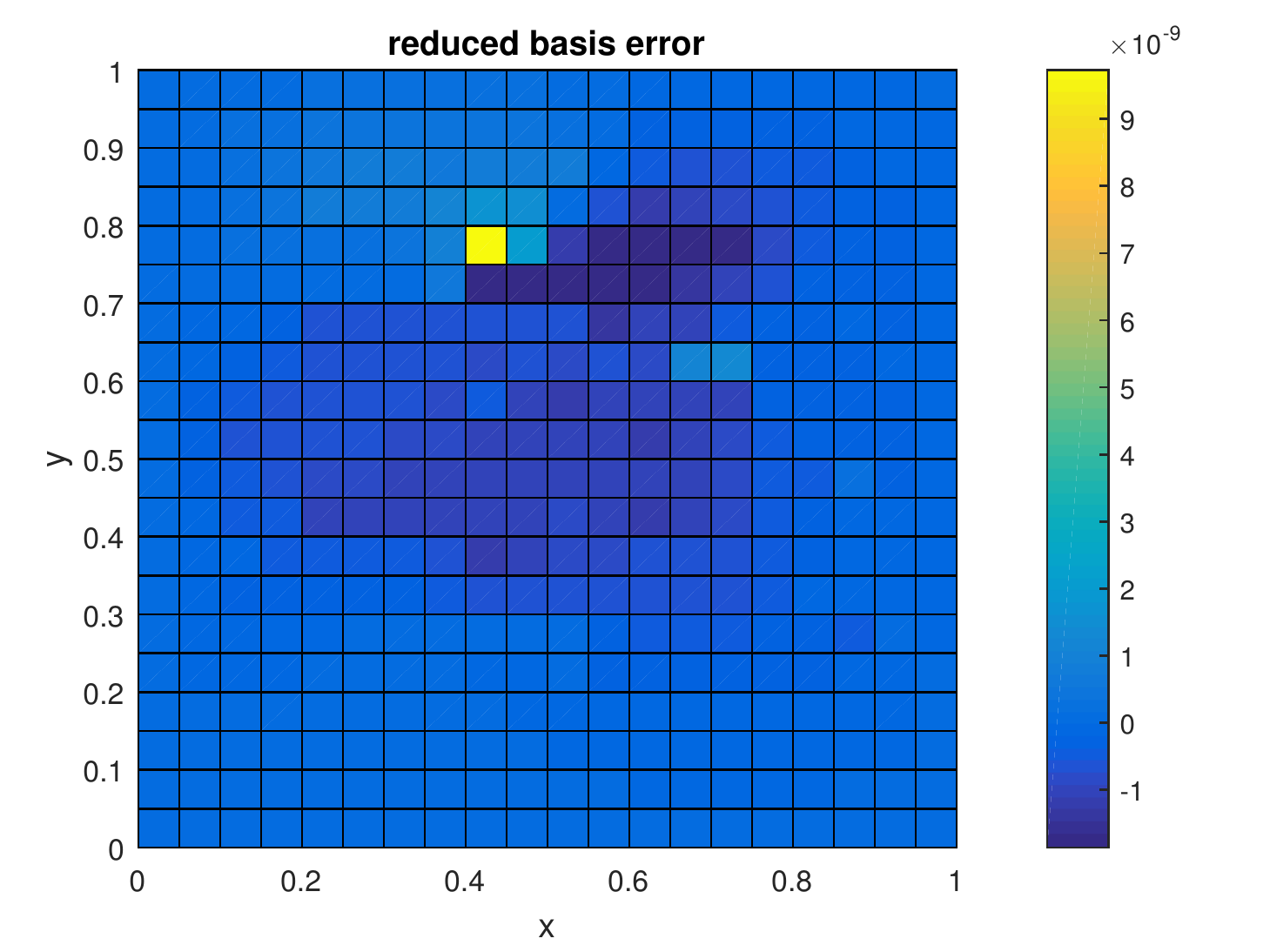}
\caption{Left: one sample solution from direct MsFEM; middle: the absolute error of the approximation by the random interpolation, which is of order $10^{-10}$; right: the absolute error of the approximation by the reduced basis method, which is of order $10^{-9}$.}\label{fig:HC_sample}
\end{figure}

In the first experiment, we set the source $b(x)=1$ and a zero Dirichlet boundary condition. In Figure~\ref{fig:HC_sample} we show one sample solution directly computed by the MsFEM, i.e. $u_H$, and the absolute error of the StoMsFEM approximating solutions $\hat{u}_H$. We can see that $\hat{u}_H - u_H$ is of the order $10^{-10}$ for the random interpolation method and $10^{-9}$ for the reduced basis method, which is negligible compared with the spatial discretization error $u_H - u$. Therefore, we can treat the approximating solution $\hat{u}_H$ as the solution $u_H$ computed directly by the MsFEM. However, their computational times are very different. Table~\ref{table:lPCgMC_cost} shows CPU times and the actual computational cost ratio. In our setting, the grid size ratio $\eta H / h = 60$, and on average we use $N_c \approx 50$ local interpolation points. Theoretically, we have $1/R = \Or((\eta H / h)^{\gamma}/N_c) = \Or(72)$ for the random interpolation method. The saving we observe is purely in the order of $\Or((\eta H / h)^{\gamma})$ because the cost to evaluate interpolants is negligible in practice. For the reduced basis method with $K_m = 2$, we have $1/R = \Or((\eta H / h)^{\gamma}/K_m^3) = \Or(450)$ theotrically, which matches what we observed numerically.  
\begin{table}[ht]
\caption{Computational cost for one sample solution(unit: s)}\label{table:lPCgMC_cost}
\centering
\begin{tabular}{c|c|c}
\hline
naive MsFEM       &  StoMsFEM(random interpolation) &		StoMsFEM(reduced basis) \\
\hline
27.38 		& 		0.0133 (1/R = 2060)		 &  		0.0814 (1/R = 336)	\\
\hline
\end{tabular}
\end{table}

To balance the spatial discretization error and sampling error as discussed in Section~\ref{sec:globalerror}, we need about $N_{\text{on}} = \Or(H^{-4}) = \Or(10^5)$ for Monte Carlo sampling. Therefore, we estimate the mean and standard deviation of $\hat{u}_H$ by both the random interpolation method and the reduced basis method on the same set of $10^5$ independent samples. The difference between these two methods is of the order $10^{-9}$, confirming again that the error introduced by the random interpolation and reduced-basis method is negligible. To compute these $10^5$ samples, it takes 1329 seconds for the random interpolation approximation and 8146 seconds for the reduced basis method. If we directly compute these samples by the MsFEM, it would take $2.7 \times 10^6$ seconds.

In the second experiment, we combine our StoMsFEM with the sparse grid SC method. Since variables $\{\xi_k\}_{k=0}^{13}$ are independent and $u_H$ is smooth with respect to $\{\xi_k\}_{k=0}^{13}$, we can implement the sparse grid SC method to estimate $\EE[u_H]$ and $\text{var}[u_H]$. The dimension-adaptive sparse grid integration is performed on the sparse grid toolbox~\cite{spalg, spdoc} and the dimension-adaptive degree is set to be $0.6$. As described in Section~\ref{sec:globalSC}, we prepare the local upscaled stiffness matrices at the local sparse grid collocation nodes in the offline stage. Since the local dimensions are small, the biggest number of local collocation nodes is only 1073. It takes about $6700$ seconds to finish the offline computation. In the online stage, each sample takes only $9.8\times 10^{-3}$ second because we only look up the precomputed dictionaries to get the local upscaled quantities. The numerical cost ratio $1/R = 2800$, which exactly matches our theoretical estimation $1/R = \Or((\eta H / h)^{\gamma}) = \Or(3600)$ for global (sparse grid) SC method.

\begin{figure}[htp]
   \begin{center}
     \includegraphics[width=0.5\textwidth]{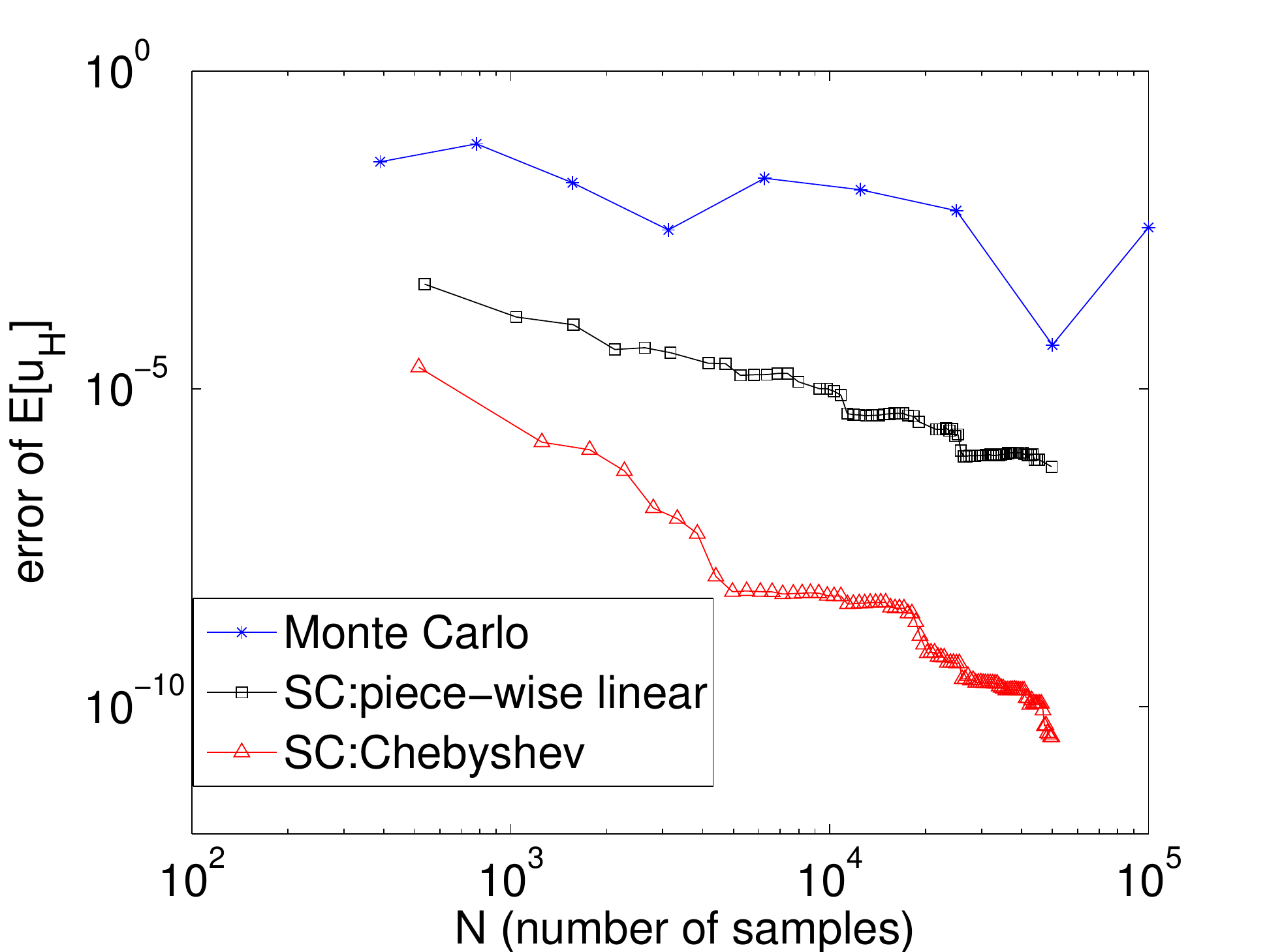}
   \end{center}
\caption{Comparison between the Monte Carlo method, the SC with the sparse grid trapezoidal rule (SC:piece-wise linear) and the sparse Clenshaw-Curtis formulas (SC:Chebyshev). It shows that sparse grid collocation with high order quadrature rules is superior in this example.}\label{fig:errorcompare}
\end{figure}

We implement the sparse grid integration using both the trapezoidal rule and the Clenshaw-Curtis formulas~\cite{gerstner1998numerical}. Taking the sparse grid integration with the Clenshaw-Curtis formulas with 50433 collocation nodes as the reference $\EE[u_H]$, we define the quadrature estimation error as:
\begin{equation}
e_{\text{quad}}\left( \mathcal{I}[u_H] \right) := \|\mathcal{I}[u_H] - \EE[u_H]\|_2.
\end{equation}
In Figure~\ref{fig:errorcompare}, we compare the performance of the MC method, the SC method with the sparse grid trapezoidal rule (piece-wise linear) and the sparse Clenshaw-Curtis formulas (Chebyshev). Due to the smoothness of $u_H$ with respect to $\vct{\xi}$, the SC with the sparse Clenshaw-Curtis formulas has the best convergence rate ($\approx 2.64$ from linear regression), the SC with the sparse grid trapezoidal rule has convergence rate about 1.5 and that of Monte Carlo method is only about 0.65. Due to the stochastic nature of the MC estimator, we can see that its estimation error oscillates while slowly decreasing. 

In Table~\ref{table:spcollocation}, we list different parts of CPU times for the most accurate SC with the sparse Clenshaw-Curtis formulas in Figure~\ref{fig:errorcompare}, which has about 49805 collocation points.
\begin{table}[ht]
\caption{Computation cost for Stochastic Collocation on sparse grid, unit: s}\label{table:spcollocation}
\centering
\begin{tabular}{c|c|c}
\hline
offline       &   online  & 	online per sample \\
\hline
6700 		& 		488		 &  		0.0098 (1/R = 2800)	\\
\hline
\end{tabular}
\end{table}
Since the standard FEM on the fine grid takes 27.38 seconds per sample, our SC method based on the sparse representation will have computational saving as long as the total collocation points $N_{\text{on}}\ge \frac{N_{\text{off}}}{1-R} = 246$, which is obviously true in our case.

\begin{figure}[htp]
   \begin{center}
     \includegraphics[width=0.35\textwidth]{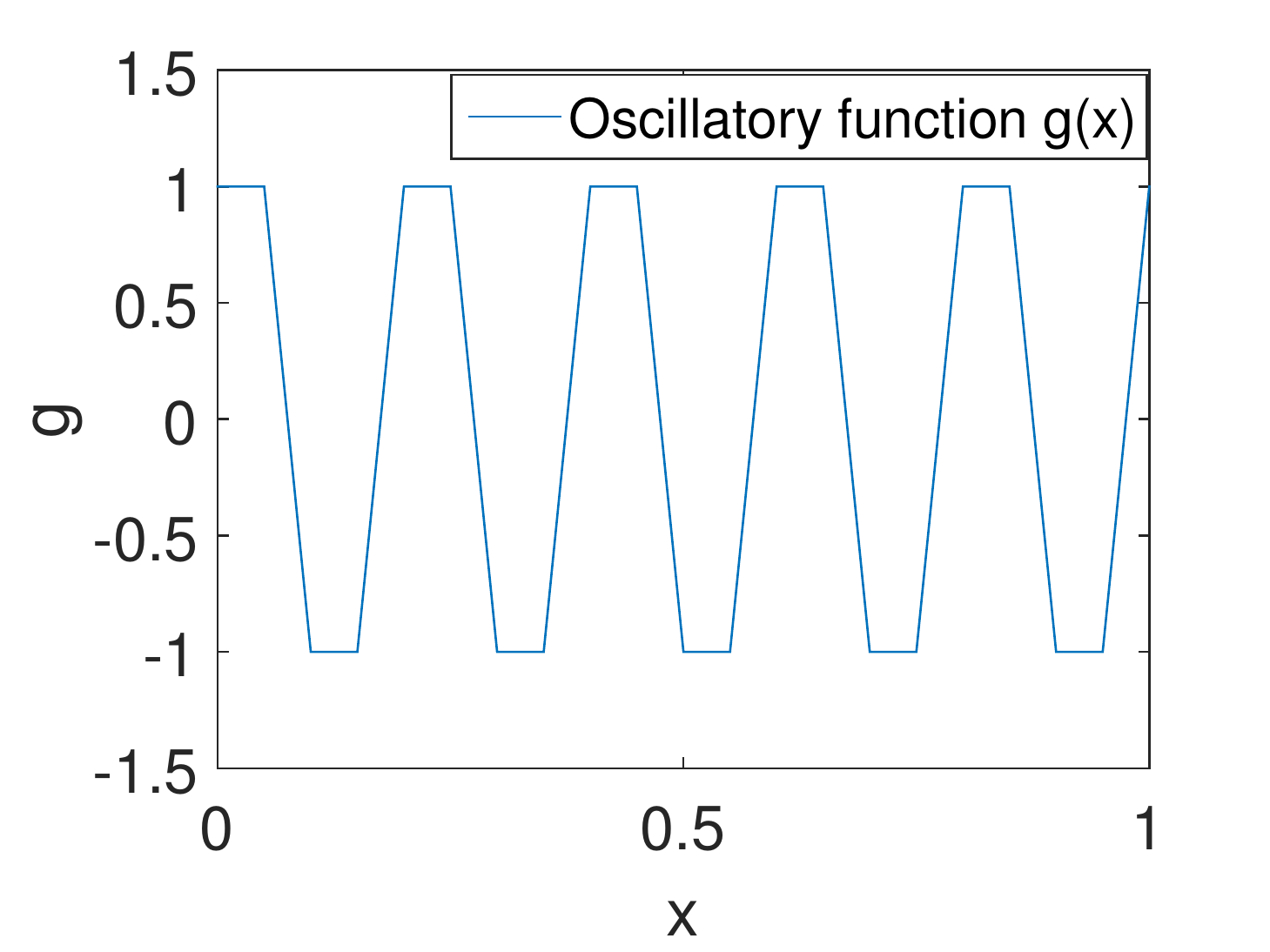}
   \end{center}
\caption{Oscillatory function applied to $x$-direction or $y$-direction}\label{fig:squarewave}
  \end{figure}

In the third experiment, we reuse the offline computation above to explore the anisotropic property of this random media. We first set 
\begin{equation} \label{eqn:xsquarewave}
\begin{split}
u(x,y,\omega)\mid_{x=0.1} &= u(x,y,\omega)\mid_{x=0.9} = g(y) \\
\vct{n}\cdot \nabla u &= 0 \quad \text{on $\partial \CalD$},
\end{split}
\end{equation}
where $g(y)$ is the oscillatory function shown in Figure \ref{fig:squarewave}. With zero-source term, we get one solution $u_H(x,y,\omega)$.
  
We then set the same zero-Neumann boundary condition and zero-source term but specify the oscillatory function in the $y$-direction
\begin{equation} \label{eqn:ysquarewave}
u(x,y,\omega)\mid_{y=0.1} = u(x,y,\omega)\mid_{y=0.9} = g(x),
\end{equation}
and get another solution $u_H(x,y,\omega)$

\begin{figure}[htp] \label{fig:anisotropic}
   \begin{center}
     {\includegraphics[width=0.35\textwidth]{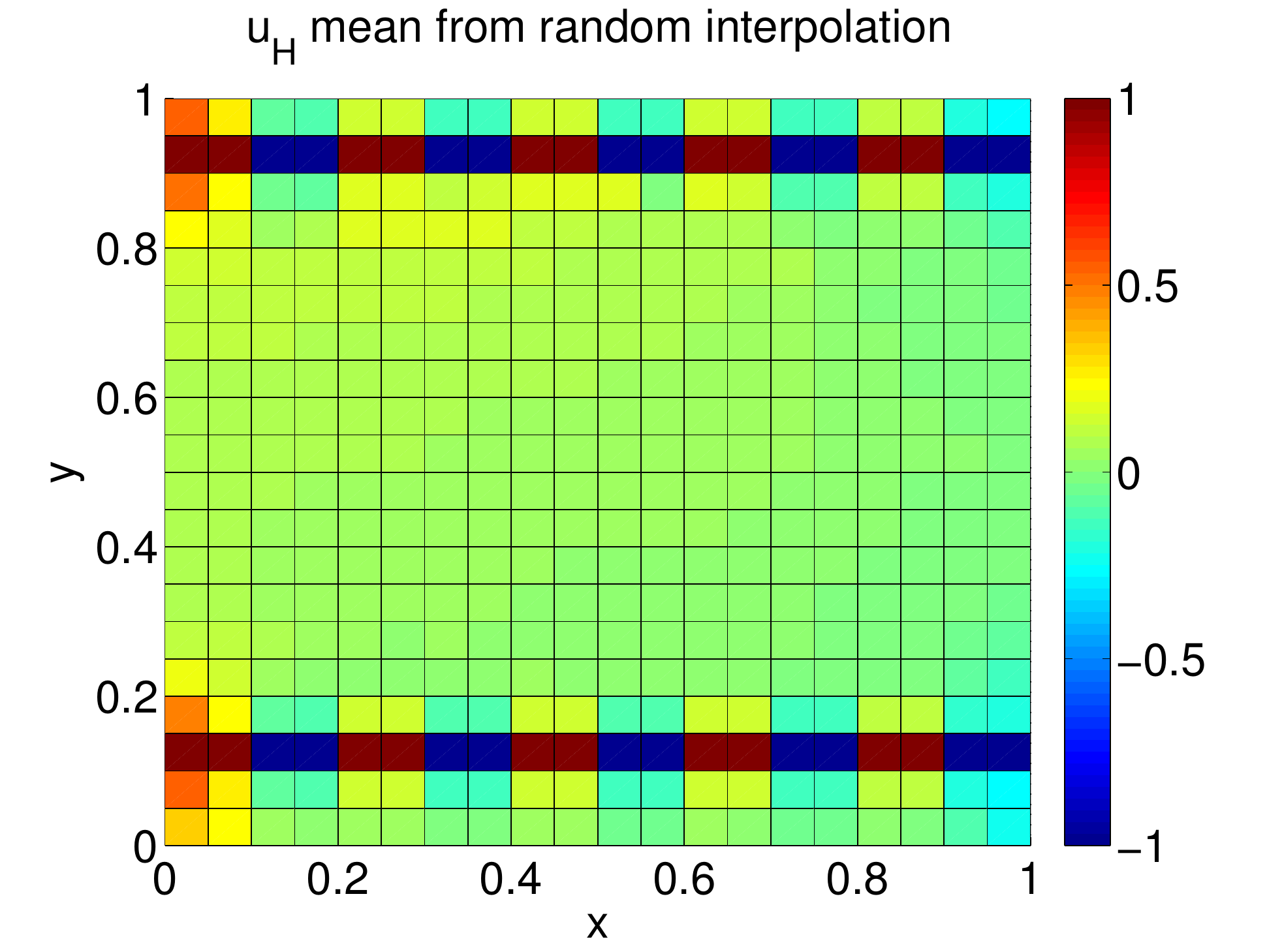}
      \includegraphics[width=0.35\textwidth]{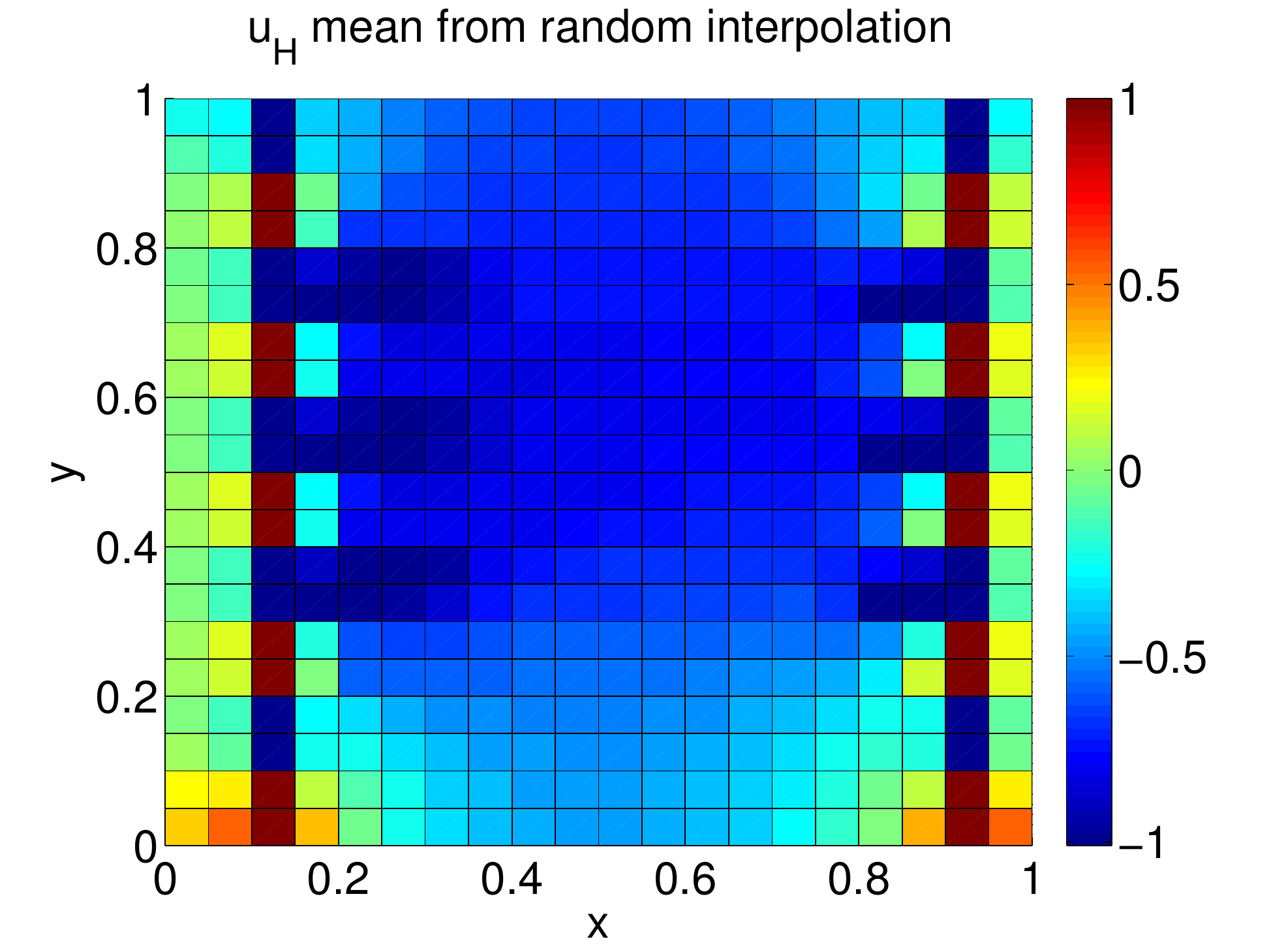}}
   \end{center}
 \caption{Left: oscillatory condition applied in $x$-direction; right: oscillatory condition applied in $y$-direction. The medium behaves as a homogeneous medium in the first setting, but shows high conductivity in the second setting.}\label{fig:HC_statistics}
\end{figure}

We compute $10^5$ samples for each example and compare their means in Figure~\ref{fig:anisotropic}. Because the results of the two approximations are visually the same, we only show the results from the random interpolation method. Since high conductivity channels are presented along the $x$-direction, the medium behaves as a homogeneous medium in the first setting, but shows high conductivity in the second setting. Note that our local upscaled quantities $\vct{S}^m$ are independent of the boundary conditions and the source functions, and thus we can reuse them for different settings.
 
\subsection{A 2d example with short correlation length}\label{subsec:guassiandecay}
\begin{figure}[htp]
\begin{center}
	\includegraphics[width = 0.3\textheight]{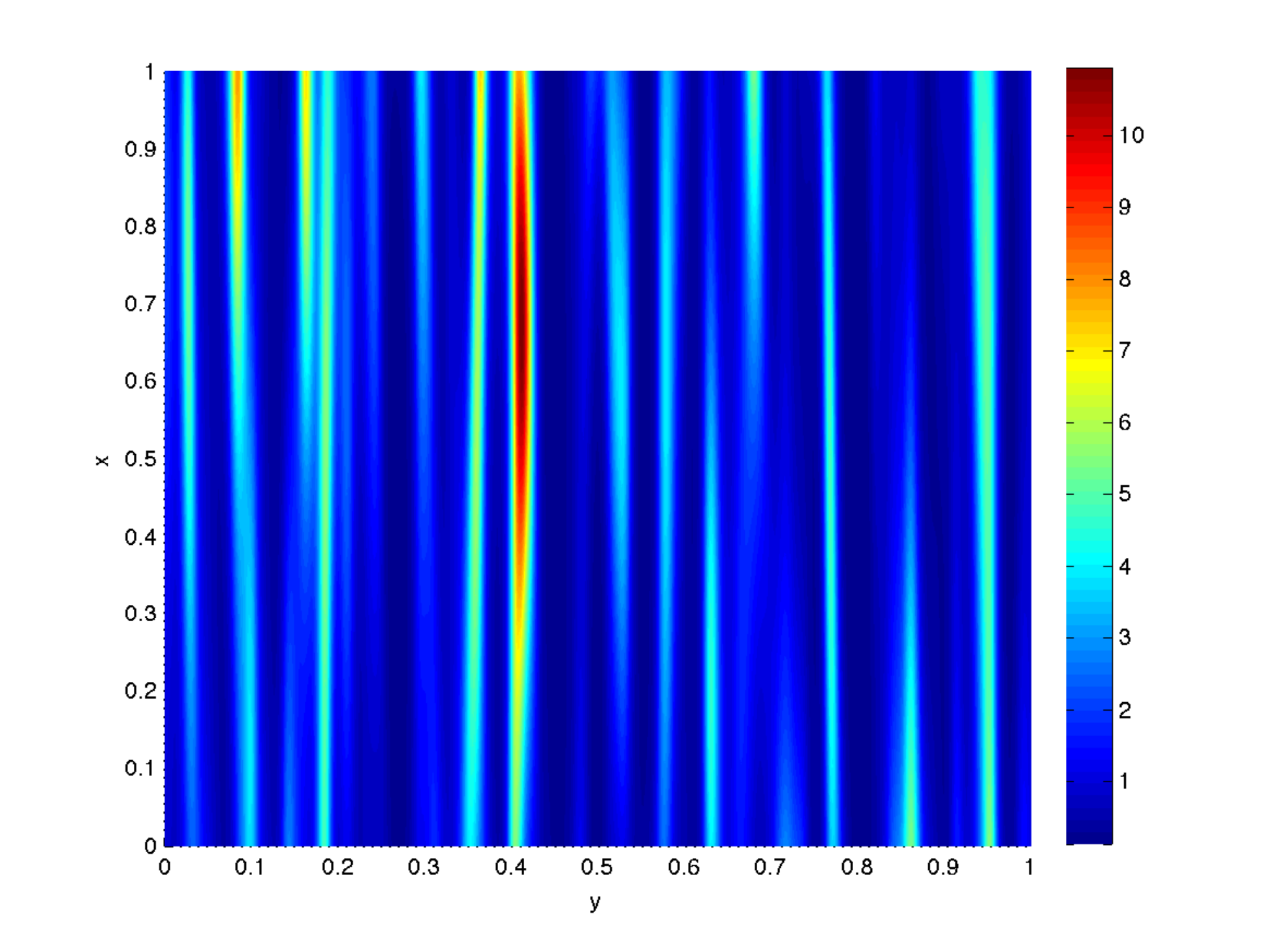}
	\includegraphics[width = 0.3\textheight]{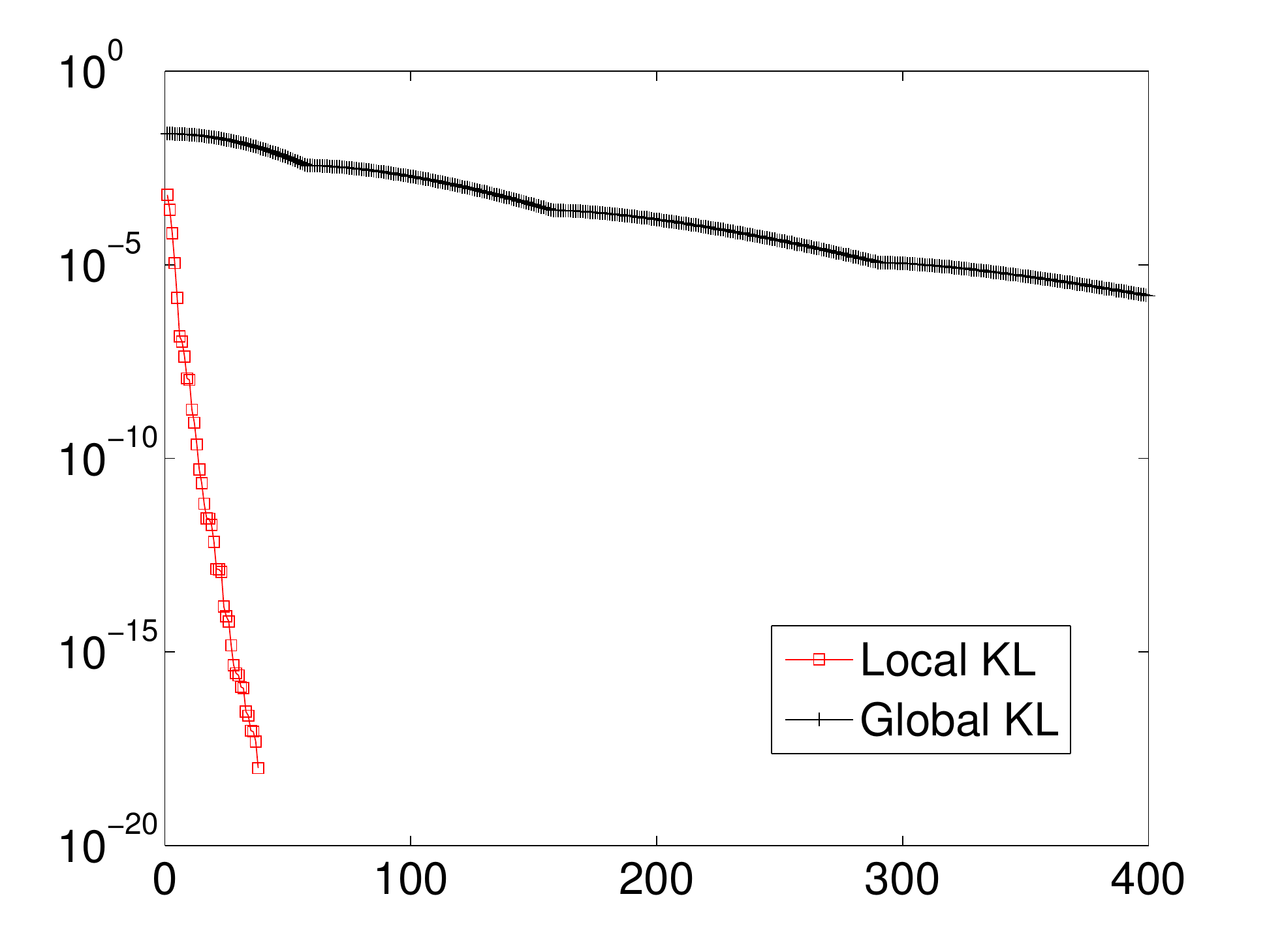}
\end{center}
\caption{Left: one sample of the anisotropic random media; right: the eigenvalues of the global and local KL expansion for covariance function $\Cov_\beta$.}
\label{fig:mediasampleEx2}
\end{figure}

In this example, we consider the two-dimensional elliptic problem~\eqref{eqn} in the physical domain $D = [0,1]^2$ with the source $b(x_1,x_2) = 2 + x_1 x_2$ and the zero Dirichlet boundary condition. The random medium $\kappa(\vct{x},\omega)$ is given as
\begin{equation}\label{eqn:mediaEx2}
	\kappa(\vct{x},\omega) = 0.1 + \exp(\beta(\vct{x},\omega))
\end{equation}
where $\beta(\vct{x},\omega)$ is a Gaussian random field with zero mean and a Gaussian covariance function
\begin{equation}\label{eqn:betaCov}
	\Cov_{\beta}(\vct{x},\vct{y}) = \exp(-\frac{|x_1 - y_1|^2}{l_1^2} - \frac{|x_2 - y_2|^2}{l_2^2}), \qquad l_1 = 1,\quad l_2 = 1/64.
\end{equation}
Here, we have different correlation lengths in $x_1$ and $x_2$ directions to model the anisotropic media. A sample of the random media is shown in Figure~\ref{fig:mediasampleEx2} (left). We can clearly see the small scales in the $x_2$ direction due to the small correlation length. 

We apply our StoMsFEM to solve this elliptic problem on a coarse mesh $\CalT_H$ with mesh size $H_x = H_y = 2^{-6}$, which does not resolve the fine scales. For the local upscaling method, we use the MsFEM with oversampling ratio $\eta = 2$ and with the oscillatory boundary condition under the Petrov-Galerkin formulation. We solve the local cell problems~\eqref{eqn:MsFEB} on a fine mesh $\CalT_h$ with mesh size $h_x = h_y = 2^{-11}$. For the parametrization on every local patch, i.e. every oversampling domain, we use the local KL expansion. Figure~\ref{fig:mediasampleEx2} (right) shows the eigenvalues for the global and local KL expansions of the Gaussian covariance function $\Cov_\beta(\vct{x}, \vct{y})$. It is obvious that the local KL expansion exhibits a much faster eigenvalue decay than the global KL expansion. In fact, to keep about 99\% of the total spectrum, the global KL requires 168 terms, whereas the local KL expansion (on the oversampling domain) requires only 4 terms. The stochastic dimensionality of the global KL expansion is 168, which is too high for most gPC based stochastic methods. On the other hand, the local stochastic dimensionality is only 4, and the random interpolation method works well in this non-affine parametrization setting.

\begin{figure}[htp]
\begin{center}
	\includegraphics[width=0.4\textwidth]{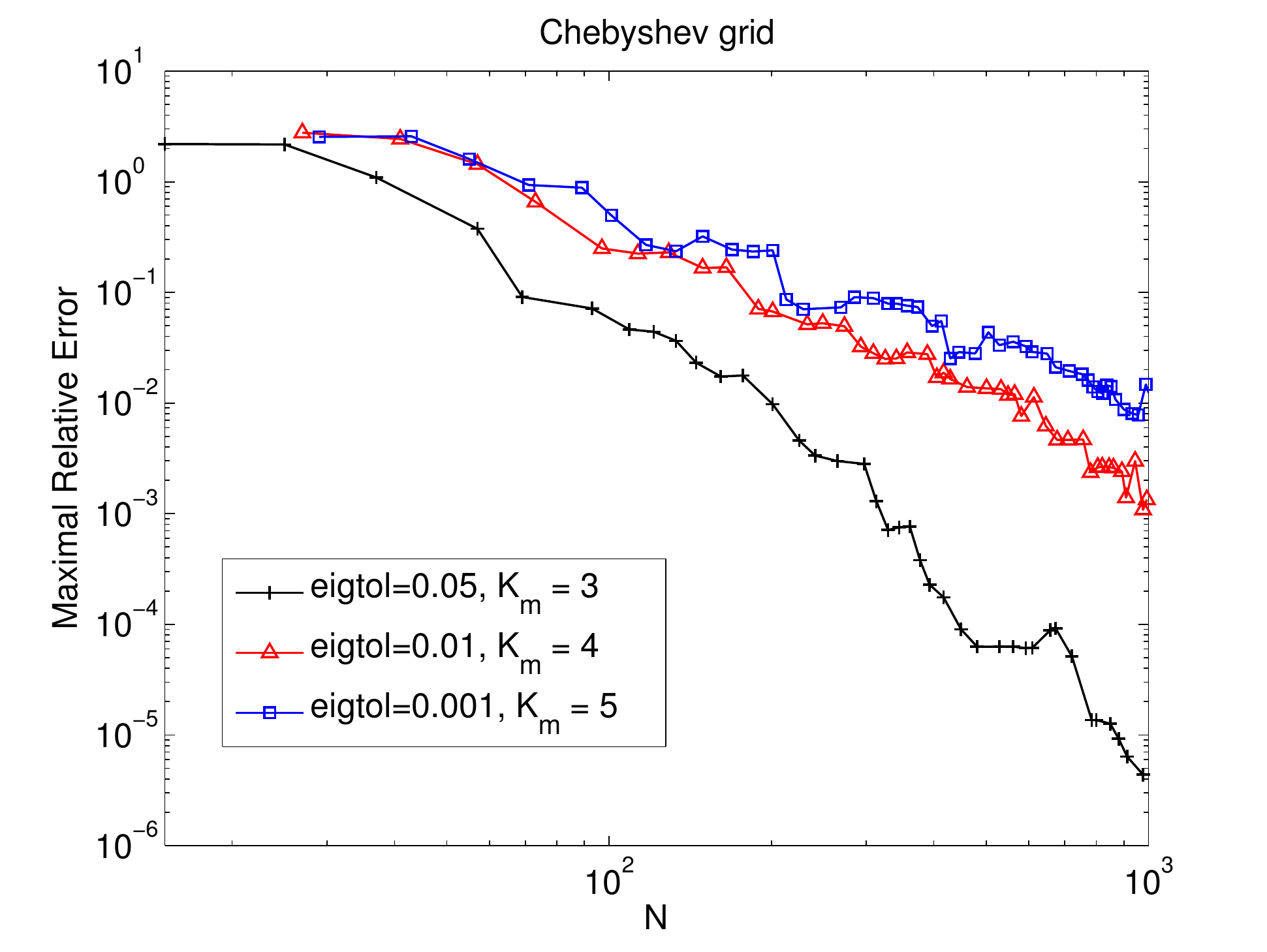}
	\includegraphics[width=0.4\textwidth]{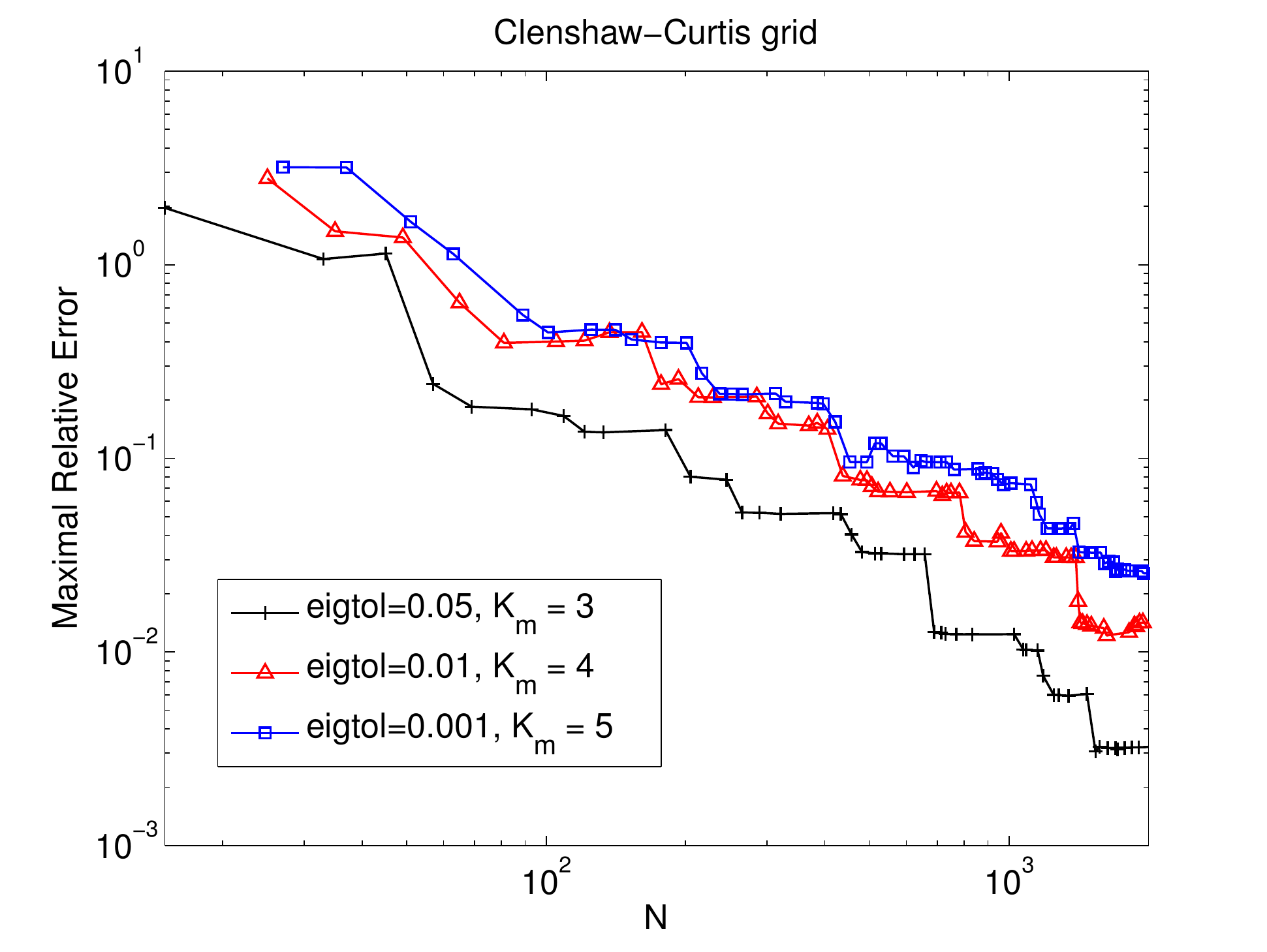}
\end{center}
\caption{Interpolation error of $\vct{S}_{1,1}^{m}$ on sparse grids versus number of interpolation nodes}
\label{fig:intererrorEx2}
\end{figure}

\begin{figure}[htp]
\begin{center}
	\includegraphics[width = 0.3\textheight]{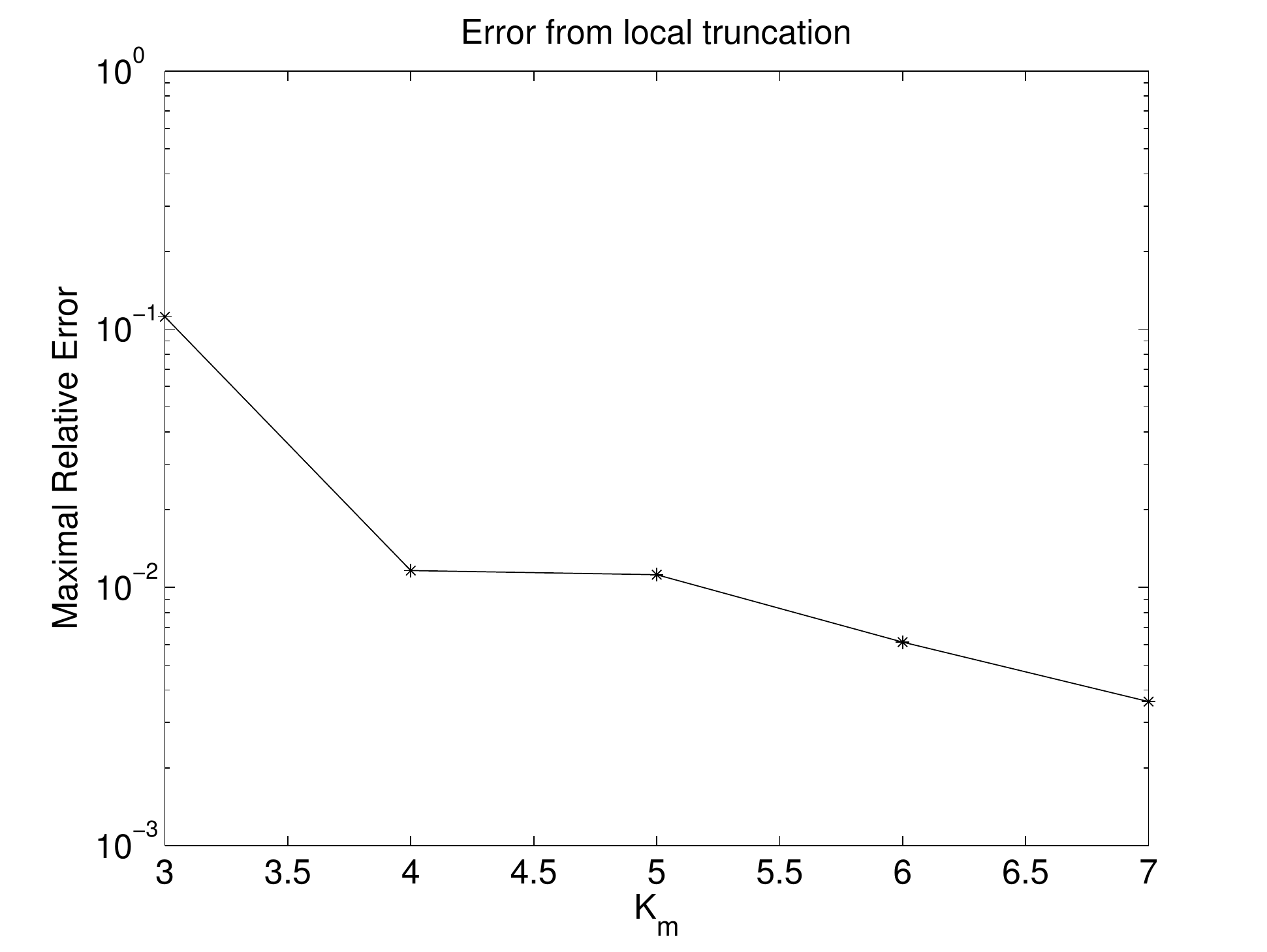}
\end{center}
\caption{Error of $\vct{S}_{1,1}^{m}$ versus the number of terms retained in the local KL expansion}
\label{fig:truncerrorEx2}
\end{figure}

In the offline stage, we construct the interpolants $\hat{\vct{S}}_m$ for the local upscaled stiffness matrices. Notice that the local parameters have a standard normal distribution, whose support is $(-\infty, +\infty)$. To make the interpolation accurate point-wisely, we construct interpolants when all the parameters lie in $[-3,3]$. We utilize the polynomial interpolation on the Chebyshev Gauss-Lobatto sparse grid~\cite{spalg, spdoc}. The interpolation error is estimated by the largest error among $10^4$ randomly drawn points. The relative interpolation error of $\vct{S}_{1,1}^{m}$ versus the number of interpolation nodes is shown in Figure~\ref{fig:intererrorEx2}. We show the results for 3 local KL expansions, which keep 95\%, 99\% and 99.9\% of the total spectrum and whose local dimensions are 3, 4, and 5 respectively. We also plot the interpolation error with the piecewise linear interpolation on the Clenshaw-Curtis sparse grid. We can clearly observe that with the same number of interpolation nodes, the high-order polynomial interpolation is more accurate than the piecewise linear interpolation. 
In Figure~\ref{fig:intererrorEx2}, we only consider the interpolation error, i.e. $\hat{\vct{S}}^m(\vct{\xi}_m) - \vct{S}^{m}(\vct{\xi}_m)$. Given a sample of the medium, denoted as $\kappa(x,\omega)$, we project it onto the local KL modes, obtain the local parameters $\vct{\xi}_m$ and truncate the small terms in the local KL expansion. The upscaled stiffness matrix $\vct{S}^{m}(\vct{\xi}_m)$ is defined based on this truncated local KL expansion. Due to this truncation, we introduce another source of error $\vct{S}^{m}(\vct{\xi}_m) - \vct{S}^{m}(\omega)$, where $\vct{S}^{m}(\omega)$ is the exact upscaled stiffness matrix based on the sample $\kappa(x,\omega)$. This truncation error is plotted in Figure~\ref{fig:truncerrorEx2} with respect to the number of terms we keep in the local KL expansion. The error is estimated by the largest error among $10^4$ randomly drawn samples. We notice that the errors decay as more terms in the local KL expansion are retained. We also note that the errors decay slower than those in Figure~\ref{fig:intererrorEx2}, which implies that the predominant contribution in the overall error $\hat{\vct{S}}^m(\vct{\xi}_m) - \vct{S}^{m}(\omega)$ is from the truncation of the local KL expansion. A theoretical result to estimate the parametrization error can be found in Corollary 2.1 in~\cite{babuska_galerkin_2004}.

We emphasize that the errors introduced by truncation of the local KL expansions will typically dominate the errors induced by random interpolation or reduced basis method when computing $\vct{S}^{m}(\omega)$. As we discussed in Section~\ref{sec:RandomInter_theory}, $\vct{S}^{m}(\vct{\xi}_m)$ is smooth and thus the random interpolation method will converge very fast. In addition, since the local dimensions are of order 1, we are able to compute reasonably high order interpolants. The fast convergence and locally low dimensionality mean that we can easily drive the interpolation error $\hat{\vct{S}}_m(\vct{\xi}_m) - \vct{S}^{m}(\vct{\xi}_m)$ below the error induced by the truncation of the local KL expansions, i.e. $\vct{S}^{m}(\vct{\xi}_m) - \vct{S}^{m}(\omega)$. Consequently, the error arising from the local KL expansion provides the leading contribution to the total error of the StoMsFEM method. Errors introduced by the parametrization, e.g., either the standard KL expansion or the local parametrization methods in Section~\ref{sec:parametrization}, should be considered as a modeling error, since they are not directly related to the StoMsFEM algorithm. Moreover, the local parametrization methods presented in this paper has smaller parametrization error compared with the popular parametrization by the global KL expansion, because the local parametrization methods allow one to capture a greater percentage of the uncertainties due to the fast eigenvalue decay in the local KL expansion.

\begin{figure}[htp]
\centering
     \includegraphics[width=0.45\textwidth]{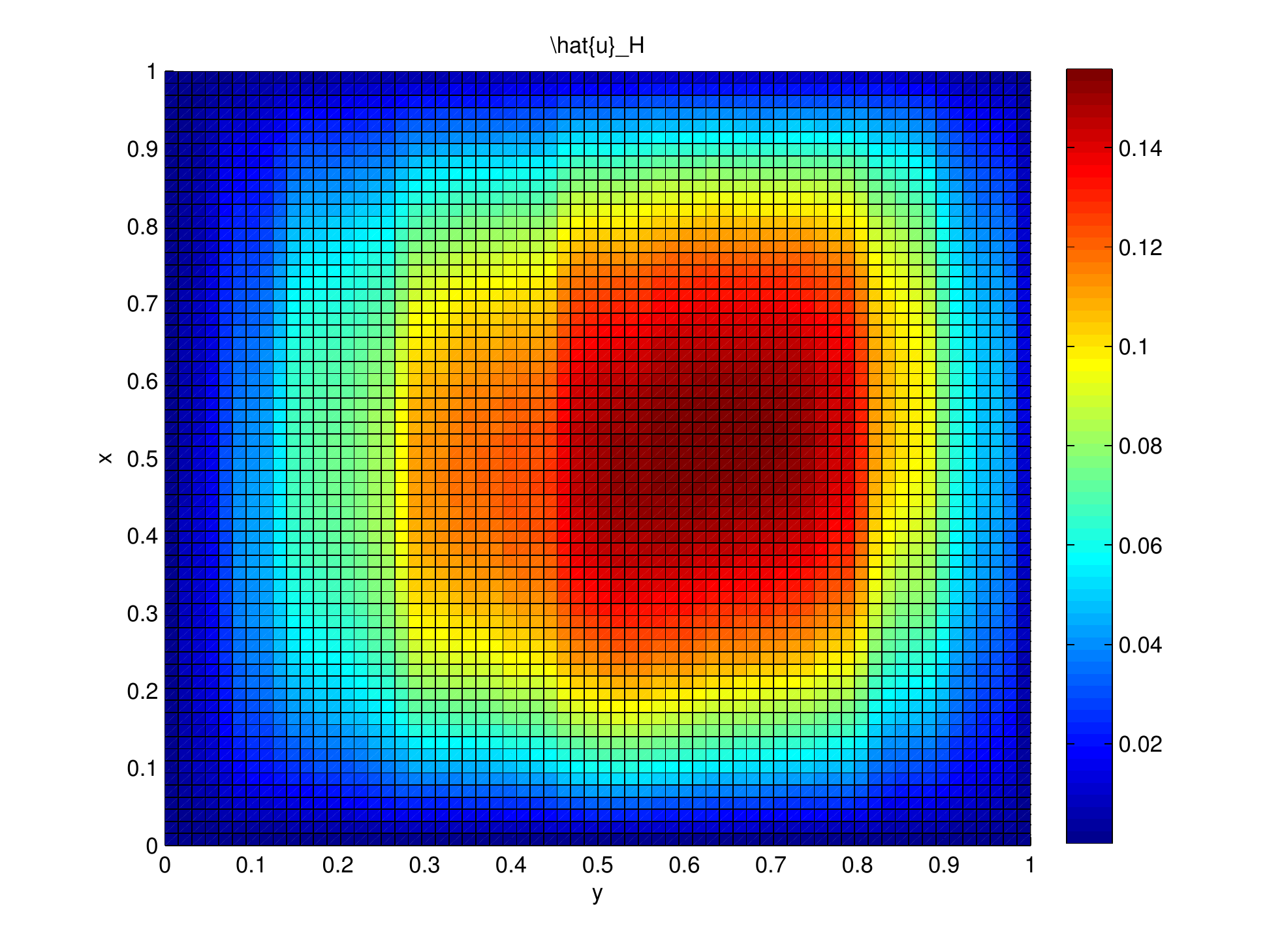}
\caption{The sample solution corresponding to the medium sample in Figure~\ref{fig:mediasampleEx2}. There are several layers in the $y$-direction (the horizontal direction). The boundary of these layers are exactly the low permeability strips in the medium sample.}\label{fig:solutionsampleEx2}
\end{figure}

In the online stage, we retain 4 terms in the local KL expansion, which keeps 99\% of the total spectrum. The interpolants for $\vct{S}^{m}(\vct{\xi}_m)$ have about 1000 interpolation nodes and the maximal relative error is below 1\%, which is smaller than the error induced by the truncated local KL expansion. For every sample, we first generate the media sample $\kappa(x,\omega)$ by the standard spectral method. Then we project it onto the local KL modes and get the local parameters $\vct{\xi}_m$. If all the local parameters lie in $[-3,3]$, we evaluate the interpolant $\hat{\vct{S}}_m$. Otherwise, we we directly solve the multiscale basis functions on the fine grid and assemble the local stiffness matrix directly from~\eqref{eqn:MsFEB} and \eqref{eqn:stiffness_local}. In our case, the probability to do interpolation for a local stiffness matrix is 0.9892, and with a very small probability 0.0108 the multiscale basis functions are required to solve on the fine grid. Finally, we solve the upscaled system~\eqref{eqn:upscaleSUF} to get the coarse grid solution $\hat{u}_H$. In Figure~\ref{fig:solutionsampleEx2}, we show the solution sample corresponding to the media sample in Figure~\ref{fig:mediasampleEx2}. We find that the error $\hat{u}_H - u_H$ induced by the interpolation is of order $10^{-5}$, while the error $u_H(x,\vct{\xi}) - u_H(x,\omega)$ induced by the truncated local KL expansion is of order $10^{-4}$. Therefore, we again confirm that the local truncation error is the main contribution of the overall error $\hat{u}_H(x,\vct{\xi}) - u_H(x,\omega)$. 

\begin{figure}[htp]
\centering
     \includegraphics[width=0.45\textwidth]{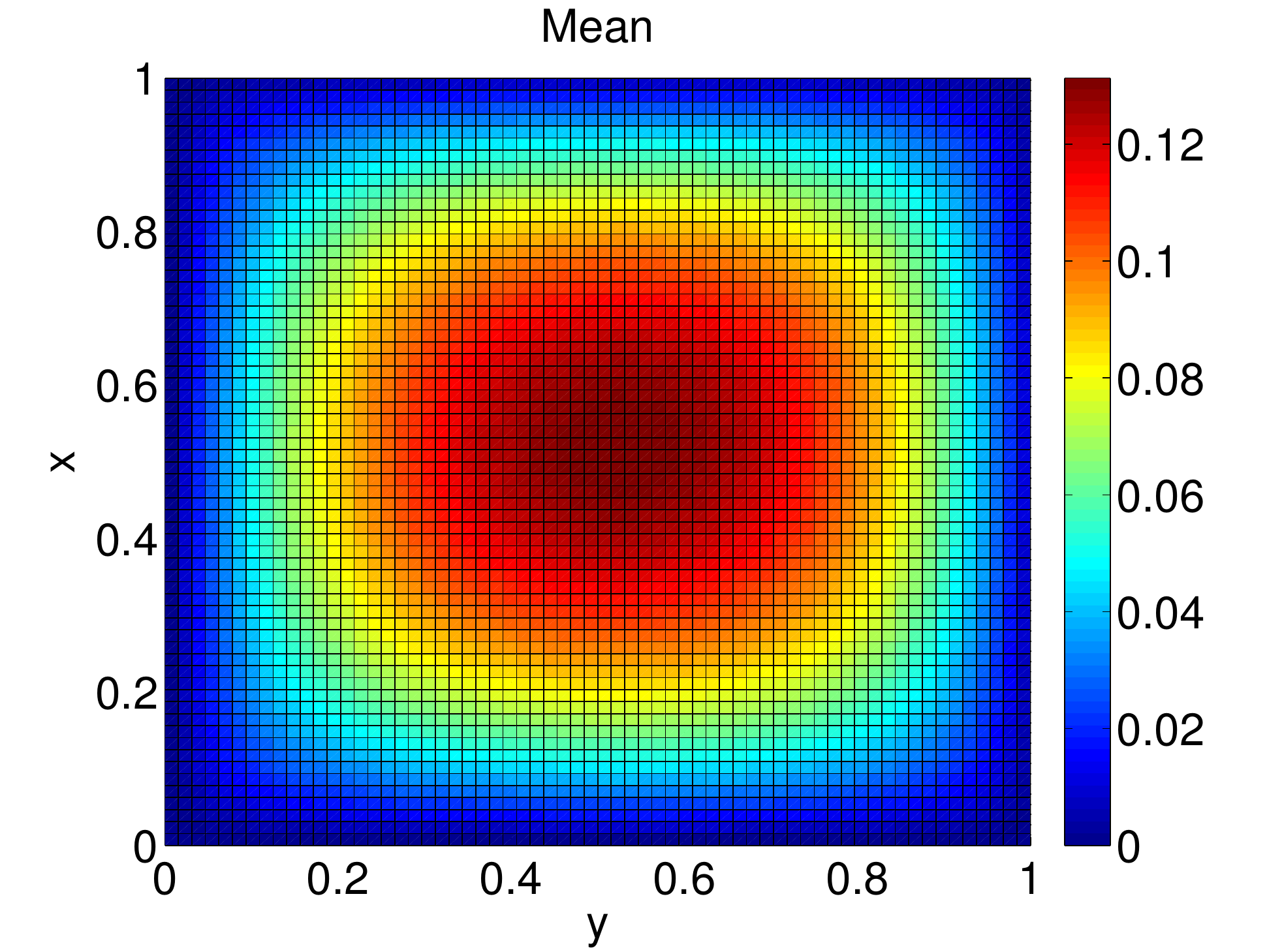}
     \includegraphics[width=0.45\textwidth]{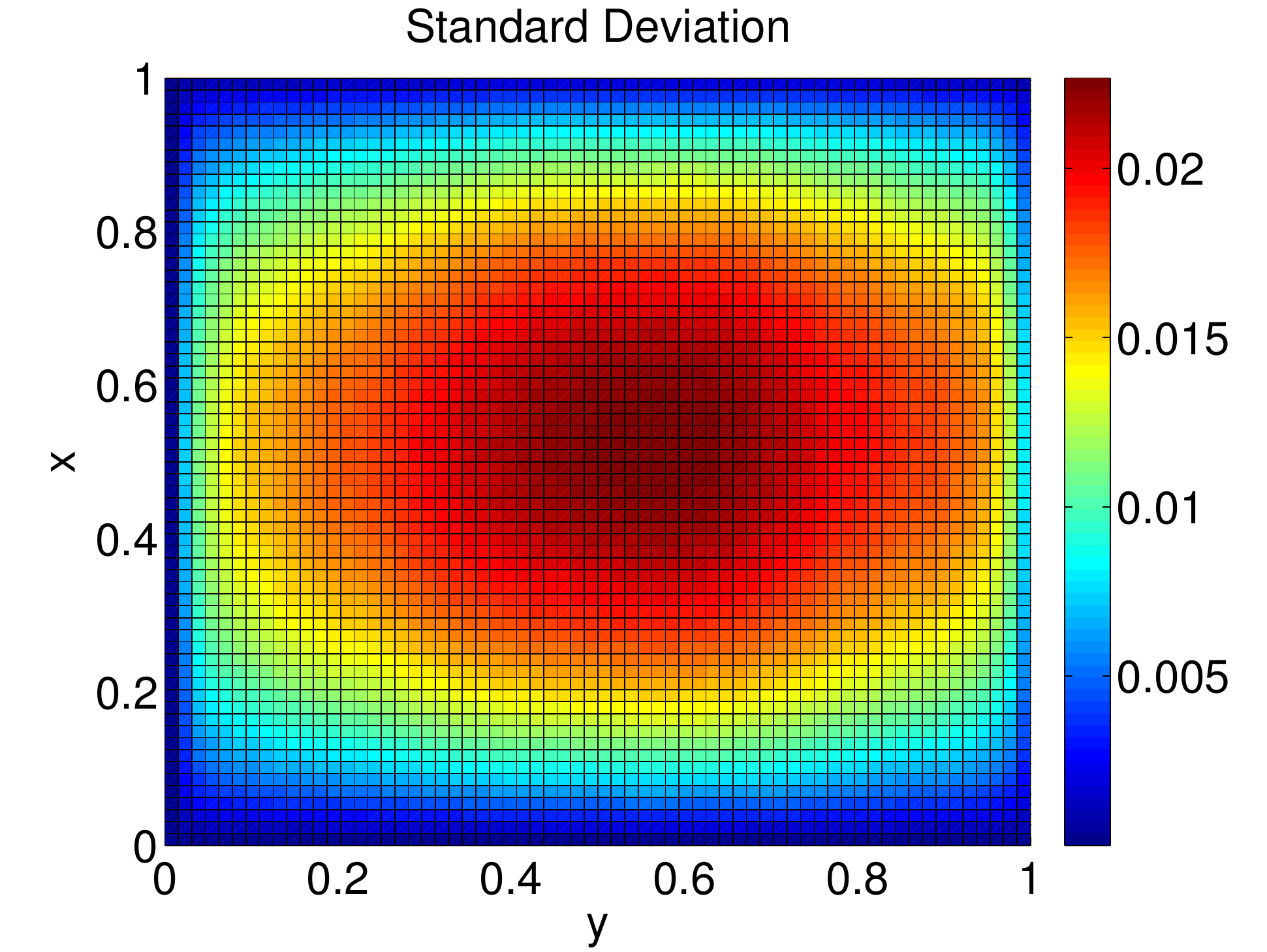}
\caption{The statistics of the numerical solutions using $10^3$ samples. Note that the layered structures in sample solutions have been ``homogenized'' when we average over samples.}\label{fig:statisticsEx2}
\end{figure}

We summarize the computational cost in Table~\ref{table:costEx2}. The offline cost for the StoMsFEM is extremely small in this example because the random field $\kappa(x,\omega)$ is translational invariant and we can construct interpolants for $\vct{S}^{m}(\vct{\xi}_m)$ on only one local domain. In the online stage, the ratio of the computational cost between the naive application of the MsFEM and the StoMsFEM is about 34. Finally, we use the Monte Carlo method as the global stochastic method to estimate the statistical properties of $\hat{u}_H(x,\omega)$. In Figure~\ref{fig:statisticsEx2} we show its mean and standard deviation estimated from $10^3$ samples. To compute these $10^3$ samples, it takes 17431 seconds for te StoMsFEM. If we directly compute these samples by naively applying the MsFEM, it would take $6.026\times 10^5$ seconds. We remark that to balance the spatial discretization error and stochastic sampling error, we need $N_{\text{on}} \approx \Or(H^{-4}) = \Or(10^7)$ samples and $10^3$ is far less than enough. This fact shows the necessity of StoMsFEM because the computational saving from StoMsFEM grows nearly with the number of samples we solve in the online stage.

\begin{table}[ht]
\caption{Computational cost for one sample solution(unit: s)}\label{table:costEx2}
\centering
\begin{tabular}{c|c|c}
\hline
naive MsFEM per sample      &  StoMsFEM (offline) & StoMsFEM (online per sample) \\
\hline
602.6104 		& 	389.3256	&   17.8173 (1/R = 34)	\\
\hline
\end{tabular}
\end{table}

\begin{remark}
Similar to the two-level MC in Remark~\ref{rem:MLMC}, the error induced by the local KL expansion can be corrected with a little more computation. We can write the coarse grid solution as
\begin{equation*}
	u_H(x, \omega) = \hat{u}_H(x,\vct{\xi}) + u_H(x, \omega) - \hat{u}_H(x,\vct{\xi}).
\end{equation*} 
From Figure~\ref{fig:solutionsampleEx2}, we can see that $\hat{u}_H$ serves as a good variance reduction for the true solution on the fine grid $u_H(x,\omega)$. With a few more samples of $u_H(x, \omega)$, we can correct the errors introduced both by the random interpolation and by the local KL expansion.
\end{remark}
  
\section{Conclusions and future work}

We proposed a stochastic multiscale finite element method (StoMsFEM) to solve random elliptic partial differential equations with a high stochastic dimension. An essential difficulty in solving this type of elliptic random PDEs is that we need to solve a huge number of sample solutions to get an acceptable statistical estimation and that the computational cost for every sample solution is already quite expensive since we need to resolve the small scale feature of the solution. The StoMsFEM saves computational cost for every sample by simultaneously upscaling the stochastic solutions in the physical space for all random samples and exploring the low stochastic dimensions of the stochastic solution within each local patch.

Moreover, we proposed two effective methods to achieve this simultaneous local upscaling. The first method is the random interpolation method that explores the high regularity of the local upscaled quantities with respect to the random variables. The second method is the reduced-order method that explores the low rank property of the multiscale basis functions within each coarse grid element. For every sample solution, our complexity analysis shows that the cost ratio between the StoMsFEM and the standard FEM on find grid is $R$, where $R = \Or(N_c (h/H)^{\gamma d})$ for the random interpolation method and $R = \Or(K_m^3 (h/H)^{\gamma d})$ for the reduced basis method. In practice, the saving is even more significant due to highly optimized fast numerical interpolation methods. In our high contrast example, we observed a factor of 2000 speed-up by the random interpolation method. 

We also analyzed different kinds of errors contributed to the final statistical estimation error. We showed that the error introduced by the interpolation or the reduced basis method is negligible, and thus we can treat our approximating solution $\hat{u}_H(x, \vct{\xi})$ as the solution $u_H(x,\vct{\xi})$ directly computed from MsFEM. We also showed that StoMsFEM optimally balances the spatial discretization error from MsFEM and the stochastic sampling error from the global stochastic methods. In comparison, the standard FEM on the fine grid wastes a lot of computational resources on resolving the small scale of the solution in order to reduce the spatial discretization error, while the total error is actually dominated by the stochastic sampling error. Therefore, to achieve the same level of estimation error, the StoMsFEM indeed offers a factor of $R$ computational saving compared with the standard FEM on a fine grid. 

In our last numerical example, we discussed the modeling error introduced by the local KL expansion when we parametrize the random medium. We showed that the errors introduced by truncating the local KL expansions will typically dominate the errors induced by the random interpolation or the reduced basis method. In our future work, we plan to combine the local parametrization step and the parametric local upscaling step together, and to optimally balance the modeling error from the local parametrization and that from the parametric local upscaling. We also briefly discussed the two-level Monte Carlo approach to achieve an $\Or(h)$ statistical estimation error. This topic will be further explored in our future work.

\vspace{0.2in}
\noindent
{\bf Acknowledgements.}
This research was in part supported by NSF Grants No. DMS-1318377 and DMS-1613861.

\bibliographystyle{plain}
\bibliography{random_elliptic}
\end{document}